\magnification=\magstep1   
\input amstex
\UseAMSsymbols
\input pictex 
\vsize=23truecm
\NoBlackBoxes
 
\parindent=18pt
  
   \font\rmk=cmr8      

\font\gross=cmbx10 scaled\magstep1 

\def\op{\text{\rm op}}
\def\grade{\operatorname{grade}}

\def\mod{\operatorname{mod}}

\def\Hom{\operatorname{Hom}}

\def\End{\operatorname{End}}
\def\Ext{\operatorname{Ext}}

\def\rad{\operatorname{rad}}
\def\add{\operatorname{add}}
\def\Ker{\operatorname{Ker}}
\def\Cok{\operatorname{Cok}}
\def\soc{\operatorname{soc}}
\def\Tr{\operatorname{Tr}}
\def\Im{\operatorname{Im}}

\def\top{\operatorname{top}}

\def\del{\operatorname{del}}
\def\des{\operatorname{des}}
\def\depth{\operatorname{depth}}

\def\pd{\operatorname{pd}}
\def\id{\operatorname{id}}

\def\fin{\operatorname{fin}}

\def\finpro{\operatorname{fin-pro}}
\def\Finpro{\operatorname{Fin-pro}}

\def\fininj{\operatorname{fin-inj}}
\def\Fininj{\operatorname{Fin-inj}}

\def\pro{\operatorname{pro}}
  \def\ss{\ssize }
\def\arr#1#2{\arrow <1.5mm> [0.25,0.75] from #1 to #2}

\def\s{\hfill \square} 
\def\snake{
  \beginpicture
    \setcoordinatesystem units <.1cm,.1cm>
\multiput{} at 0 0  3 0  /
  \setquadratic
\plot 0 0  0.5 0.1  1 0.4 
           1.5 .5  2 0
           3 -.5  4 0 
           4.5 .5  5 0.4 
           5.5 0.1  6 0 /
\endpicture}

\vglue1.5truecm
\centerline{\gross The finitistic dimension of a Nakayama algebra.}
                     \bigskip
\centerline{Claus Michael Ringel}
                \bigskip\medskip

\noindent {\narrower Abstract:  \rmk Let $\ss A$ be an artin algebra,
G\'elinas has introduced an interesting upper bound for the
finitistic dimension $\ss\finpro A$ of $\ss A$,
namely the delooping level $\ss\del A$.  We assert that $\ss \finpro A\ =\ \del A$ for
any Nakayama algebra $\ss A$.  This yields also a new proof
that the finitistic dimension of $\ss A$ and its opposite algebra are equal,
as shown recently by Sen. If $\ss S$ is a simple module, let $\ss e(S)$ be
the minimum of the projective dimension of $\ss S$ and of its injective envelope
(one of these numbers has to be finite); and
$\ss e^*(S)$ 
the minimum of the injective dimension of $\ss S$ and of its projective cover.
Then the finitistic dimension of $\ss A$ is the maximum of the numbers $\ss e(S)$,
as well as the maximum of the numbers $\ss e^*(S)$.

Using suitable syzygy modules, we construct a permutation $\ss h$
of the simple modules $\ss S$ such that $\ss e^*(h(S)) = e(S).$ In particular, this shows 
for $\ss z\in \Bbb N$ that the 
number of simple modules $\ss S$ with $\ss e(S) = z$ is equal to
the number of simple modules $\ss S'$ with $\ss e^*(S') = z$.
	\medskip
\noindent 
\rmk Key words. Nakayama algebra. Finitistic dimension. Delooping level. 
Desuspending level. The grade of a module.
The homological permutation. Catalan combinatorics.
	\medskip
\noindent 
2020 Mathematics Subject classification. 
  Primary 16G10,   Secondary  16G20, 16G70, 16E05, 16E10, 16E65. 
\par}
	\medskip\bigskip
{\bf 1. Introduction.}
	\medskip 
Let $A$ be an artin algebra. The modules $M$ to be considered are usually left 
$A$-modules. Often (but not always) we will assume in addition that $M$ is finitely generated.
We denote by $A^{\op}$ the opposite algebra of $A$ (and stress
that the right $A$-modules are just the (left) $A^{\op}$-modules).

The projective dimension of a module $M$ will be denoted by $\pd M$, its injective
dimension by $\id M$. We write $\finpro A$ for the supremum of the projective
dimension of finitely generated left modules of finite projective dimension, and
similarly, $\fininj A$ for the supremum of the injective
dimension of finitely generated left modules of finite injective dimension.
Note that $\finpro A$ is called the (small) {\it finitistic dimension} of $A$ 
and $\fininj A$ is just the finitistic dimension of the opposite algebra 
$A^{\op}.$ 
The finitistic dimension was considered first by Serre [S] in 1955. In 1960, its relevance was
stressed by Bass [B]. 
It lies at the center of the classical homological conjectures.

We will focus the attention to Nakayama algebras:
An algebra $A$ is called {\it Nakayama} provided all indecomposable modules 
are serial (that means: they have a unique composition series).
The module category of a Nakayama algebra is well understood and many questions
which are difficult to deal with in general, can easily be answered for Nakayama
algebras. But, actually, only very recently, it was shown by Sen [S2]
that $\finpro A = \fininj A$ for a Nakayama algebra $A$. This seems to be 
an important observation. We will provide
here a new proof as well as some further information on $\finpro A$. 
	\medskip
Being interested in the finitistic dimension of an algebra $A$, G\'elinas [Ge] has
introduced a new invariant,  
the delooping level $\del A$. If $M$ is a module,
let $\Omega M$ be its (first) {\it syzygy module} (it is
the kernel of a projective cover $PM \to M$, and sometimes also called the ``loop'' module
of $M$) and  let $\Sigma M$ be its (first) {\it suspension
module} (the cokernel of an injective envelope $M \to IM$). 
The {\it delooping level} $\del S$ of
a simple module $S$ is the smallest number $d\ge 0$ such that $\Omega^d S$
is a direct summand of $P\oplus \Omega^{d+1} M$, where $M,\ P$ are finitely generated modules with $P$ projective
(and $\del S = \infty$, if such a number $d$ does not exist). 
Similarly, one may define the {\it desuspending level}
$\des S$ of $S$ as  the smallest number $d\ge 0$ such that 
$\Sigma^d S$ is a direct summand of $I\oplus \Sigma^{d+1} M$
for some finitely generated modules $M,\ I$ with $I$ injective 
(and $\des S = \infty$, if such a $d$ does not exist). 
By definition, $\del A$ is the maximum of $\del S$, where $S$ runs through the simple
modules, and $\des A$ is the corresponding maximum of the numbers $\des S$.
Of course, $\des A = \del A^{\op}$. 

If $A$ is a Nakayama algebra, we will show that $\finpro A = \del A$, thus, altogether 
$$
 \finpro A = \fininj A = \del A = \des A.
$$
	
The proof provides an explicit formula for the finitistic dimension of a Nakayama algebra $A$, namely
$$
 \fin\pro A = \max\nolimits_S \min\{\pd S, \pd IS\},
$$
where $S$ runs through the simple modules (note that we will show that 
at least one of $\pd S,\ \pd IS$ is finite). It may happen that both $S$ and $IS$ 
have finite, but different, projective dimension, then the formula shows that we may delete 
$S$ on the right when looking for the maximum! Also, we may delete on the right 
all the torsionless simple modules (namely, if $S$ is torsionless and simple, then $\pd IS = 0$). 

Let us write $e(S) = \min\{\pd S, \pd IS\}$ and $e^*(S) = \min\{\id S, \id PS\}$ for
any simple module $S$, thus $\finpro A = \max_S e(S) = \max_S 
e^*(S).$ We show in sections 4 and 5 that for any natural number $z$, there is a bijection $h_z$
between the set of simple modules $S$ with $e(S) = z$ and the set of simple modules $S'$
with $e^*(S') = z.$ These bijections $h_z$ can be combined in order to obtain a
permutation $h$ of the simple modules. If $S$ is a simple module, then 
$$
 h(S) = h_{e(S)}(S) = \top \Omega^{e(S)} NS,
$$
where $NS = S$ if $e(S)$ is odd and $NS = IS$ if $e(S)$ is even (see also
the description of $h$ in terms of $\tau,\ \phi$ and $\gamma$ mentioned in 4.4).
We call $h = h_A$ the {\it 
homological permutation} for $A$. 
	\bigskip
{\bf Outline of the paper.} Section 2 deals with arbitrary artin algebras. We
will review the definition of the 
delooping level of a module, as introduced by G\'elinas. We will show in 2.3:
If $X$ is a submodule of a finitely generated module $Y$, then $\del X \le \pd Y.$
As a consequence, we get: if any simple module is
a submodule of a finitely generated module of finite projective dimension, then 
$\del A \le \finpro A$, see Theorem 2.4.

Sections 3 to 5 and the appendices restrict 
the attention to Nakayama algebras. Following Madsen [M],
we show that if $A$ is Nakayama, then any simple module $S$ is a submodule of
a finitely generated module of finite projective dimension, see 3.6.
As a consequence, we can apply the general considerations of section 2, see 3.7. 
The assertions (1) to (8) in 3.7 may be
considered as the heart of the first part of the paper.

The second part of the paper, the sections 4 and 5, are devoted to
the homological permutation $h$ of a Nakayama algebra. 
The simple modules $S$ with $e(S) = 0$
are just the torsionless simple modules, and $h$ sends the
socle of an indecomposable projective-injective module $B$ to the top of $B$. 
For $z\ge 1$, the bijection $h_z$ comes equipped with long exact sequences
(the corresponding ties or bonds),
namely minimal projective resolutions of $T$ or $IT$ and 
minimal injective coresolutions of $S$ or $PS$, where $h_z(T) = S.$
As we will see in 4.4,
the bijections $h_z$ can be described in terms of the Auslander-Reiten translations $\tau$ and 
$\tau^-$, and the functions $\soc P(-)$ and $\top I(-)$. 

There are several appendices which expand the discussion. 
In Appendix A, we consider the functions $\psi = \tau^-\top I(-)$ 
and $\gamma  = \tau\soc P(-)$ on the set of simple modules for a cyclic Nakayama algebra. 
We use covering theory and deal with corresponding monotone 
endo\-functions of $\Bbb Z$. In this appendix A, we also will draw the attention to the
resolutuion quiver and the coresolution quiver of a Nakayama algebra. 

In Appendix B, we provide an outline of Sen's $\epsilon$-construction which allowed him
to prove several interesting results using induction, for example the equality of the 
finitistic dimension of a Nakayama algebra $A$ and its opposite algebra [S2], but also,
very recently, the equality $\finpro A = \del A$, see [S3]. In the center of this approach
lies the fact that there are canonical bijections between several sets of modules. All these
sets have cardinality $r$, where $r$ is the minimal number of relations needed to define $A$. 
Appendix C is devoted to these bijections. 

In the final Appendix D, we discuss two further invariants of a simple module $S$ 
which are related to the finitistic dimension, and which are considered by G\'elinas, 
namely $\grade S$ as well as the set of
numbers $t$ such that $\mho^t\Omega^tS$ is torsionless. 

Let us stress that the discussions in sections 3, 4 and 5 and in the Appendices A, B, C are
purely combinatorial. Of course, the terminology which we use comes from the theory of
rings and modules, and concerns homological and categorical invariants. However, for a 
Nakayama algebra $A$, these invariants can be expressed in combinatorial
terms, say starting from the Auslander-Reiten quiver of $A$ (this is a translation quiver).
In particular, the homological permutation defined in 4.2 should be seen in the light of 
(Catalan) combinatorics, see 4.4.
	\medskip
{\bf Acknowledgment.} The author wants to thank Vincent G\'elinas, Emre Sen and Xiao-Wu Chen
as well as a referee
for fruitful suggestions
concerning the presentation of the results. 
	\bigskip\bigskip
{\bf 2. The delooping level, the desuspending level.}
	\medskip 
{\bf 2.1. Definitions.} Let $M$ be a finitely generated module. 
The {\it delooping level\ } 
$\del M$ of $M$ is the minimal number $d\ge 0 $ such that $\Omega^d M$ 
belongs to $\add {}_AA\oplus \Omega^{d+1}M'$ for some finitely generated module $M'$.
(Note that if $\Omega^d M$ 
belongs to $\add {}_AA\oplus \Omega^{d+1}M'$ for some module $M'$,
then $\Omega^{d+1} M$ 
belongs to $\add {}_AA\oplus \Omega^{d+2}M'$.) This definition is due to G\'elinas [Ge].
As we have mentioned already, the {\it delooping level} of the algebra $A$ is defined as
$\del A = \max_S \del S$, where $S$ are the simple modules.

There is the dual concept of the {\it desuspending level\ } $\des M$ of a 
finitely generated module
$M$, this is the minimal number $d\ge 0$ such that $\Sigma^d M$ belongs to 
$\add D(A_A)\oplus \Sigma M'$ for some finitely generated module $M'$. 
Also, we put $\des A = 
\max_S \del S$, where $S$ are the simple modules. Clearly, $\des A = \del A^{\op}.$

In addition, let $\Finpro A$ be the supremum of $\pd N$, where $N$ has
finite projective dimension, but is not necessarily finitely generated
($\Finpro A$ is usually called the {\it big finitistic dimension}), 
as well as $\Fininj A$, the supremum of $\id N$, where $N$ has
finite injective dimension, but not necessarily finitely generated. 
	\bigskip
{\bf 2.2. Proposition (G\'elinas).}
{\it Let $M$ be a finitely generated module. 
If there is a (not necessarily finitely generated) module $N$ with finite
injective dimension $d\ge 1$ such that $\Ext^d(M,N) \neq 0$, then 
$d \le \del M.$}
	\medskip
Proof.
Let $N$  be a  module with $d = \id N$ and $1 \le d  < \infty.$
Then $\Ext^{d+1}(X,N) = 0$ for all modules $X$. 
Let $M$ be a finitely generated 
module and assume that $\del M < d.$ By definition, $\Omega^{d-1} M$ belongs to 
$\add {}_AA\oplus\Omega^d M',$ for some finitely generated module $M'$.
It follows that $\Ext^d(M,N) = \Ext^1(\Omega^{d-1}M,N)$
belongs to $\add\Ext^1({}_AA\oplus\Omega^d M',N)$, but 
$$
 \Ext^1({}_AA\oplus \Omega^d M' ,N) =  \Ext^1(\Omega^d M',N) = \Ext^{d+1}(M',N) = 0,
$$
and therefore $\Ext^d(M,N) = 0$, a contradiction. 
$\s$
	\medskip
{\bf Corollary:} {\it $\Fininj A \le \del A.$} 
	\medskip
Proof. We have to show: If $N$ is an arbitrary module with 
finite injective dimension $d$, then 
$d \le \del S $ for some simple module $S$. This is clear for $d = 0.$ 
Thus, we can assume that $d\ge 1.$ Since $\id N = d,$
there is a simple module $S$ with $\Ext^d(S,N) \neq 0.$ The Proposition
yields $d\le \del S.$  $\s$
	\bigskip
{\bf 2.2$'$. The dual assertions.} 
	\medskip
{\bf Proposition.}
{\it Let $M$ be a finitely generated module. 
If there is a (not necessarily finitely generated) module $N$ with finite
projective dimension $d\ge 1$ which satisfies $\Ext^d(N,M) \neq 0$, then 
$d \le \des M.$}  $\s$
	\medskip
{\bf Corollary:} {\it $\Finpro A \le \des A.$} $\s$
	\bigskip
{\bf 2.3. Proposition.} 
{\it If $X$ is a submodule of a finitely generated module
$Y$, then}
$$
  \del X \le \pd Y.
$$
	\medskip
Proof. Let $\pd Y = d$. We can assume that $d < \infty$. We apply the horseshoe lemma
to the exact sequence $0 \to X \to Y \to Y/X \to 0$ and 
obtain an exact sequence 
$$
 0 \to \Omega^d X \to P'\oplus \Omega^d Y\to \Omega^d (Y/X) \to 0,
$$ 
where $P'$ is projective.  
Since $\Omega^d Y$ is projective, we see that the middle term $E = 
P'\oplus \Omega^d Y$ is projective and the
projective cover $P(\Omega^d(Y/X))$ of $\Omega^d(Y/X)$ is a direct summand of $E$,
say $E = P''\oplus P(\Omega^d(Y/X))$
for some projective module $P''$. Therefore 
$$
 \Omega^d X = P''\oplus \Omega(\Omega^d (Y/X)) = P''\oplus \Omega^{d+1}(Y/X).
$$ 
This shows that $\del X \le d.$
$\s$
	\medskip
{\bf Remark.} {\it The inequality
$$ 
 \del X \le \min\{\pd Y \mid Y \text{\ finitely generated and\ }X \subseteq Y\}
$$ 
may be proper.}
For example,
if $A$ is representation-finite, and $X$ is any finitely generated module, then
$\del X$ is, of course, finite. 
On the other hand, 
if $X$ is an injective module with infinite projective dimension,
then $X \subseteq Y$ implies that also $\pd Y$ is infinite (since $X$ is a direct
summand of $Y$), thus the right hand side is $\infty$. A typical example of a 
representation finite algebra with an injective module with infinite projective
dimension is the radical-square-zero algebra $A$ whose quiver has two vertices,
say $1$ and $2$, and two arrows: a loop at $1$ and an arrow $2 \to 1$. The simple
module $S(2)$ is injective and has infinite projective dimension (since $\Omega S(2) =
S(1)$ and $\Omega S(1) = S(1)$), Also, there are just five isomorphism classes of
indecomposable modules. $\s$
	\bigskip
{\bf 2.3$'$. The dual assertion.}
	\medskip
{\bf Proposition.} 
{\it If $X$ is a factor module of a finitely generated module
$Y$. Then}
$$
  \des X \le \id Y.
$$ 
\vglue-.5cm
$\s$
	\bigskip
{\bf 2.4. Theorem.} {\it Assume that every simple module $S$ is a submodule of a 
finitely generated module $M_S$ of finite projective dimension. Let $d = \max_S \pd M_S$.
Then}
$$
     \Fininj A \le \del A \le d \le \finpro A.
$$
	\medskip
Proof. The first inequality is Corollary 2.2.
According to Proposition 2.3, we have $\del S \le \pd M_S$, thus
$\del A = \max_S \del S \le \max \pd M_S = d.$ Finally, we have of course
$\pd M_S \le \finpro A,$ thus $d = \max \pd M_S \le \finpro A.$
$\s$
	\medskip
Recall that an algebra is called {\it left Kasch} provided any simple module 
occurs as a left ideal. Thus, $A$ is left Kasch if and only if any
simple module is a submodule of a module of projective dimension $d = 0.$
The case $d = 0$ of 2.4 is therefore just the well-known assertion that 
$\Fininj A = 0$ for any Kasch algebra $A$.
	\bigskip
{\bf 2.4$'$. The dual assertion.} 
	\medskip
{\bf Theorem.} {\it Assume that every simple module $S$ is a factor module of a 
finitely generated module $N_S$ of finite injective dimension. Let $d' = \max_S \id N_S$.
Then}
$$
  \Finpro A \le \des A \le d' \le \fininj A.
$$
\vskip-1truecm
$\s$
	\bigskip
{\bf 2.5. Combination.} 
{\it Assume that every simple module $S$ is a submodule of a 
finitely generated module $M_S$ of finite projective dimension and also  
a factor module of a 
finitely generated module $N_S$ of finite injective dimension. 
Let $d = \max_S \pd M_S,$ and  $d' = \max_S \id N_S$.
Then}
$$
   \finpro A =  \fininj A =  \Finpro A = \Fininj A = \del A = \des A = d = d'.
$$
\vskip-1truecm
$\s$
	\bigskip
{\bf 2.6. Remark.} Clearly,
2.5 applies to a large class of algebras, but the remaining parts of the paper
are devoted just to Nakayama algebras. As one of the referees has stressed, 
2.5 can be used to obtain a direct proof of the following result (see G\'elinas [G1], 1.18); it
concerns the so called Gorenstein algebras of Gorenstein dimension $d$: 

{\it Let $A$ be an artin algebra with minimal injective cogenerator
$Q$. Assume that $\pd Q = d < \infty$ and $\id {}_AA = d' < \infty$. Then}
$$
   \finpro A =  \fininj A =  \Finpro A = \Fininj A = \del A = \des A = d = d'.
$$
Namely, every simple module $S$ is a submodule of $M_S = Q$ and a  
factor module of $N_S = {}_AA$. By assumption,  $\pd Q = d$ and $\id {}_AA = d'$. $\s$
	\bigskip
{\bf 3. Nakayama algebras.}
	\medskip
Let $A$ be a Nakayama algebra. 
Let us note that if $M$ is an indecomposable module, then $\Omega M$ and
$\Sigma M$ are either zero or indecomposable again. 
An indecomposable module $M$ of finite projective dimension will be said to be
{\it even,} or {\it odd} provided its projective dimension is even, or odd, 
respectively. 
	\medskip
A connected Nakayama algebra $A$ is called {\it linear} provided there is a 
simple projective module, otherwise $A$ is called {\it cyclic.}
	\medskip
{\bf 3.1. The maps $\psi$ and $\gamma$.} We denote by $D\Tr$ and $\Tr D$ the Auslander-Reiten
translations in $\mod A$. 
If $M$ is indecomposable and not projective, we write $\tau M = D\Tr M$; if
$M$ is indecomposable and not injective, we write $\tau^- M = \Tr D M$
(let us stress that we apply $\tau$ only to modules
which are indecomposable and not projective, and $\tau^-$ only to modules which are
indecomposable and not injective).

If $S$ is a simple module, let $\psi S = \tau^-\top IS$ provided $\top IS$ is not injective.
There is the dual concept: If $S$ is a simple module, let $\gamma S = \tau\soc PS$
provided $\soc PS$ is not projective. 
[The functions $\psi$ and $\gamma$ are important tools for dealing in particular with cyclic 
Nakayama algebras; their study goes back to Gustafson [Gu]. They were used in many papers. 
The appendix A of this paper will focus the attention to further properties of 
$\psi$ and $\gamma$. 
In particular, one can define the coresolution quiver (or $\psi$-quiver) of 
a cyclic Nakayama algebra $A$.
The $\psi$-paths considered in 3.4 are just the paths in the
coresolution quiver. Similarly, using $\gamma$, we obtain the resolution quiver (or
$\gamma$-quiver) of a cyclic Nakayama $A$.] 
	\medskip
{\bf Lemma.} (a) {\it A simple module $U$ is in the image of $\psi$ if and only if $\pd U \ge 2$.
A simple module $T$ is in the image of $\gamma$ if and only if $\id T \ge 2.$}
	\smallskip
(b) {\it If $M$ is an indecomposable module, and $m\in \Bbb N$.  Then either $\pd M < 2m$ or else
$\top \Omega^{2m}M = \gamma^m \top M$. Similarly, either $\id M < 2m$ or else
$\soc \Sigma^{2m}M = \psi^m \soc M.$}
	\medskip 
Proof. (a) We show the second assertion (the first follows by duality).

Let $T$ be simple with $\id T \ge 2$. Then $\Sigma T$ is not injective, thus we have a proper inclusion
$\Sigma T \subset I\Sigma T.$ We denote by $S$ the top of $I\Sigma T$.  
Now all the indecomposable modules $M$ which properly include $\Sigma T$ have to be projective. In particular, 
$I\Sigma T$ is projective and therefore $I\Sigma T = PS.$ 
Since $\soc PS = \soc \Sigma T$ is not projective, we can apply $\tau$.
By definition, $\gamma S = \tau\soc PS = 
\tau \soc I  \Sigma T = \tau\tau^-T = T.$ This shows that $T$ is in the image of $\gamma.$

Conversely, let $S$ be any simple module and $T = \gamma S$. We want to show that $\Sigma T$ is
not injective. Now $\Sigma T$ and $PS$ have the same socle, thus both are submodules of
the projective module $IPS$. Since $IPS$ is serial, its submodules are pairwise
comparable. If $PS \subseteq \Sigma T$, then $IT$ has a submodule $X$ of length $|PS|+1$
which maps onto $PS$. But this is impossible since $PS$ is projective. It follows that there
is the proper inclusion $\Sigma T \subset PS$. This shows that $\Sigma T$ is not
injective, thus $\id T \ge 2.$
	\smallskip
(b) See [R2], section 3, Corollary. 
$\s$  
	\bigskip
{\bf 3.2.} Nakayama algebras have been investigated quite thoroughly by
Madsen [M]. We will use many of his results, in particular:
	\medskip
{\bf Lemma  (Madsen).} {\it Let $A$ be any Nakayama algebra and $M$ an indecomposable module.
\item{\rm (a)} If $X$ is a subfactor of $M$, and $M$ is odd, then $X$ is odd and 
 $\pd X \le \pd M.$ 
\item{\rm (b)} If $X$ is a subfactor of $M$, and $X$ is even, then $M$ is even and 
 $\pd M \le \pd X.$
\item{\rm (c)} A simple module $S'$ is a composition factor of
$\Omega^2S$ if and only if $\psi(S') = S$ (and then the multiplicity of 
$S'$ in $\Omega^2S$ is $1$).\par} 
	\medskip
Proof. For (a) and (b), see [M], 2.2. For (c), see [M] 3.1 (note that 
the assumption that $A$ is a cyclic Nakayama algebra is not needed). 
$\s$
	\medskip
{\bf Corollary.} {\it Let $X \subseteq Y$ be indecomposable modules. If $Y$ is odd,
then $\pd X \le \pd Y$. If $Y$ is even, then $\pd IX \le \pd Y.$}
	\medskip
Proof. The first assertion is a special case of (a). Thus, assume that $Y$ is even. We
look at the inclusion $Y \subseteq IY = IX$ and use (b): it asserts that $\pd IX \le \pd Y.$
$\s$ 
	\bigskip
{\bf 3.3. Maximum property of odd modules.} {\it Let $A$ be any Nakayama algebra. 
An indecomposable module $M$ is odd if
and only if all composition factors of $M$ are odd;  and then $\pd M$ is the maximum of
$\pd S$, where $S$ is a composition factor of $M$.}
	\medskip
Proof. Let $M$ be indecomposable. First, assume that $M$ is odd. According to
3.2 (a), all composition factors $S$ of $M$ are odd and $\pd S \le \pd M$.
Of course, at least one of the composition factors must have $\pd S = \pd M$ (since
the class of modules of projective dimension smaller $\pd M$ is closed under
extensions). Thus $\pd M$ is odd. 
Conversely, assume that all composition factors of $M$ are odd. According to [M]
Proposition 4.1, it follows that $\pd M$ is the maximum of
$\pd S$, where $S$ is a composition factor of $M$.
$\s$
	\medskip
{\bf Remark.} There is the following consequence: {\it Let $A$ be a Nakayama algebra. For any
odd number $1 \le i \le \finpro A$, there is a simple module $S$ with $\pd S = i.$}
(Note that the odd simple modules of a cyclic Nakayama algebra are nicely displayed in the 
$\psi$-quiver; see Appendix A.8 and Proposition 3.5.)
	\medskip
Proof. There is a (finitely generated) module $M$ with $\pd M = d = \finpro A$. Thus, for any
number $0 \le i \le d$, there is an indecomposable module $M_i$ with $\pd M_i = i$,
namely a direct summand of $\Omega^{d-i}M$. If $i$ is odd, 3.3 asserts that $M_i$ has
a composition factor $S$ with $\pd S = i.$ $\s$
	\bigskip
The next two subsections 3.4 and 3.5 are devoted to cyclic Nakayama algebras.
	\medskip
{\bf 3.4. The $\psi$-paths (and the $\gamma$-paths).} Let $\Lambda$ be a cyclic Nakayama algebra.
A {\it $\psi$-path} of cardinality $m$ is a sequence $(S_1,S_2,\dots,S_m)$ of simple 
modules with $S_{i+1} = \psi S_i$ for $1\le i < m,$ it starts in $S_1$ and ends in $S_m$.
If $\psi(S) = T$, we say that $S$ is a {\it $\psi$-predecessor} of $T$,
If $\gamma(S) = T$, we say that $S$ is a {\it $\gamma$-predecessor} of $T$,
Similarly, a {\it $\gamma$-path} of cardinality $m$ 
is a sequence $(S_1,S_2,\dots,S_m)$ of simple 
modules with $S_{i+1} = \gamma S_i$ for $1\le i < m,$ it starts in $S_1$ and ends in $S_m$.
	
For $A$ a cyclic Nakayama algebra, and $S$ a simple module, we define $a(S)$ and $a'(S)$ as
follows: Let $a(S)$ be the supremum of the cardinality of the 
$\psi$-paths ending in $S$, and dually, let $a'(S)$ be the supremum of the cardinality of the 
$\gamma$-paths ending in $S$.
Thus, $a(T)$ is the maximum of $1+a(S)$ with $S$ a $\psi$-predecessor of $T$
and $a'(T)$ is the maximum of $1+a'(S)$ with $S$ a $\gamma$-predecessor of $T$.
	
We say that a simple module $S$ is {\it $\psi$-cyclic} provided $a(S) = \infty$,
or, equivalently, provided there is some number $e\ge 1$ with $\psi^e S = S.$
Similarly, $S$ is {\it $\gamma$-cyclic} provided $a'(S) = \infty$, thus
provided there is some number $e\ge 1$ with $\gamma^e S = S.$
	
Finally, let $a(A)$ be the maximum of $a(S)$, where $S$ is simple and not $\psi$-cyclic. 
Note that $a(A)$ is also 
the maximum of $a'(S)$, where $S$ is simple and not $\gamma$-cyclic, see A.1 in the Appendix.
	
If $a(A) = 0$, then all simple modules are torsionless, thus all modules
are torsionless, thus $A$ is selfinjective. For the topics discussed in the paper,
we usually may assume that $A$ is not self-injective. 
	\medskip
{\bf 3.5. Proposition.} {\it Let $A$ be a cyclic Nakayama algebra, and $S$ a simple module. 

\item{\rm (1)}
We have $a(S)$ finite if and only if $\pd S$ is odd, and then $\pd S = 2a(S)-1$,
in particular, $\pd S < 2a(S).$
If $a(S)$ is infinite, then $\pd IS$ is even and $\pd IS \le 2a(A).$

\item{\rm (2)}
Similarly, $a'(S)$ finite if and only if $\id S$ is odd, and then $\id S = 2a'(S)-1$,
in particular, $\id S < 2a(S).$
If $a'(S)$ is infinite, then $\id PS$ is even and $\id PS \le 2a(A).$\par}

	\medskip 
Proof. First, we show: {\it If $a(S)$ finite, then $\pd S = 2a(S)-1.$}
We use induction on $a(S)$. If $a(S) = 1$, then $S$ has no 
$\psi$-predecessor, thus $\rad PS$ is projective, thus $\pd S = 1$.
Now assume that $a(S) \ge 2.$ Then $S$ has a $\psi$-predecessor, thus
$\rad PS$ is not projective, and therefore $\pd PS \ge 2$. As a consequence,
$\Omega^2 S \neq 0.$  According to 3.2 (c),
the composition factors of $\Omega^2 S$ are the $\psi$-predecessors
$S'$ of $S$. If $S'$ is a $\psi$-predecessor of $S$, then $a(S') \le a(S)-1$,
and by induction $\pd S' = 2a(S')-1 \le 2a(S)-3.$ And there is a predecessor
$S''$ of $S$ with $a(S'') = a(S)-1$, thus with $\pd S'' = 2a(S'')-1 = 2a(S)-3.$ 
This shows that the composition factors of $\Omega^2 S$ have odd 
projective dimension with maximum $2a(S)-3$. The maximum property 3.3
asserts that $\pd \Omega^2 S = 2a(S)-3.$
Thus, $\pd S = 2+\pd \Omega^2 S = 2a(S)-1.$ 
	
Conversely, let us show: {\it If $\pd S$ is odd, then 
$a(S)$ is finite.} The proof is again by induction. If $\pd S = 1$, then
$S$ cannot have a $\psi$-predecessor, thus $a(S) = 1$.
If $\pd S$ is odd and not $1$, then $\pd \Omega^2 S = \pd S-2$ is odd, again.
Let $S'$ be a predecessor of $S$. By 3.2 (c), $S'$ is a composition factor of
$\Omega^2 S$. By the maximum property 3.3, $\pd S'$ is odd and 
$\pd S' \le \pd \Omega^2 S = \pd S -2.$ By induction, $a(S')$ is finite.
Since $a(S')$ is finite for all predecessors $S'$ of $S$, also $a(S)$ is finite.

Now assume that $a(S)$ is infinite. 
If $IS$ is projective, then its projective dimension is zero, therefore even. 
Thus, we
suppose that $IS$ is not projective and show that $\Omega IS$ has odd projective
dimension. The following picture displays part (of the universal covering) of the
Auslander-Reiten quiver of $A$, with arrows going from left to right:
$$
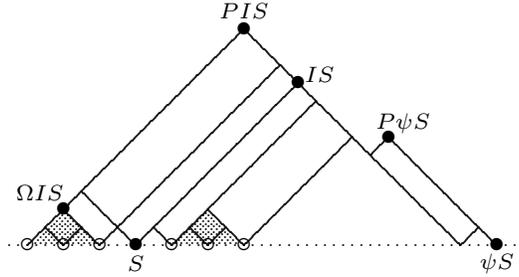
{\beginpicture
    \setcoordinatesystem units <.2cm,.2cm>
\multiput{} at 0 0  26 12  /
\plot 0 0  12 12  24 0  25 1 /
\plot 2 2  4 0  14 10 /
\plot 1 1  2 0  3 1 /
\plot 5 1  6 0  15 9 /
\plot 7 1  8 0  16 8 /
\plot 9 1  10 0  11 1 /
\plot 10 2  12 0  18 6 /
\plot 3 3  5 1 /
\setshadegrid span <.3mm>
\vshade 0 0 0  <z,z,,> 2 0 2  <z,z,,> 4 0 0 /
\vshade 8 0 0  <z,z,,> 10 0 2  <z,z,,> 12 0 0 /
\put{$\ss S$} at 6 -1 
\plot 19 5  20 6  26 0 /
\put{$\ss \psi S$} at 26 -1 
\put{$\ss IS$} at 16.2 9.5 
\put{$\ss \Omega IS$} at .7 2.9  
\put{$\ss PIS$} at 12 13   
\multiput{$\bullet$} at 6 0  15 9  2 2  26 0  12 12  20 6 /
\put{$\ss P\psi S$} at 20.8 6.8    
\multiput{$\circ $} at 0 0  2 0  4 0  8 0  10 0 12 0  /
\setdots <1mm>
\plot -1 0  28 0 /
\endpicture}
$$
On the left and on the right of $S$ (at the lower boundary of he shaded areas) 
are the $\psi$-predecessors $S'$ of $S$ which are different from $S$.
For these modules $S'$, we have $a(S') < \infty$, thus, 
as we have seen already, 
they have odd projective dimension. 
The module $\Omega IS$ has a filtration using such
modules, thus 3.3 asserts that $\pd\Omega IS$ is odd, therefore 
$\pd IS$ is even. Since $\pd \Omega IS \le 2a(A)-1$, we have $\pd IS \le 2a(A).$
This completes the proof of (1).

The assertions (2) follow by duality. 
$\s$
	\bigskip
{\bf 3.6. Theorem.}
{\it Let $A$ be any Nakayama algebra and $S$ a simple module.
Then $S$ is odd or $IS$ is even.} (Thus, at least one of $S, IS$ has finite
projective dimension.) 
	\medskip
Proof. We can assume that $A$ is connected. 
We can assume that $S$ is not odd, thus $\pd S$ is even or infinite.
If $\pd S$ is even, then Lemma 3.2 (b) asserts that $\pd IS$ is even.
If $\pd S$ is infinite, then $A$ has to be cyclic, and 
3.5 (a) asserts first that $a(S)$ cannot be finite and then
that $\pd IS$ is even.
$\s$
	\medskip
{\bf Dual assertion.} 
{\it Let $A$ be any  Nakayama algebra and $S$ a simple module.
Then $\id S$ is odd or $\id PS$ is even.} (Thus, at least one of $S, PS$ has finite
injective dimension.) $\s$
	\bigskip
{\bf 3.7. Summary.} 
	\medskip
{\bf (1)} {\it Let $A$ be a Nakayama algebra. Then  
any simple module is
a submodule of an indecomposable module with finite projective dimension,
and a factor module of an indecomposable module with finite injective dimension.}
	\medskip
Proof. Let $S$ be simple. 
Theorem 3.6 asserts that $\pd S$ or $\pd IS$ is finite, thus $S$ is a submodule of an
indecomposable module with finite projective dimension. 
By duality, $S$ or $PS$ has finite injective dimension, thus $S$ is a factor module of 
an indecomposable module with finite injective dimension.
$\s$
	\medskip
{\bf (2) Theorem.}
{\it Let $A$ be a Nakayama algebra. If $S$ is simple, let $M_S = S$ provided $S$ has finite
projective dimension, and $M_S = IS$ otherwise,
and let $N_S = S$ provided $S$ has finite injective dimension, and $N_S = PS$ otherwise. 
Let $d = \max_S \pd M_S$ and $d' = \max_S \id N_S.$ Then}
$$ 
    \finpro A = \fininj A = \del A = \des A = d = d'.
$$
	\medskip
Proof. We need to know that all the modules $M_S$ have finite projective dimension, 
and all the modules $PS$ have finite injective dimension. Then we can apply 2.5.
Theorem 3.6 and its dual 
assert that all the modules $M_S$ have finite
projective dimension and that all the modules $N_S$ have finite injective dimension.
$\s$

	\medskip
{\bf (3)} {\it Let $A$ be a Nakayama algebra and $X$ an indecomposable module. 
let us write}
$$
 e(X) = \min\{\pd X, \pd IX\}.
$$

(a) {\it We have
$$
 e(X) = \min\{\pd Y \mid Y \text{\ \it finitely generated and\ }X \subseteq Y\} 
$$
and}
$$
 \del X \le e(X).
$$

(b) {\it If $S$ is simple, then}
$$ 
 e(S) < \infty.
$$
	\medskip
Proof. (a) 
Let $X \subseteq Y$, where $Y$ is finitely generated and $\pd Y$ minimal.
Since $X$ has simple socle, we can assume that $Y$ is indecomposable. Namely, let $Y = \bigoplus Y_i$
with all $Y_i$ indecomposable. Then $X$ is a submodule of some $Y_i$ and $\pd Y_i \le \pd Y.$

If $Y$ is odd, then $X\subseteq Y$ shows that $\pd X \le  \pd Y$, see 3.2 (a). Thus the 
minimality of $\pd Y$ yields $\pd X = \pd Y$. If $Y$ is even, then $Y \subseteq IY = IX$ shows that
$\pd IX \le \pd Y$, see 3.2 (b). Thus, the minimality of $\pd Y$ yields $\pd IX = \pd Y.$ 
The inequality $\del X \le e(X)$ has been shown in 2.3. 

(b) For the assertion $e(S) < \infty$, see Theorem 3.6. $\s$ 

	\medskip
Since $e(X) =  \min\{\pd Y \mid Y \text{\ \it finitely generated and\ }X \subseteq Y\}$, 
G\'elinas has suggested to call the
number $e(X)$ the {\it embedding projective dimension} of $X$.
	\bigskip
{\bf (4) Theorem.} {\it Let $A$ be a Nakayama algebra. Then 
$$
 \fin\pro A = \max\nolimits_S e(S)
$$
were $S$ runs through the simple modules.}

	\medskip
Proof. As we have mentioned in (3), 
we have $\del X \le e(X)$ for all modules $X$. Thus $\del A = \max_S \del S \le
\max_S e(S)$. If $S$ is simple, then, 
according to (1), $\pd S$ or $\pd IS$ is finite, thus 
$e(S) = \min\{\pd S, \pd IS\} \le \finpro A.$
This shows that
$$
 \finpro A = \del A = \max\nolimits_S \del S \le \max\nolimits_S e(S)
 \le \finpro A.
$$
\vglue-0.5cm
$\s$
	\bigskip
{\bf (5)} We can rewrite the definition of $e(S)$ as follows:
$$
{\beginpicture
    \setcoordinatesystem units <.8cm,.4cm>
\multiput{} at 0 0   /
\put{$e(S) = \left\{ \matrix \min\{\pd S, \pd IS\} \cr
                       \pd S \cr
                       \pd IS 
                     \endmatrix 
              \right.$} at 0 0 
\put{if} at 3.2 0
\put{$S$ is odd, $IS$ is even,\strut} [l] at 4 1
\put{$IS$ is not even,\strut} [l] at 4 0
\put{$S$ is not odd.\strut} [l] at 4 -1
\endpicture}
$$
Note that if $IS$ is not even, then $S$ has to be odd, see Theorem 3.6, thus 
the definition covers all possible cases for the pair $(S,IS)$ and is well-defined.
	\medskip
Proof. Let $S$ be simple. 
If $IS$ is not even, then $\pd S \le \pd IS$
(this is clear, if $\pd IS = \infty$, whereas for $IS$ odd, we use the first part of
Corollary 3.2), therefore 
$\pd S = \min\{\pd S,\pd IS\} = e(S).$
If $S$ is not odd, then $\pd IS \le \pd S$
(this is clear, if $\pd S = \infty$, whereas. if $S$ is even, 
then we use the second part of 
Corollary 3.2). Therefore
$\pd IS = \min\{\pd S,\pd IS\} = e(S).$ 
$\s$
	\medskip 
{\bf  Corollary 1.} {\it Let $A$ be a Nakayama algebra. If $S$ is simple, let}
$$
{\beginpicture
    \setcoordinatesystem units <.8cm,.4cm>
\multiput{} at 0 0   /
\put{$f(S) = \left\{ \matrix 0 \cr
                       \pd S \cr
                       \pd IS 
                     \endmatrix 
              \right.$} at 0 0 
\put{if} at 2.2 0
\put{$S$ is odd, $IS$ is even,\strut} [l] at 3 1
\put{$IS$ is not even,\strut} [l] at 3 0
\put{$S$ is not odd.\strut} [l] at 3 -1
\endpicture}
$$
{\it Then $f(S) = e(S)$ for all simple modules $S$ such that $IS$ is not even or 
$S$ is not odd, thus $f(S)\le e(S)$ for all simple modules $S$. Also, }
$$
 \finpro A = \max\nolimits_S f(S).
$$
	\medskip
Proof of Corollary 1.
As we have seen at the beginning of (5), $f(S) = e(S)$ if $IS$ is not even or $S$ is not odd. 
If  $S$ is odd and $IS$ is even, then $f(S) = 0 \le 
\min\{\pd S,\pd IS\} = e(S).$
It follows that $f(S) \le e(S)$ for all simple modules $S$. 

In order to show that $\fin\pro A = \max_S f(S),$ we note that
on the right hand side of the formula (4), we may delete all the terms
$\min\{\pd S,\pd IS\}$, where $S$ is a simple module which is odd such that $IS$ is even. Namely, 
if $S$ is odd and $IS$ is even, then $\pd S \neq \pd IS$,
thus $\min\{\pd S,\pd IS\} < \max \{\pd S,\pd IS\} \le \finpro A$. 
$\s$ 
	\medskip
{\bf Corollary 2.} {\it Let $A$ be a Nakayama algebra and $z = \finpro A.$
If $z$ is odd, then there is $S$ simple with $\pd S = z$ and  $\pd IS$ odd or infinite.
If $z$ is even, then there is $S$ simple with $\pd IS = z$ and  $\pd S$ even or infinite.}
	\medskip
The proof follows from Corollary 1 and Theorem 3.6 as follows.
If all simple modules are odd, then all non-zero modules are odd by the maximum
principle, but this is impossible, since the projective modules are even. 
Thus $\finpro A = \pd S$ for some $S$ with $IS$ not even (we call it the first case), or
$\finpro A = \pd IS$ for some $S$ with $S$ not odd (the second case).
In the first case, $\pd IS$ is odd or infinite. Also, since $IS$ is not even,
$S$ has to be odd, thus $z = \pd S$ is odd, according to 3.6.
In the second case, $\pd S$ is even or infinite. Since $S$ is not odd, 3.6 asserts that
$z = \pd IS$ is even. Conversely, of course, if $z$ is odd, we cannot be in the second case,
thus we are in the first case. And if $z$ is even, then we have to be in the second case. 
$\s$
	\bigskip
{\bf (6) (Madsen).}  
{\it Let $A$ be a connected Nakayama algebra. If there exists a simple module $S$ such that
$S$ is even or $IS$ is odd, then $A$ has finite global dimension.}
	\medskip
Proof. If $\pd S$ is even, then see Madsen [M], Theorem 3.1.
Thus assume that $IS$
is odd. Then all factor modules of $IS$ are odd, according to 3.2 (a).
Assume that $X$ is a non-zero factor module of $IS$
which is not injective, thus there is an inclusion $X \subset Y$ with $Y/X$ simple.
Then $Y$ is projective and $\pd Y/X = 1+\pd X$ is even. Thus, we haven a simple module
of even projective dimension and as mentioned already,  
$A$ has finite global dimension. We therefore may assume that all factor modules
of $IS$ are injective. In particular, the simple factor module of $IS$ is injective,
thus $A$ is a linear Nakayama algebra and thus has finite global dimension. 
$\s$
	\medskip
{\bf Corollary.}
{\it  Let $A$ be a connected Nakayama algebra with infinite global dimension.
Let $z = \finpro A.$
If $z$ is odd, then there is $S$ simple with $\pd S = z$ and  $\pd IS$  infinite.
If $z$ is even, then there is $S$ simple with $\pd IS = z$ and  $\pd S$  infinite.}
	\medskip
Proof: This follows from (6) and Corollary 2 in (5). Namely,
Since $A$ has infinite global dimension, $\pd IS$ cannot be odd and $\pd S$ cannot
be even.
$\s$
	\bigskip
{\bf (7)} 
Looking at the finitistic dimension of an algebra, the algebras of infinite
global dimension are those of real interest. For algebras of infinite
global dimension, the function $f$ defined in (5) can be defined slightly
different:
	\medskip 
{\it Let $A$ be a connected Nakayama algebra with infinite global
dimension. Then the function $f$ considered in {\rm (5)} coincides with the following
function $g$:}
$$
{\beginpicture
    \setcoordinatesystem units <.8cm,.4cm>
\multiput{} at 0 0   /
\put{$g(S) = \left\{ \matrix 0 \cr
                       \pd S \cr
                       \pd IS 
                     \endmatrix 
              \right.$} at 0 0 
\put{if} at 2.2 0
\put{$\pd S$, $\pd IS$ both are finite,} [l] at 3 1
\put{$\pd IS$ is infinite,} [l] at 3 0
\put{$\pd S$ is infinite.} [l] at 3 -1 
\endpicture}
$$
(according to Theorem 3.6, at least one of $\pd S$ and $\pd IS$ is finite).
	\medskip
{\bf Corollary.} {\it Let $A$ be a connected Nakayama algebra with infinite global
dimension. Then }
$$
 \finpro A = \max\nolimits_S g(S).
$$
\vglue-0.5truecm
 $\s$
	\medskip
{\bf Remark.} If $A$ has finite global dimension, then the function $g$ is the zero
function. Thus, in this case,  $g$ cannot be used for calculating $\finpro A$
($\finpro A$ is just the global dimension of $A$, and the global dimension of 
a cyclic Nakayama algebra is at least 2).
	\medskip
Using (6), we obtain a short proof for the following
result ([S2], Corollary 4.13):
	\medskip
{\bf (8) (Sen)} {\it  Let $A$ be a connected Nakayama algebra which is 
Gorenstein and has infinite global dimension. Then $\finpro A$ is even.}
	\medskip
Proof. If 
$\finpro A$ is odd, the Corollary in (6) provides 
an indecomposable injective module $I$
with $\pd I$ infinite, a contradiction to the assumption that $A$ is Gorenstein.
$\s$
	\bigskip
{\bf 3.8. Examples. The possible pairs $(\pd S, \pd IS).$} 
We have seen in 3.2: (a) if $IS$ is odd, then $S$ is odd and 
$\pd S \le \pd IS$;\ (b) if $S$ is even, then $IS$ is even and $\pd S \ge \pd IS$.
Theorem 3.6 shows that it is impossible that both $S$ and $IS$ have infinite
projective dimension. 
{\it All the remaining possibilities do occur.} Here are some examples (for a general
discussion, we refer to section 4.7).
	\medskip
\item{$\star$} $\pd S < \pd IS$, both odd: see $S = S_7.$
\item{$\star$} $\pd S = \pd IS$, both odd: see $S = S_6.$
\item{$\bullet$} $\pd S$ odd, $\pd IS = \infty$: see $S = S_3.$
\item{$\bullet$} $\pd S < \pd IS$, with  $\pd S$ odd, $\pd IS$ even;
   see $S = S_8$ and $S = S_{10}.$
\item{$\bullet$} $\pd S > \pd IS$, with  $\pd S$ odd, $\pd IS$ even: 
   see $S = S_1$ and $S = S_9.$
\item{$\bullet$} $\pd S = \infty$, $\pd IS$ even: see $S = S_4.$
\item{$\star$} $\pd S > \pd IS$, both even: see $S = S_5.$
\item{$\star$} $\pd S = \pd IS$, both even: see $S = S_2.$
	\smallskip
\noindent
According to 3.7 (6), the cases marked by a star $\star$ can occur only if $A$ has
finite global dimension. For the two cases which occur both for finite and for
infinite global dimension, we have provided two examples.
	\medskip
The following pictures exhibit Auslander-Reiten quivers (the
Auslander-Reiten quiver of a cyclic Nakayama algebra 
lives on a cylinder --- the dashed lines left and right have to be
identified). A vertex of an Auslander-Reiten quiver is the isomorphism class $[M]$
of an indecomposable module $M$. But instead of drawing the vertex $[M]$, 
we insert here the 
corresponding value $\pd M$ (this is a natural number or $\infty$). The lower boundary
of the Auslander-Reiten quiver of a Nakayama algebra consists of the simple
modules: those which are of interest here, are encircled and labeled $S_1,\dots, S_{10}$.

$$
{\beginpicture
    \setcoordinatesystem units <.4cm,.4cm>

\put{\beginpicture
\multiput{} at 0 -1  6 3 /
\setdots <1mm>
\plot 0 0  6 0 /
\setdots <.5mm>
\plot 0 1  3 4  6 1  /
\plot 0 1  1 0  4 3 /
\plot 1 2  3 0  5 2 /
\plot 2 3  5 0  6 1 /

\multiput{$\ss 0$} at 1 2  2 3  3 4   /
\multiput{$\ss 1$} at 0 1  1 0  5 0  6 1 /
\multiput{$\ss 2$} at 2 1  3 0  3 2  4 1  4 3  5 2  /
\setdashes <1mm>
\plot 0 -0.5  0 2 /
\plot 6 -0.5  6 2 /
\put{$S_1$} at 5 -1
\put{$S_2$} at 3 -1
\put{$S_8$} at 1 -1
\multiput{$\bigcirc$} at 1 0  3 0  5 0 /
\endpicture} at 5 0
\put{\beginpicture
\multiput{} at 0 -1  6 3.5 /
\setdots <1mm>
\plot 0 0  6 0 /
\setdots <.5mm>
\plot  0 0  2 2  4 0  5 1  6 0 /
\plot  1 1  2 0  3 1 /
\multiput{$\ss 0$} at 1 1  2 2  5 1 /
\multiput{$\ss 1$} at  4 0 /
\multiput{$\ss 2$} at 0 0  6 0  /
\multiput{$\ss 3$} at 2 0  3 1 /
\setdashes <1mm>
\plot 0 0.5  0 2 /
\plot 6 0.5  6 2 /
\put{$S_5$} at 0 -1
\put{$S_6$} at 2 -1
\multiput{$\bigcirc$} at 0 0  2 0 / 

\endpicture} at 15 0

\put{\beginpicture
\multiput{} at 0 -1  8 3.5 /
\setdots <1mm>
\plot 0 0  8 0 /
\setdots <.5mm>
\plot  0 0  2 2  4 0  6 2 /
\plot  1 1  2 0  5 3  8 0 /
\plot 4 2  6 0  7 1 /
\multiput{$\ss 0$} at 1 1  2 2  4 2  5 3  /
\multiput{$\ss 1$} at  0 0  4 0  8 0 /
\multiput{$\ss 2$} at 2 0   3 1    /
\multiput{$\ss 3$} at 5 1  6 0  6 2  7 1    /
\setdashes <1mm>
\plot 0 0.5  0 2 /
\plot 8 .5  8 2 /
\put{$S_7$} at 4 -1
\put{$\bigcirc$} at 4 0 

\endpicture} at 25 0 
\put{\beginpicture
\multiput{} at 0 -2  4  2   /
\setdots <1mm>
\plot 0 -1  4 -1 /
\setdots <.5mm>
\plot  0 0  2 2  4 0  /
\plot  1 1  2 0  3 1 /
\plot 0 0  1 -1  2 0  3 -1  4 0 /
\multiput{$\ss 0$} at 1 1  2 2 /
\multiput{$\ss 1$} at 1 -1 /
\multiput{$\ss \infty$} at 0 0  2 0  4 0  3 1  3 -1 /
\setdashes <1mm>
\plot 0 0.5  0 2 /
\plot 4 0.5  4 2 /
\plot 0 -0.5  0 -1.5 /
\plot 4 -0.5  4 -1.5 /
\put{$S_3$} at 1 -2
\put{$S_4$} at 3 -2
\multiput{$\bigcirc$} at 1 -1  3 -1  /
\endpicture} at 7 -8

\put{\beginpicture
\multiput{} at 0 -1  6 3 /
\setdots <1mm>
\plot 0 -3  6 -3 /
\setdots <.5mm>
\plot 0 1  3 4  6 1  /
\plot 0 1  1 0  4 3 /
\plot 1 2  3 0  5 2 /
\plot 2 3  5 0  6 1 /
\plot 0 -1  1 0  2 -1  3 0  4 -1  5 0  6 -1 /
\plot 0 -1  1 -2  2 -1  3 -2  4 -1  5 -2  6 -1 /
\plot 0 -3  1 -2  2 -3  3 -2  4 -3  5 -2  6 -3 /
\multiput{$\ss 0$} at 1 2  2 3  3 4   /
\multiput{$\ss 1$} at 2 -3  3 -2  4 -3 /
\multiput{$\ss 2$} at 2 1  3 0  3 2  4 1  4 3  5 2  /
\multiput{$\ss \infty$} at 0 1  6 1  1 0   5 0  0 -1  2 -1  4 -1  6 -1  1 -2  5 -2  0 -3  6 -3 /
\setdashes <1mm>
\plot 0 -3.5  0 2 /
\plot 6 -3.5  6 2 /
\put{$S_{10}$} at 4 -4
\put{$S_{9}$} at 2 -4
\multiput{$\bigcirc$} at  2 -3  4 -3 /
\endpicture} at 17 -7

\endpicture}
$$
	\medskip
As we have mentioned, if $S$ is $\psi$-cyclic, then $\pd S \ge \pd IS$. But 
$\pd S \ge \pd IS$ may also happen for $S$ odd (and then $\pd S > \pd IS$). 
It happens, of course, in
case $S$ is torsionless, since then $\pd S  > 0$ and $\pd IS = 0.$ 
Here is an example
of $S$, with $\pd S = 3,\ \pd IS = 2,$ thus $S$ is 
odd, not torsionless, and $\pd S > \pd IS.$
$$
{\beginpicture
    \setcoordinatesystem units <.4cm,.4cm>
\multiput{} at 0 -1  12 2.2    /
\setdots <1mm>
\plot 0 0  12 0 /
\setsolid
\setdots <0.5 mm>
\plot 0 0  2 2  4 0  6 2  8 0  10 2  12 0 /
\plot 1 1  2 0  3 1 /
\plot 5 1  6 0  8 2  10 0  11 1 /
\setdashes <1mm>
\plot 0 0.5  0 1.5 /
\plot 12 0.5  12 1.5 /
\multiput{$\ss 0$} at 1 1  2 2  5 1  6 2  8 2  10 2 /
\multiput{$\ss 1$} at 4 0  8 0 /
\multiput{$\ss 2$} at 6 0  7 1  11 1 /
\multiput{$\ss 3$} at 9 1  10 0 /
\multiput{$\ss 4$} at 0 0  12 0 /
\multiput{$\ss 5$} at 2 0  3 1 /
\multiput{$\bigcirc$} at 0 0  6 0  12 0 /
\put{$S$} at 10 -.8 
\endpicture}
$$
	\medskip
In addition, we see the following: Proposition 3.5 asserts 
that if $a(S)$ is infinite, then $\pd IS$ is even, but the converse does not hold: 
$\pd IS$ may be even, whereas $a(S)$ is finite, as the examples
$S_1$ and $S_8$ show. 
	\bigskip
{\bf 3.9. Remarks.}
We do not know whether the inequality $\del S \le e(S)$
can be proper. But
{\it there is always a simple module $S$ with 
$\del S = \finpro A$, and if $\del S = \finpro A$,
then $\del S = e(S).$} \ 
Proof. The first assertion follows directly from $\max_S \del S = \finpro A$.
If $ \del S = \finpro A$, then the inequalities
$$ 
  \del S \le e(S) \le \finpro A
$$
yield the equality $\del S = \min\{\pd S,\pd IS\}.$
$\s$
	\medskip
One also may compare $f(S)$ as introduced in 3.7 (5) with $\del S.$ 
For the simple module $S_8$ mentioned in 3.8, we have $f(S_8) = 0$, since $S_8$
is odd and $IS_8$ is even. On the other hand, since $S_8$ is not torsionless, 
$\del S_8 \ge 1$ (since $\pd S_8 = 1$, we have $\del S_8 = 1$).

	\bigskip
{\bf 3.10. Proposition.} {\it Let $A$ be a cyclic Nakayama algebra. Then}
$$
 2a(A)-1 \le \finpro A \le 2a(A).
$$
	\medskip
Proof.
In order to show that $2a(A)-1 \le \finpro A$, 
take a simple module $S$ with $a(S) = a(A).$ According to 3.5 (1), we have
$\pd S = 2a(S)-1 = 2a(A)-1.$

For the second inequality, let $d = \finpro A$. 
There is a module $M$ with $\pd M = d$. 
If $d$ is odd, then 3.3 yields a composition
factor $S$ of $M$ with $\pd S = d.$ 
According to 3.5 (1), we have $\pd S = 2a(S)-1.$ 
Therefore $d = \pd S = 2a(S)-1 \le 2a(A)-1 < 2a(A).$ 

If $d$ is even and at least 2, then $\pd \Omega M = d-1$
and 3.3 yields a compositon factor of $\Omega M$ with $\pd S = d-1.$ 
According to 3.5 (1), we have $\pd S = 2a(S)-1.$  
Therefore $d -1 = \pd S = 2a(S)-1 \le 2a(A)-1,$ thus $d\le 2a(A).$
$\s$
	\bigskip
{\bf 3.11. Historical Remarks.} Some
important properties of Nakayama algebras which we use in
order to apply the general theory are collected in 3.5. One should
be aware that both assertions of 3.5 can be found (at least implicitly) in the literature.
The odd modules have been studied quite carefully by Madsen. For the calculation 
of $\pd IS$, where $S$ is a $\psi$-cyclic simple module, we should also refer to Shen,
see Lemma 3.5 (2) of [Sh3].  
	\bigskip

{\bf 4. The homological permutation.}
	\medskip
Let $A$ be a Nakayama algebra and $\Cal E$ the set of (isomorphism classes of)
simple modules. We are going to construct a permutation $h\:\Cal E \to \Cal E$
related to projective resolutions and injective coresolutions.
The essential step in the construction is the following 
Theorem.
	\medskip
{\bf 4.1. Theorem.} {\it Let $A$ be a Nakayama algebra. Let $S, T$ be simple modules
and $z\in \Bbb N_0$.} 

(a) {\it If $z$ is even, $\Omega^z IT = PS$, and $\pd T \ge z$, then $\Sigma^z PS = IT$
and $\id S \ge z.$}

(b) {\it If $z$ is odd, $\Omega^z T = PS$, and $\pd IT \ge z$, then $\Sigma^z S = IT$
and $\id PS \ge z.$}
	\medskip
The proof will be given in section 5 (first, 
in case $A$ is a linear Nakayama algebra, see 5.1;
then, using covering theory, also for $A$ being a cyclic Nakayama algebra, see 5.4). 
	\bigskip
{\bf 4.2. Reformulations.}
Recall that if $S$ is a simple module, then $e(S) = \min\{\pd S,\pd IS\}$. 
Dually, we have defined $e^*(S) = \min\{\id S,\id PS\}$. For any natural number $z$,
let $\Cal E(z)$ be the set of simple modules $S$ with $e(S) = z$, and
$\Cal E^*(z)$ the set of simple modules $S$ with $e^*(S) = z$. 
	\medskip
If $S$ is a simple module, let 
$$
\align
  NS &= \left\{ \matrix S & \text{\ \ if $e(S)$ is odd,} \cr
              IS & \text{\quad if $e(S)$ is even.} \endmatrix \right.  \cr
  N^*S &= \left\{ \matrix S & \text{\ \, if $e^*(S)$ is odd,} \cr
              PS & \text{\quad if $e^*(S)$ is even.} \endmatrix \right.  \cr
\endalign
$$
Also, let
$$
 h(S) = \top \Omega^{e(S)} NS, \quad \text{and} \quad 
 h^*(S) = \soc \Sigma^{e^*(S)} N^*S.
$$
	\medskip
{\bf Corollary.} {\it The map $h$ is a permutation of $\Cal E$ with inverse $h^*$. If 
$S$ is a simple module, then}
$$
 e^*(h(S)) = e(S).
$$
We call $h = h_A$ the {\it homological permutation} of $\Cal E$.
	\medskip
Note that for $A$ being self-injective, we have $h(S) = \top IS$ and
$h^*(S) =  \soc PS$, for any simple module $S$. Thus, in this case, {\it $h$ is 
the Nakayama permutation of $A$,} as considered by Nakayama himself. 
(Recall that Nakayama [N, Lemma 2], has characterized the 
self-injective algebras $B$ 
by the existence of a permutation $h$ of the set of simple modules 
such that $\soc P(h(S)) = S$, or, equivalently, that $h(S) = \top IS$.
This permutation $h$, but sometimes also its inverse $h^*$,
is usually called the Nakayama permutation of $B$.) 
	\medskip 
We can reformulate the Corollary as follows: {\it Let $z$ be a natural number.
If $S$ belongs to $\Cal E(z)$, then $h(S)$ belongs to $\Cal E^*(z)$,
and $h$ yields a bijection $h_z\:\Cal E(z) \to \Cal E^*(z)$ with inverse 
given by $h^*$.} In particular,
it follows that {\it the sets $\Cal E(z)$ and $\Cal E^*(z)$ have the same cardinality.}
	\medskip
Proof. 4.1 asserts that $h$ maps $\Cal E(z)$ into $\Cal E^*(z)$. By duality, 
$h^*$ maps $\Cal E^*(z)$ into $\Cal E(z)$. We refer again to 4.1 in order to see that
$h^*h(S) = S$ and $hh^*(S) = S$ for all $S\in \Cal E$. 
$\s$
	\bigskip
{\bf 4.3. The ties.} 
The permutation $h$ is, of course, determined by the set of pairs $(T,h(T))$ with $T\in
\Cal E$. For such a pair $(T,S) = (T,h(T))$, the number $e(T) = e^*(S)$
plays an important role. Let us stress that the pairs $(T,S) = (T,h(T))$ come equipped with 
addditional information, the corresponding ties: 
If $e(T) = 0$, thus also $e^*(S) = 0$, there is given 
an indecomposable module which is both projective and injective, namely $IT = PS$.
If for $e(T) > 0$, thus $e^*(S) > 0$, there are given two long exact sequences. 
Namely, if $e(T) = 2t \ge 2$ is even, there is given a minimal projective resolution
of $IT$ and a minimal injective coresolution of $PS$ (both are long exact sequences
starting with $PS$ and ending with $IT$), and we use, in addition, the inclusion map 
$T \to IT$ and the projection $PS \to S$:
$$
{\beginpicture
    \setcoordinatesystem units <1.5cm,.7cm>
\put{\beginpicture
\multiput{$\ss 0$\strut} at 0.3 0  5.7 0 /
\put{$\ss S$\strut} at 1 -1
\put{$\ss T$\strut} at 4.9 1
\put{$\ss IT$\strut} at 4.9 0 
\put{$\ss P_0(IT)$\strut} at 4 0
\put{$\cdots$\strut} at 3.1 0
\put{$\ss P_{2t-1}(IT)$\strut} at 2 0 
\put{$\ss \Omega^{2t}IT$\strut} at 1 0
\arr{0.45 0}{0.65 0}
\arr{1.3 0}{1.5 0}
\arr{2.5 0}{2.7 0}
\arr{3.4 0}{3.6 0}
\arr{4.4 0}{4.6 0}
\arr{5.2 0}{5.4 0}
\arr{1 -.2}{1 -.7}
\arr{4.9 .7}{4.9 .2}
\multiput{$\bigcirc$\strut} at 1 -1  4.9 1 /
\endpicture} at 0 0
\put{\beginpicture
\multiput{$\ss 0$\strut} at 0.3 0  5.7 0 /
\put{$\ss S$\strut} at 1 -1
\put{$\ss T$\strut} at 4.9 1
\put{$\ss \Sigma^{2t}PS$\strut} at 4.95 0 
\put{$\ss I_{2t-1}(PS)$\strut} at 3.85 0
\put{$\cdots$\strut} at 2.9 0
\put{$\ss I_{0}(PS)$\strut} at 1.9 0 
\put{$\ss PS$\strut} at 1 0
\arr{0.45 0}{0.65 0}
\arr{1.25 0}{1.45 0}
\arr{2.3 0}{2.5 0}
\arr{3.2 0}{3.4 0}
\arr{4.35 0}{4.55 0}
\arr{5.3 0}{5.5 0}
\arr{1 -.2}{1 -.7}
\arr{4.9 .7}{4.9 .25}
\multiput{$\bigcirc$\strut} at 1 -1  4.9 1 /
\endpicture} at 0 -2
\endpicture}
$$
If $e(T)$ is odd, there is given
a minimal projective resolution of $T$ (thus a long exact sequence
starting with $PS$ and ending with $T$), as well as 
a minimal injective coresolution of $S$ (thus a long exact sequence
starting with $S$ and ending with $IT$),
and we use, in addition, the projection $PS \to S$ and the inclusion map 
$T \to IT$.
$$
{\beginpicture
    \setcoordinatesystem units <1.5cm,.7cm>
\put{\beginpicture
\multiput{$\ss 0$\strut} at 0.3 0  5.7 0 /
\put{$\ss S$\strut} at 1 -1
\put{$\ss T$\strut} at 4.9 0 
\put{$\ss P_0(T)$\strut} at 4 0
\put{$\cdots$\strut} at 3.1 0
\put{$\ss P_{2t-1}(T)$\strut} at 2 0 
\put{$\ss \Omega^{2t}T$\strut} at 1 0
\arr{0.45 0}{0.65 0}
\arr{1.3 0}{1.5 0}
\arr{2.5 0}{2.7 0}
\arr{3.4 0}{3.6 0}
\arr{4.4 0}{4.6 0}
\arr{5.2 0}{5.4 0}
\arr{1 -.2}{1 -.7}
\multiput{$\bigcirc$\strut} at 1 -1  4.9 0 /
\endpicture} at 0 0
\put{\beginpicture
\multiput{$\ss 0$\strut} at 0.3 0  5.7 0 /
\put{$\ss T$\strut} at 4.9 1
\put{$\ss \Sigma^{2t}S$\strut} at 4.95 0 
\put{$\ss I_{2t-1}(S)$\strut} at 3.85 0
\put{$\cdots$\strut} at 2.9 0
\put{$\ss I_{0}(S)$\strut} at 1.9 0 
\put{$\ss S$\strut} at 1 0
\arr{0.45 0}{0.65 0}
\arr{1.25 0}{1.45 0}
\arr{2.3 0}{2.5 0}
\arr{3.2 0}{3.4 0}
\arr{4.35 0}{4.55 0}
\arr{5.3 0}{5.5 0}
\arr{4.9 .7}{4.9 .25}
\multiput{$\bigcirc$\strut} at 1 0  4.9 1 /
\endpicture} at 0 -1.4
\endpicture}
$$

Alternatively, we can say that 
for any pair $(T, h(T))$, there are given specified walks in the Auslander-Reiten quiver, 
starting with
$T$ and ending in $h(T)$: For $e(T) = 0$, there is given 
the composition of the sectional path $T \to IT$
with the sectional path $IT = PS \to S$. For $e(T) > 0$, 
there is given the walk from $T$ to  $S = h(T)$ 
via the long exact sequences which we have mentioned. 
Some examples are outlined in 4.5 and 4.6. 
A more detailled discussion of these ties for small values $z = e(T)$ will be found in 4.7. 
	\medskip
{\bf Remark.} Let us stress the following: Let $T$ be a simple module.
{\it If $z = e(T)$ is odd, then clearly $\Ext^z(T,h(T)) \neq 0$}
(since $h(T) = \top\Omega^z T$
and $\Ext^z(T,\top \Omega^zT) = \Ext^z(T,\Omega^z T)$). {\it However, for $z = e(T)$ even, we may have $\Ext^z(T,h(T)) = 0.$} Here is a typical example with $z = 2$: Take for $T$
the simple module with index 4, thus $h(T)$ is the simple module with index $2$ and
clearly $\Ext^2(T,h(T)) = 0.$
$$
{\beginpicture
    \setcoordinatesystem units <.4cm,.4cm>
\put{\beginpicture
\multiput{} at 0 0  8 2  /
\setdots <1mm>
\plot 0 0  8 0 /
\setsolid 
\plot 0 0  2 2  4 0  6 2  8 0 /
\plot 1 1  2 0  4 2  6 0  7 1  /
\put{$\ss 1$} at 0 -1
\put{$\ss 2$} at 2 -1
\put{$\ss 3$} at 4 -1
\put{$\ss 4$} at 6 -1
\put{$\ss 5$} at 8 -1
\put{$\bullet$} at 6 0 
\endpicture} at 0 0 
\put{$h_A = (1,3,5)(2,4)$} at 11 0 
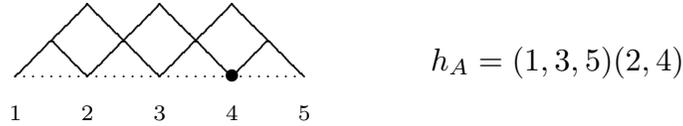
\endpicture}
$$
	\bigskip
{\bf 4.4. Combinatorial description.} Let us repeat the steps of the algorithm 
which yields $h_A$ for a given
Nakayama algebra $A$. Let $\Gamma$ be the Auslander-Reiten quiver of $A$. The
south-east arrows of $\Gamma$ correspond to irreducible maps which are epimorphisms, and the
north-east arrows correspond to irreducible maps which are monomorphisms. 

(1) Given a vertex
$x$, the {\bf top} $\top x$ of $x$ is the end of the longest south-east path starting in $x$; the 
{\bf projecive cover} $Px$ of $x$ is the start of the longest south-east path $p_x$ ending in $x$;
the {\bf injective envelope} $Ix$ of $x$ is the end of the longest north-east path starting in $x$.
If $x$ is not a projective vertex, the {\bf syzygy module} $\Omega x$ is the kernel of $p_x$ 
(thus there is a north-east path from $\Omega x$ to
$Px$ of length at least $1$). 

(2) Using inductively $P$ and $\Omega$, we can define the function 
$\pd\:\Gamma_0 \to \Bbb N_0\cup\{\infty\}$.  
For any simple module $S$, define $e(S) = 
\min\{\pd S,\pd IS\}$. 

(3) Now, for any simple module $S$,
let $NS = S$ in case $e(S)$ is odd, and $NS = IS$ in case $e(S)$ is even.
Finally let $h_A(S) = \top\Omega^{e(S)}NS$.
	\medskip 
Once we know the values $e(S)$, for $S$ simple, 
there is an alternative description of $h_A$, 
using the Auslander-Reiten translation $\tau$ as well as the 
functions  $\phi = \top I(-)$ and $\gamma = \tau\soc P(-)$
(but note that it is our convention to apply $\tau$ only to indecomposable 
{\bf non-projective}
modules, thus $\tau$ and $\gamma$ are not always defined). 
	\medskip
{\bf Proposition.} 
{\it If $S$ is a simple module with $e(S) = 2t$, then $h(S) = \gamma^t\phi(S)$.
If $S$ is a simple module with $e(S) = 2t+1$, then $h(S) = \gamma^t\tau(S)$.}
	\medskip
Proof. If $\pd IS = 0$, we have $e(S) = 0$ and $h(S) = \top IS = \phi S = h_0(S).$
If $\pd IS \ge 1,$ then also $\pd S \ge 1$ and $\top\Omega S = \top \Omega IS = \tau S.$
First, let $e(S) = 2t$. Then $h(S) = \top \Omega^z IS = \gamma_t\top IS = \gamma^t\phi S.$
Second, let $e(S) = 2t+1$. Then $h(S) = \top\Omega^{2t}\Omega S = \top\Omega^{2t}\tau S
= \gamma^t\tau S.$
$\s$
	\medskip
We write $h_z = h|\Cal E(z)$. Thus, the Proposition asserts that for all $t\ge 0$, 
$$
 h_{2t}  = \gamma^t\phi\ |\ \Cal E(2t), \quad \text{and} \quad
 h_{2t+1}  = \gamma^t\tau\ |\ \Cal E(2t+1).
$$
	\medskip
{\bf Remark.} 
Let us assume that $T$ is simple and that $z = \pd T = \pd IT$. {\it If $z$ is odd,
then $\Omega^zT = \Omega^z IT$, whereas for $z$ even, we many have $\Omega^zT \neq \Omega^zIT$.}
Namely, if $z = 2t+1$, then
$$
 \top \Omega^zT = \gamma^t\top \Omega T \quad\text{and}\quad
 \top \Omega^zIT = \gamma^t\top \Omega IT,
$$
and $\top \Omega T = \tau T = \top \Omega IT.$ On the other hand, 
as in  Remark 4.3, we look at the example
$$
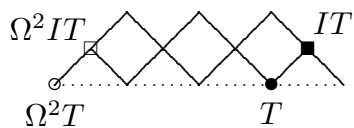
{\beginpicture
    \setcoordinatesystem units <.4cm,.4cm>
\put{\beginpicture
\multiput{} at 0 0  8 2  /
\setdots <1mm>
\plot 0 0  8 0 /
\setsolid 
\plot 0 0  2 2  4 0  6 2  8 0 /
\plot 1 1  2 0  4 2  6 0  7 1  /
\put{$T$} at 6 -.8
\put{$IT$} at 7.7 1.7
\put{$\Omega^2 T$} at 0 -.8
\put{$\Omega^2 IT$} at -.2 1.6
\put{$\bullet$} at 6 0 
\put{$\sssize\blacksquare$} at 7 1 
\put{$\circ$} at 0 0
\put{$\sssize \square$} at 1 1 
\endpicture} at 0 0 
\endpicture}
$$
Note that for $\pd T =  \pd IT$ odd, a minimal projective resolution of $IT$
as well as a minimal injective coresolution of $P(h(T))$ provide additional interesting
ties between $T$ and $h(T)$. 

	\bigskip
{\bf 4.5. Linear Nakayama algebras. Some examples and remarks.}
We want to provide some examples which 
illustrate the homological permutation and the corresponding ties, and we will add some
remarks. 
	\medskip
We start with the following example of a linear Nakayama algebra $A$. The simple modules
are labeled $1,2,\dots,5$
$$
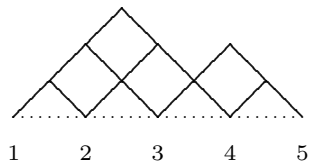
{\beginpicture
    \setcoordinatesystem units <.4cm,.4cm>
\multiput{} at 0 0  8 3  /
\setdots <1mm>
\plot 0 0  8 0 /
\setsolid 
\plot 0 0  3 3  6 0  7 1 /
\plot 1 1  2 0  4 2 /
\plot 2 2  4 0  6 2  8 0 /
\put{$\ss 1$} at 0 -1
\put{$\ss 2$} at 2 -1
\put{$\ss 3$} at 4 -1
\put{$\ss 4$} at 6 -1
\put{$\ss 5$} at 8 -1
\endpicture}
$$
Let us indicate ties for the simple modules $T$. Here, $T$ is always marked by a bullet;
for $T = 2$ and $T = 4$, we draw just one of the walks in the Auslander-Reiten quiver (the
one which uses a minimal projective resolution of $T$).
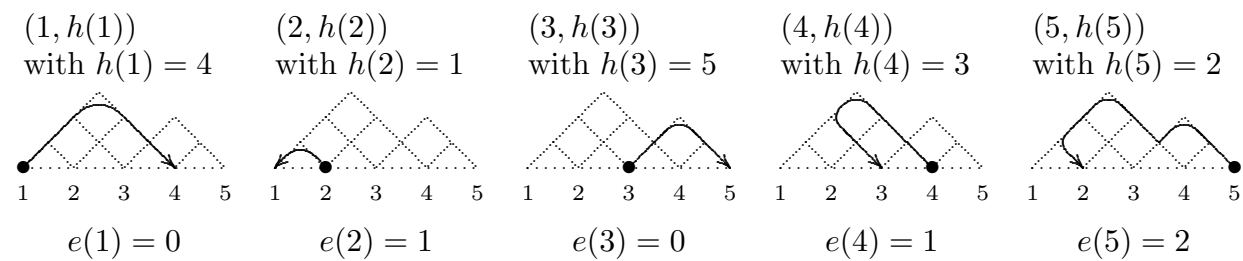
$$
{\beginpicture
    \setcoordinatesystem units <.28cm,.28cm>
\put{\beginpicture
\multiput{} at 0 0  8 3  /
\put{$(1,h(1))$}  [l] at 0 5.5
\put{with $h(1) = 4$} [l] at 0 4
\setdots <1mm>
\plot 0 0  8 0 /
\setdots <.5mm>
\plot 0 0  3 3  6 0  7 1 /
\plot 1 1  2 0  4 2 /
\plot 2 2  4 0  6 2  8 0 /
\put{$\ss 1$} at 0 -1
\put{$\ss 2$} at 2 -1
\put{$\ss 3$} at 4 -1
\put{$\ss 4$} at 6 -1
\put{$\ss 5$} at 8 -1
\setsolid
\setquadratic
\plot 0 0  1 1   2 2  3 2.5  4 2  5 1  6 0 /
\arr{5.9 0.1}{6 0}
\put{$\bullet$} at 0 0 
\put{$e(1) = 0$}  at 4 -3
\endpicture} at 0 0
\put{\beginpicture
\multiput{} at 0 0  8 3  /
\put{$(2,h(2))$}  [l] at 0 5.5
\put{with $h(2) = 1$} [l] at 0 4

\setdots <1mm>
\plot 0 0  8 0 /
\setdots <.5mm>
\plot 0 0  3 3  6 0  7 1 /
\plot 1 1  2 0  4 2 /
\plot 2 2  4 0  6 2  8 0 /
\put{$\ss 1$} at 0 -1
\put{$\ss 2$} at 2 -1
\put{$\ss 3$} at 4 -1
\put{$\ss 4$} at 6 -1
\put{$\ss 5$} at 8 -1
\setsolid
\setquadratic
\plot 2 0  1 0.7  0 0  /
\arr{0.1 0.1}{0 0}
\put{$\bullet$} at 2 0 
\put{$e(2) = 1$}  at 4 -3
\endpicture} at 10 0

\put{\beginpicture
\multiput{} at 0 0  8 3  /
\put{$(3,h(3))$}  [l] at 0 5.5
\put{with $h(3) = 5$} [l] at 0 4
\setdots <1mm>
\plot 0 0  8 0 /
\setdots <.5mm>
\plot 0 0  3 3  6 0  7 1 /
\plot 1 1  2 0  4 2 /
\plot 2 2  4 0  6 2  8 0 /
\put{$\ss 1$} at 0 -1
\put{$\ss 2$} at 2 -1
\put{$\ss 3$} at 4 -1
\put{$\ss 4$} at 6 -1
\put{$\ss 5$} at 8 -1
\setsolid
\setquadratic
\plot 4 0  5 1  5.5 1.5  6 1.7  6.5 1.5  7 1  8 0 /
\arr{7.9 0.1}{8 0}
\put{$\bullet$} at 4 0 
\put{$e(3) = 0$}  at 4 -3
\endpicture} at 20 0
\put{\beginpicture
\multiput{} at 0 0  8 3  /
\put{$(4,h(4))$}  [l] at 0 5.5
\put{with $h(4) = 3$} [l] at 0 4
\setdots <1mm>
\plot 0 0  8 0 /
\setdots <.5mm>
\plot 0 0  3 3  6 0  7 1 /
\plot 1 1  2 0  4 2 /
\plot 2 2  4 0  6 2  8 0 /
\put{$\ss 1$} at 0 -1
\put{$\ss 2$} at 2 -1
\put{$\ss 3$} at 4 -1
\put{$\ss 4$} at 6 -1
\put{$\ss 5$} at 8 -1
\setsolid
\setquadratic
\plot 6 0  5 1       3.5 2.5  
           3  2.7    2.5 2.5  
           2.2 2   2.5 1.5  
           3 1  4 0 /
\arr{3.9 0.1}{4 0}
\put{$\bullet$} at 6 0 
\put{$e(4) = 1$}  at 4 -3
\endpicture} at 30 0
\put{\beginpicture
\multiput{} at 0 0  8 3  /
\put{$(5,h(5))$}  [l] at 0 5.5
\put{with $h(5) = 2$} [l] at 0 4
\setdots <1mm>
\plot 0 0  8 0 /
\setdots <.5mm>
\plot 0 0  3 3  6 0  7 1 /
\plot 1 1  2 0  4 2 /
\plot 2 2  4 0  6 2  8 0 /
\put{$\ss 1$} at 0 -1
\put{$\ss 2$} at 2 -1
\put{$\ss 3$} at 4 -1
\put{$\ss 4$} at 6 -1
\put{$\ss 5$} at 8 -1
\setsolid
\setquadratic
\plot 8 0  7 1  6.5 1.5  6 1.7  5.5 1.5  5.2 1.2  5 1  /
\plot  5 1  4 2   3.5 2.5 
           3  2.7  2.5 2.5  2 2 
  1.5 1.5  
           1.2 1   1.5 0.5  
           1.7 0.3  2 0 /
\arr{1.9 0.1}{2 0}
\put{$\bullet$} at 8 0 
\put{$e(5) = 2$}  at 4 -3

\endpicture} at 40 0
\endpicture}
$$
The permutation $h$ is, in cycle notation, the permutation $(14352).$
	\medskip
Here are the permutations $h_A$ which occur for a linear Nakayama algebra $A$ with 4 simple
modules.
$$
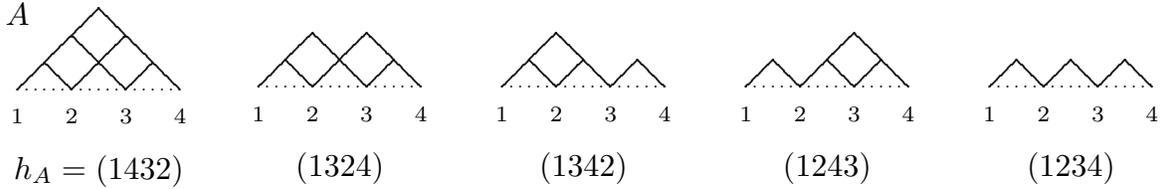
{\beginpicture
    \setcoordinatesystem units <.3cm,.3cm>
\put{\beginpicture
\put{$A$} at 0 2.8
\multiput{} at 0 0  6 3 /
\setdots <1mm>
\plot 0 0  6 0 /
\setsolid
\plot 0 0  3 3  6 0 /
\plot 1 1  2 0  4 2 /
\plot 2 2  4 0  5 1 /
\put{$\ss 1$} at 0 -1
\put{$\ss 2$} at 2 -1
\put{$\ss 3$} at 4 -1
\put{$\ss 4$} at 6 -1
\put{$h_A = (1432)$} at 3 -3 
\endpicture} at 0 0
\put{\beginpicture
\multiput{} at 0 0  6 3 /
\setdots <1mm>
\plot 0 0  6 0 /
\setsolid
\plot 0 0  2 2 /
\plot 4 2   6 0 /
\plot 1 1  2 0  4 2 /
\plot 2 2  4 0  5 1 /
\put{$\ss 1$} at 0 -1
\put{$\ss 2$} at 2 -1
\put{$\ss 3$} at 4 -1
\put{$\ss 4$} at 6 -1
\put{$(1324)$} at 3 -3 
\endpicture} at 9 0
\put{\beginpicture
\multiput{} at 0 0  6 3 /
\setdots <1mm>
\plot 0 0  6 0 /
\setsolid
\plot 0 0  2 2  /
\plot 5 1  6 0 /
\plot 1 1  2 0  3 1  /
\plot 2 2  4 0  5 1 /
\put{$\ss 1$} at 0 -1
\put{$\ss 2$} at 2 -1
\put{$\ss 3$} at 4 -1
\put{$\ss 4$} at 6 -1
\put{$(1342)$} at 3 -3 
\endpicture} at 18 0
\put{\beginpicture
\multiput{} at 0 0  6 3 /
\setdots <1mm>
\plot 0 0  6 0 /
\setsolid
\plot 0 0  1 1 /
\plot 4 2   6 0 /
\plot 1 1  2 0  4 2 /
\plot 3 1   4 0  5 1 /
\put{$\ss 1$} at 0 -1
\put{$\ss 2$} at 2 -1
\put{$\ss 3$} at 4 -1
\put{$\ss 4$} at 6 -1
\put{$(1243)$} at 3 -3 
\endpicture} at 27 0
\put{\beginpicture
\multiput{} at 0 0  6 3 /
\setdots <1mm>
\plot 0 0  6 0 /
\setsolid
\plot 0 0  1 1  2 0  3 1  4 0  5 1  6 0 /
\put{$\ss 1$} at 0 -1
\put{$\ss 2$} at 2 -1
\put{$\ss 3$} at 4 -1
\put{$\ss 4$} at 6 -1
\put{$(1234)$} at 3 -3 
\endpicture} at 36 0
\endpicture}
$$
	\medskip
{\bf Remarks.} We start with some basic observation concerning linear Nakayama algebras $A$. 
We assume that $\Cal E = \{1,2,\cdots,n\}$
with $\tau x = x\!-\!1$ for $2\le x \le n$.
	\smallskip
(1) {\it If $e(x) = 0$, then $h_A(x) > x.$ If $e(x) > 0$, then $h_A(x) < x.$}
	\smallskip
(2) As a consequence, {\it $h_A$ is a fixed-point-free permutation of $\{1,2,\cdots,n\}$.}
	\smallskip
(3) {\it Already the pairs $(x,h_A(x))$ with $h(x) > x$
determine the algebra $A$}: namely, they
determine the peaks, thus the Kupisch series of $A$ (a {\it peak} of $A$ is by definition an
indecomposable module which is projective and injective).
In particular, $A$ is determined by $h_A$.
	\medskip
(4) Note that the linear Nakayama algebras with $n\ge 1$ simple modules correspond bijectively 
to the Dyck paths of length $2(n-1)$. Namely, the roof of $A$ is (say after a rotation by
225 degrees) a Dyck path $D_A$: 
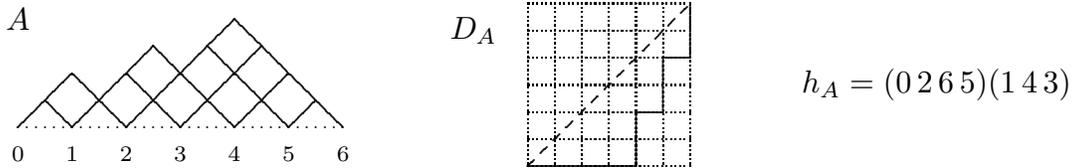
$$
{\beginpicture
    \setcoordinatesystem units <.3cm,.3cm>
\put{\beginpicture
\multiput{} at 0 0  12 4  /
\setdots <1mm>
\plot 0 0  12 0 /
\setsolid
\plot 0 0  2 2  4 0  8 4  12 0 /
\plot 1 1  2 0  5 3  8 0  10 2 /
\plot 4 2  6 0  9 3 /
\plot 7 3  10 0  11 1 /
\put{$\ss 0$} at  0 -1
\put{$\ss 1$} at  2 -1
\put{$\ss 2$} at  4 -1
\put{$\ss 3$} at  6 -1
\put{$\ss 4$} at  8 -1
\put{$\ss 5$} at  10 -1
\put{$\ss 6$} at  12 -1
\put{$A$} at 0 4
\endpicture} at 0 0 
\put{\beginpicture
\multiput{} at 0 0  6 6  /
\setdots <.5mm>
\setplotarea x from 0 to 6, y from 0 to 6
\grid {6} {6}
\setsolid 
\plot 0 0  4 0  4 2  5 2  5 4  6 4  6 6 /
\setdashes <1mm>
\plot 0 0  6 6 /
\put{$D_A$} at -2 5
\endpicture} at 14.5 0 
\put{$h_A = (0\,2\,6\,5)(1\,4\,3)$} at 28 0 
\endpicture}
$$ 
The number of Dyck paths of length $2m$ (thus the number of linear Nakayama algebras
with $m+1$ simple modules)
is the Catalan number $C_{m} = \frac1{m+1}\binom {2m}m$: 
we are in the center of Catalan combinatorics.
There are several well-known procedures to attach to Dyck paths permutations, but usually
one attaches to a Dyck path of length $2m$ a permutation of $m$ objects, not of $m+1$
objects, as we do it here. We are not aware that the correspondence $D_A \mapsto h_A$
has yet been considered in combinatorics (but the number of papers on  
Catalan combinatorics seems to be too large to be overlooked) --- if it has, our 
discussion may add at least a homological interpretation. 
	\medskip
(5) M. Rubey has added the homological permutation for linear Nakayama algebras 
to the database FindStat, see [RS]. 
	\bigskip
{\bf 4.6. Cyclic Nakayama algebras. Some examples (and remarks).} 
Let $A$ be the following cyclic Nakayama algebra with simple modules with
labels $1,2,3,4$
$$
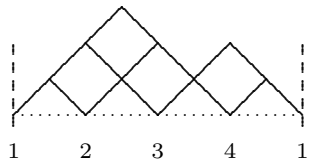
{\beginpicture
    \setcoordinatesystem units <.4cm,.4cm>
\multiput{} at 0 0  8 3  /
\setdots <1mm>
\plot 0 0  8 0 /
\setsolid 
\plot 0 0  3 3  6 0  7 1 /
\plot 1 1  2 0  4 2 /
\plot 2 2  4 0  6 2  8 0 /
\put{$\ss 1$} at 0 -1
\put{$\ss 2$} at 2 -1
\put{$\ss 3$} at 4 -1
\put{$\ss 4$} at 6 -1
\put{$\ss 1$} at 8 -1
\setdashes <1mm>
\plot 0 -.3  0 2 /
\plot 8 -.3  8 2 /

\endpicture}
$$
Here are ties for the simple modules $S$
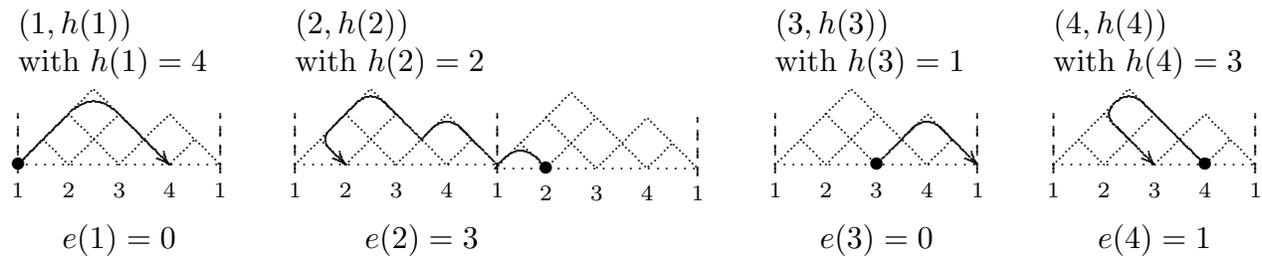
$$
{\beginpicture
    \setcoordinatesystem units <.28cm,.28cm>
\put{\beginpicture
\multiput{} at 0 0  8 3  /
\put{$(1,h(1))$}  [l] at 0 5.5
\put{with $h(1) = 4$} [l] at 0 4
\setdots <1mm>
\plot 0 0  8 0 /
\setdots <.5mm>
\plot 0 0  3 3  6 0  7 1 /
\plot 1 1  2 0  4 2 /
\plot 2 2  4 0  6 2  8 0 /
\put{$\ss 1$} at 0 -1
\put{$\ss 2$} at 2 -1
\put{$\ss 3$} at 4 -1
\put{$\ss 4$} at 6 -1
\put{$\ss 1$} at 8 -1
\setsolid
\setquadratic
\plot 0 0  1 1   2 2  3 2.5  4 2  5 1  6 0 /
\arr{5.9 0.1}{6 0}
\put{$\bullet$} at 0 0 
\setlinear
\setdashes <1mm>
\plot 0 -.3  0 2 /
\plot 8 -.3  8 2 /
\put{$e(1) = 0$}  at 4 -3
\endpicture} at -19 0

\put{\beginpicture
\multiput{} at 0 0  8 3  /
\put{$(2,h(2))$}  [l] at 0 5.5
\put{with $h(2) = 2$} [l] at 0 4
\setdots <1mm>
\plot 0 0  8 0 /
\setdots <.5mm>
\plot 0 0  3 3  6 0  7 1 /
\plot 1 1  2 0  4 2 /
\plot 2 2  4 0  6 2  8 0 /
\put{$\ss 1$} at 0 -1
\put{$\ss 2$} at 2 -1
\put{$\ss 3$} at 4 -1
\put{$\ss 4$} at 6 -1
\put{$\ss 1$} at 8 -1
\setsolid
\setquadratic
\plot 8 0  7 1  6.5 1.5  6 1.7  5.5 1.5  5.2 1.2  5 1  /
\plot  5 1  4 2   3.5 2.5 
           3  2.7  2.5 2.5  2 2 
  1.5 1.5  
           1.2 1   1.5 0.5  
           1.7 0.3  2 0 /
\arr{1.9 0.1}{2 0}
\setlinear
\setdashes <1mm>
\plot 0 -.3  0 2 /
\plot 8 -.3  8 2 /
\put{$e(2) = 3$}  at 5 -3
\endpicture} at -8 0

\put{\beginpicture
\multiput{} at 0 0  8 3  /
\put{} [l] at 0 4.7
\setdots <1mm>
\plot 0 0  8 0 /
\setdots <.5mm>
\plot 0 0  3 3  6 0  7 1 /
\plot 1 1  2 0  4 2 /
\plot 2 2  4 0  6 2  8 0 /
\put{$\ss 2$} at 2 -1
\put{$\ss 3$} at 4 -1
\put{$\ss 4$} at 6 -1
\put{$\ss 1$} at 8 -1

\setsolid
\setquadratic
\plot 2 0  1 0.7  0 0  /
\put{$\bullet$} at 2 0 
\setlinear
\setdashes <1mm>
\plot 8 -.3  8 2 /
\put{}  at 4 -3
\endpicture} at 0.05 -.55

\put{\beginpicture
\multiput{} at 0 0  8 3  /
\put{$(3,h(3))$}  [l] at 0 5.5
\put{with $h(3) = 1$} [l] at 0 4
\setdots <1mm>
\plot 0 0  8 0 /
\setdots <.5mm>
\plot 0 0  3 3  6 0  7 1 /
\plot 1 1  2 0  4 2 /
\plot 2 2  4 0  6 2  8 0 /
\put{$\ss 1$} at 0 -1
\put{$\ss 2$} at 2 -1
\put{$\ss 3$} at 4 -1
\put{$\ss 4$} at 6 -1
\put{$\ss 1$} at 8 -1
\setsolid
\setquadratic
\plot 4 0  5 1  5.5 1.5  6 1.7  6.5 1.5  7 1  8 0 /
\arr{7.9 0.1}{8 0}
\put{$\bullet$} at 4 0 
\setlinear
\setdashes <1mm>
\plot 0 -.3  0 2 /
\plot 8 -.3  8 2 /
\put{$e(3) = 0$}  at 4 -3
\endpicture} at 11 0

\put{\beginpicture
\multiput{} at 0 0  8 3  /
\put{$(4,h(4))$}  [l] at 0 5.5
\put{with $h(4) = 3$} [l] at 0 4
\setdots <1mm>
\plot 0 0  8 0 /
\setdots <.5mm>
\plot 0 0  3 3  6 0  7 1 /
\plot 1 1  2 0  4 2 /
\plot 2 2  4 0  6 2  8 0 /
\put{$\ss 1$} at 0 -1
\put{$\ss 2$} at 2 -1
\put{$\ss 3$} at 4 -1
\put{$\ss 4$} at 6 -1
\put{$\ss 1$} at 8 -1

\setsolid
\setquadratic
\plot 6 0  5 1       3.5 2.5  
           3  2.7    2.5 2.5      
           2.2 2   2.5 1.5  
           3 1  4 0 /
\arr{3.9 0.1}{4 0}
\put{$\bullet$} at 6 0 
\setlinear
\setdashes <1mm>
\plot 0 -.3  0 2 /
\plot 8 -.3  8 2 /
\put{$e(4) = 1$}  at 4 -3

\endpicture} at 22 0

\endpicture}
$$
The permutation $h$ is, in cycle notation, the permutation $(143)(2).$
	\medskip
{\bf Remarks.} (1) 
As we see, {\it the permutation $h_A$ may have fixed points.}
	\smallskip 
(2) {\it The permutation $h_A$ does not determine the algebra $A$.} There are many cyclic
Nakayama algebras say with two simple modules, but of course just two permuations of $\{1,2\}$.
For example, the algebras 
$$
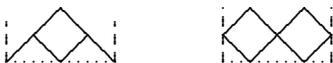
{\beginpicture
    \setcoordinatesystem units <.3cm,.3cm>
\put{\beginpicture
\multiput{} at 0 0  4 3 /
\setdots <1mm>
\plot 0 0  4 0 /
\setsolid
\plot 0 0  2 2  4 0 /
\plot 1 1  2 0  3 1 /
\setdashes <1mm>
\plot 0 0  0 1.5 /
\plot 4 0  4 1.5 /
\endpicture}  at 0 0
\put{\beginpicture
\multiput{} at 0 0  4 3 /
\setdots <1mm>
\plot 0 0  4 0 /
\setsolid
\plot 0 1  1 2  3 0  4 1  3 2  1 0  0 1 /
\setdashes <1mm>
\plot 0 0  0 2 /
\plot 4 0  4 2 /
\endpicture}  at 8 0
\endpicture}
$$
both have as homological permutation the identity permutation $h_A = (1)(2)$.
	\bigskip
{\bf 4.7. The sets $\Cal E(z)$.}
If $z$ is a natural number, then $\Cal E(z)$ is the set of simple modules
$S$ with $e(S) = z$ (thus, for any natural number $t$, the set $\Cal E(2t)$ 
consists of all simple modules $S$ with $2t = \pd IT \le \pd T$, and 
the set $\Cal E(2t+1)$ 
consists of all simple modules $S$ with $2t+1 = \pd T \le \pd IT$).
Of course, the sets $\Cal E(z)$ form a partion of $\Cal E$. 
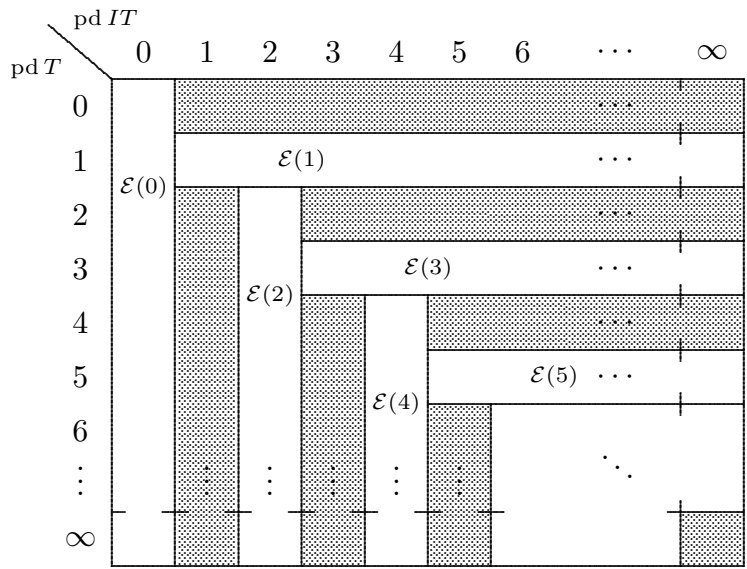
$$
{\beginpicture
    \setcoordinatesystem units <.7cm,.6cm>
\multiput{} at -2 1  12 -9.5 /
\plot 0 0  10  0  10 -9  9 -9 /
\plot 0 0  0 -9  9 -9  9 -7.8  /
\plot 8.8 -8  10 -8  /

\put{$\ddots$} at 8 -7 

\plot 1 0  1 -9 /
\plot 1 -1  10 -1 /
\plot 1 -2  10 -2 /

\plot 2 -2  2 -9 /
\plot 3 -2  3 -9 /

\plot 3 -3  10 -3 /
\plot 3 -4  10 -4 /

\plot 4 -4  4 -9 /
\plot 5 -4  5 -9 /

\plot 5 -5  10 -5 /
\plot 5 -6  10 -6 /

\plot 6 -6  6 -9 /

\multiput{$0$} at 0.5 0.5  -0.5 -0.5 /
\multiput{$1$} at 1.5 0.5  -0.5 -1.5 /
\multiput{$2$} at  2.5 0.5  -0.5 -2.5 /
\multiput{$3$} at  3.5 0.5  -0.5 -3.5 /
\multiput{$4$} at  4.5 0.5  -0.5 -4.5 /
\multiput{$5$} at  5.5 0.5  -0.5 -5.5 /
\multiput{$6$} at  6.5 0.5  -0.5 -6.5 /

\multiput{$\cdots$} at  8 0.5  8 -.5  8 -1.5  8 -2.5  8 -3.5  8 -4.5  8 -5.5 /
\multiput{$\vdots$} at  -0.5 -7.25 1.5 -7.25  2.5 -7.25  3.5 -7.25  4.5 -7.25  5.5 -7.25 /
\multiput{$\infty$} at  9.5 0.5  -0.5 -8.5 /
\setshadegrid span <.3mm>
\vshade 1 -1 0  10 -1 0    /
\vshade 3 -3 -2  10 -3 -2    /
\vshade 5 -5 -4  10 -5 -4    /

\vshade 1 -9 -2  2 -9 -2    /
\vshade 3 -9 -4  4 -9 -4    /
\vshade 5 -9 -6  6 -9 -6    /

\hshade -9  9 10   -8 9 10 /

\plot 9 0  9 -0.2 /
\plot 9 -.8  9 -1.2 /
\plot 9 -1.8  9 -2.2 /
\plot 9 -2.8  9 -3.2 /
\plot 9 -3.8  9 -4.2 /
\plot 9 -4.8  9 -5.2 /
\plot 9 -5.8  9 -6.2 /

\plot 0 -8  0.2 -8 /
\plot 0.8 -8  1.2 -8 /
\plot 1.8 -8  2.2 -8 /
\plot 2.8 -8  3.2 -8 /
\plot 3.8 -8  4.2 -8 /
\plot 4.8 -8  5.2 -8 /
\plot 5.8 -8  6.2 -8 /

\plot 0 0  -1 1 /
\put{$\ss\pd T$} at -1.2 .2
\put{$\ss\pd IT$} at -.1 1.1

\put{$\ss \Cal E(0)$} at .5 -2 
\put{$\ss \Cal E(2)$} at 2.5 -4
\put{$\ss \Cal E(4)$} at 4.5 -6 
\put{$\ss \Cal E(1)$} at 3 -1.5 
\put{$\ss \Cal E(3)$} at 5 -3.5 
\put{$\ss \Cal E(5)$} at 7 -5.5 

\endpicture}
$$

In case $A$ is connected and has infinite global dimension, then, for $T$ simple, 
$\pd T$ cannot be even,
$\pd IT$ cannot be odd, thus it is sufficient to look at the following display:
$$
{\beginpicture
    \setcoordinatesystem units <.8cm,.7cm>
\multiput{} at -2 1  9 -6.5 /
\plot 0 0  7  0  7 -6  6 -6 /
\plot 0 0  0 -6  6 -6  6 -4.8 /
\plot 5.8 -5  7 -5  /
\plot 1 0  1 -6 /
\plot 1 -1  7 -1 /
\plot 2 -1  2 -6 /
\plot 2 -2  7 -2 /

\plot 3 -2  3 -6 /

\plot 3 -3  7 -3 /

\multiput{$0$} at 0.5 0.5  /
\multiput{$1$} at -0.5 -0.5 /
\multiput{$2$} at  1.5 0.5  /
\multiput{$3$} at -0.5 -1.5 /
\multiput{$4$} at  2.5 0.5  /
\multiput{$5$} at -0.5 -2.5 /
\multiput{$6$} at  3.5 0.5   /
\multiput{$7$} at -0.5 -3.5 /

\multiput{$\cdots$} at  5 0.5  5 -.5  5 -1.5  5 -2.5  /
\multiput{$\vdots$} at  -0.5 -4.2  .5 -4.2  1.5 -4.2  2.5 -4.2  /
\multiput{$\infty$} at  6.5 0.5  -0.5 -5.5 /

\plot 6 0  6 -0.2 /
\plot 6 -.8  6 -1.2 /
\plot 6 -1.8  6 -2.2 /
\plot 6 -2.8  6 -3.2 /

\plot 0 -5  0.2 -5 /
\plot 0.8 -5  1.2 -5 /
\plot 1.8 -5  2.2 -5 /
\plot 2.8 -5  3.2 -5 /

\plot 0 0  -1 1 /
\put{$\ss\pd T$} at -1.2 .2
\put{$\ss\pd IT$} at -.1 1.1

\put{$\ss \Cal E(0)$} at .5 -1 
\put{$\ss \Cal E(2)$} at 1.5 -2
\put{$\ss \Cal E(4)$} at 2.5 -3 
\put{$\ss \Cal E(1)$} at 2 -.5 
\put{$\ss \Cal E(3)$} at 3 -1.5 
\put{$\ss \Cal E(5)$} at 4 -2.5

\setshadegrid span <.3mm>
\hshade -6  6 7   -5 6 7 /

\put{$\ddots$} at 5 -4 

\endpicture}
$$
	\medskip
Let us return to the general case where $A$ has arbitrary global dimension.
We show: {\it all cases $(p,q) = (\pd S, \pd IS)$ which are not shaded, actually can occur.}
	\medskip
Proof. 
	\smallskip
Of course, the case $(0,0)$ occurs only if $A$ has a connected component with linear quiver
(the number of connected components with linear quiver is just the number
of simple modules $S$ with $\pd S = 0$). On the other hand, 
the cases $(z,\infty)$ and $(\infty,z)$ 
can occur only in case the global dimension of $A$ is infinite.

In the cases $(z,\infty)$ and $(\infty,z)$, the examples provided are cyclic Nakayama 
algebras (of course). If $(p,q)\neq (0,0)$ and $p\neq \infty,$ and $q\neq \infty$, 
we provide translation quivers which are the Auslander-Reiten quivers of a non-cyclic
Nakayama algebra, and which also can be part of the Auslander-Reiten quiver of a
cyclic Nakayama algebra.
	\medskip
{\bf The case $(0,0)$.}
If $A$ is a Nakayama algebra with a simple projective module
$S$, then $\pd S = 0,\ \pd IS = 0.$ 
	\medskip
{\bf The case $(t,0)$ with $t\ge 1$.}
$$
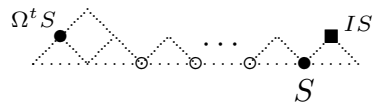
{\beginpicture
    \setcoordinatesystem units <.3cm,.3cm>
\multiput{} at -2 0  10  2  /
\setdots <1mm>
\plot -2 0  10 0 /
\setdots <.5mm>
\plot -2 0  0 2  2 0  3 1  4 0 /
\plot -1 1  0 0  1 1 /
\plot 6 0  7 1  8 0 9 1  10 0 /
\put{$\cdots$} at 5 .7
\put{$\sssize \blacksquare$} at 9 1 
\multiput{$\bullet$} at -1 1  8 0 /
\put{$S$}$ at 8 -1 
\put{$\ss IS$}$ at 10 1.5
\multiput{$\circ$} at 6 0  4 0  2 0  /
\put{$\ss\Omega^t S$} at -2 1.7
\endpicture}
$$
	\medskip
{\bf The case $(\infty,0)$.}
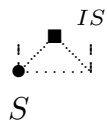
$$
{\beginpicture
    \setcoordinatesystem units <.4cm,.4cm>
\multiput{} at 0 0  2  1  /
\setdots <1mm>
\plot 0 0  2 0 /
\setdots <.5mm>
\plot 0 0  1 1  2 0 /
\put{$\sssize \blacksquare$} at 1 1
\put{$\bullet$} at 0 0
\put{$S$}$ at 0 -1 
\put{$\ss IS$}$ at 2 1.5
\setdashes <1mm>
\plot 0 0  0 1 /
\plot 2 0  2 1 /
\endpicture}
$$
	\medskip
{\bf The case $(1,1).$}
$$
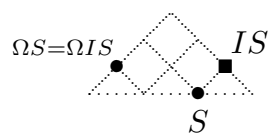
{\beginpicture
    \setcoordinatesystem units <.3cm,.3cm>
\multiput{} at 0 0  6  3  /
\setdots <1mm>
\plot 0 0  6 0 /
\setdots <.5mm>
\plot 0 0  3 3  6 0 /
\plot 1 1  2 0  4 2 /
\plot 2 2  4 0  5 1 /
\put{$\sssize \blacksquare$} at 5 1
\put{$\bullet$} at 4 0
\put{$\dsize \bullet$} at 1 1
\put{$S$}$ at 4 -1 
\put{$IS$}$ at 5.9 1.9
\put{$\ss\Omega S = \Omega IS$} at -1 1.7
\endpicture}
$$
	\medskip
{\bf The case $(z,z)$ with $z\ge 2.$}
$$
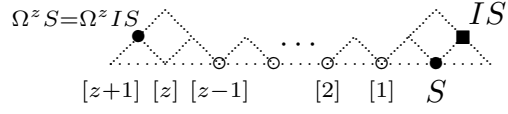
{\beginpicture
    \setcoordinatesystem units <.3cm,.3cm>
\multiput{} at 0 0  14  2  /
\setdots <1mm>
\plot 0 0  14 0 /
\setdots <.5mm>
\plot 0 0  2 2  4 0  5 1  6 0 /
\plot 8 0  9 1  10 0  12 2  14 0 /
\plot 1 1  2 0  3 1 /
\plot 11 1  12 0  13 1 /
\put{$\sssize \blacksquare$} at 13 1
\put{$\bullet$} at 12 0 
\put{$\dsize \bullet$} at 1 1
\put{$S$}$ at 12 -1 
\put{$IS$}$ at 13.9 1.9
\put{$\ss\Omega^z S = \Omega^z IS$} at -1.3 1.7
\put{$\ss [1]$} at 10 -1
\put{$\ss [2]$} at 8 -1
\put{$\ss [z-1]$} at 4 -1
\put{$\ss [z]$} at 2 -1
\put{$\ss [z+1]$} at 0 -1
\multiput{$\circ$} at 10 0  8 0  6 0  4 0  /
\put{$\cdots$} at 7 .7
\endpicture}
$$
In this example, we have $\Omega^i S = \Omega^i IS = [i]$, for $1\le i \le z\!-\!1,$ and
$\Omega^zS = \Omega^zIS = P[z].$
	\medskip
If $z = \pd S = \pd IS$ is odd, then always $\Omega^z S = \Omega^z IS$. However, if
$z = \pd S = \pd IS$ is even,
then we may have that $\Omega^z S \neq \Omega^z IS$,
as the following example shows (with $z = 2t$):
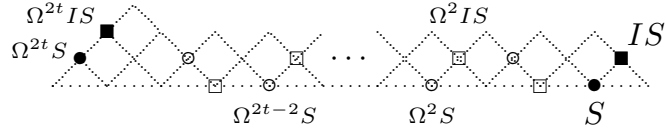
$$
{\beginpicture
    \setcoordinatesystem units <.3cm,.3cm>
\multiput{} at 8 0  30  2  /
\setdots <1mm>
\plot 8 0  30 0 /
\setdots <.5mm>
\plot 8 0  11 3  12 2 /
\plot 9 1  10 0  12 2  14 0  16 2  18 0 /
\plot 10 2  12 0  14 2  16 0  18 2 /
\plot 20 0  22 2  24 0  26 2  28 0  29 1 /
\plot 20 2  22 0  24 2  26 0  28 2  30 0 /

\multiput{$\sssize \blacksquare$} at 29 1  10 2 /
\multiput{$\bullet$} at 28 0  9 1 /
\multiput{$\sssize \square$} at   14 0  17 1  26 0  23 1  /
\put{$S$}$ at 28 -1 
\put{$IS$}$ at 29.9 1.9
\put{$\ss \Omega^{2}S$} at 22 -1
\put{$\ss \Omega^{2t-2}S$} at 16.2 -1
\multiput{$\circ$} at   13 1  16 0  25 1  22 0   /
\put{$\cdots$} at 19 1 
\put{$\ss \Omega^{2}IS$} at 23 2.7

\put{$\ss \Omega^{2t}S$} at 7.5 1.4
\put{$\ss \Omega^{2t}IS$} at 8.3 2.7
\endpicture}
$$
	\medskip 
{\bf The case $(2t+u,2t)$ with $u \ge 1$ and $t\ge 1.$ Let $v = 2t+u$.}
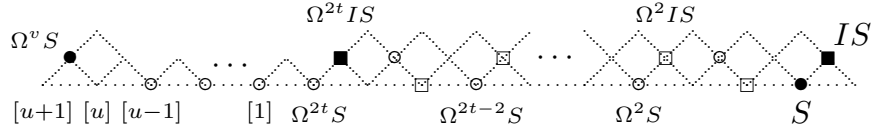
$$
{\beginpicture
    \setcoordinatesystem units <.3cm,.3cm>
\multiput{} at 0 0  30  2  /
\setdots <1mm>
\plot 0 0  30 0 /
\setdots <.5mm>
\plot 0 0  2 2  4 0  5 1  6 0 /
\plot 1 1  2 0  3 1 /
\plot 8 0  9 1  10 0  12 2  14 0  16 2  18 0 /
\plot 11 1  12 0  14 2  16 0  18 2 /
\plot 20 0  22 2  24 0  26 2  28 0  29 1 /
\plot 20 2  22 0  24 2  26 0  28 2  30 0 /

\multiput{$\sssize \blacksquare$} at 29 1  11 1 /
\put{$\bullet$} at 28 0 
\put{$\dsize \bullet$} at 1 1
\multiput{$\sssize \square$} at 11 1  14 0  17 1  26 0  23 1  /
\put{$S$}$ at 28 -1 
\put{$IS$}$ at 29.9 1.9
\put{$\ss\Omega^v S$} at -.3 1.7
\put{$\ss \Omega^{2}S$} at 22 -1
\put{$\ss \Omega^{2t-2}S$} at 16.2 -1
\put{$\ss \Omega^{2t}S$} at 10.2 -1
\put{$\ss [1]$} at 8 -1
\put{$\ss [u-1]$} at 4 -1
\put{$\ss [u]$} at 2 -1
\put{$\ss [u+1]$} at 0 -1
\multiput{$\circ$} at 10 0  8 0  6 0  4 0  13 1  16 0  25 1  22 0   /
\put{$\cdots$} at 19 1 
\put{$\ss \Omega^{2}IS$} at 23 2.7
\put{$\ss \Omega^{2t}IS$} at 11 2.7
\put{$\cdots$} at 7 .7
\endpicture}
$$
where, on the right, starting with $\Omega^{2t}S$ and ending with $\tau^-S$, there
are $3t+2$ simple modules (thus altogether there are $3t+u+3$ simple modules).
Note that $\Omega^{2t}S$ is projective, thus $\pd S = 2t$. Of course, also
$\Omega^vS = \Omega^{2t+u}S = \Omega^u(\Omega^{2t}S)$ is projective, thus $\pd IS = v = 2t+u.$
	\medskip
{\bf The case $(\infty,2t)$ with $t\ge 1.$}
$$
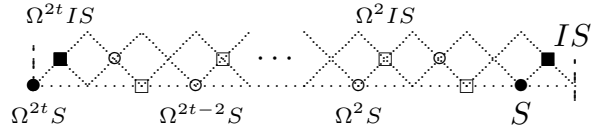
{\beginpicture
    \setcoordinatesystem units <.3cm,.3cm>
\multiput{} at 10 0  30  2  /
\setdots <1mm>
\plot 10 0  30 0 /
\setdots <.5mm>
\plot 10 0  12 2  14 0  16 2  18 0 /
\plot 11 1  12 0  14 2  16 0  18 2 /
\plot 20 0  22 2  24 0  26 2  28 0  29 1 /
\plot 20 2  22 0  24 2  26 0  28 2  30 0 /

\multiput{$\sssize \blacksquare$} at 29 1  11 1 /
\multiput{$\bullet$} at 28 0  10 0 /
\multiput{$\sssize \square$} at 11 1  14 0  17 1  26 0  23 1  /
\put{$S$}$ at 28 -1 
\put{$IS$}$ at 29.9 1.9
\put{$\ss \Omega^{2}S$} at 22 -1
\put{$\ss \Omega^{2t-2}S$} at 16.2 -1
\put{$\ss \Omega^{2t}S$} at 10.2 -1
\multiput{$\circ$} at 10 0   13 1  16 0  25 1  22 0   /
\put{$\cdots$} at 19 1 
\put{$\ss \Omega^{2}IS$} at 23 2.7
\put{$\ss \Omega^{2t}IS$} at 11 2.7
\setdashes <1mm>
\plot 10 0  10 1.5 /
\plot 30 -.5  30 1 /
\endpicture}
$$
It is easy to check that $\Omega^{2t}S$ is $\Omega$-periodic. Thus $\pd S = \infty$.
	\medskip
{\bf The case $(2t+1,2t+2)$ with $t\ge 0.$}  
$$
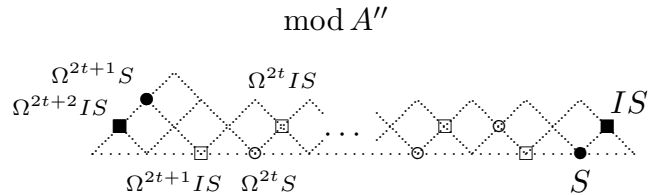
{\beginpicture
    \setcoordinatesystem units <.3cm,.3cm>
\multiput{} at -2 0  18 0  /
\setdots <1mm>
\plot -2 0  18 0 /
\setdots <.5mm>
\plot -1 1  0 0  2 2 /
\plot -2 0 1 3  4 0  6 2  6.5 1.5 /
\plot 0 2  2 0  4 2  6 0  6.5 0.5 /
\plot 8.5  0.5  10 2  12 0  14 2  16 0  17 1 /
\plot 8.5  1.5  10 0  12 2  14 0  16 2  18 0 /

\put{$\cdots$} at 7.3 .7

\multiput{$\bullet$} at 16 0  0 2 /
\multiput{$\circ$} at  0 2  4 0  13 1  10 0  /
\multiput{$\sssize \square$} at -1 1 14 0  11 1  2 0  5 1 /
\multiput{$\sssize \blacksquare$} at 17 1   -1 1 /
\put{$\mod A''$} at 7 5
\put{$S$} at 16 -1 
\put{$IS$} at 17.8 1.9 
\put{$\ss \Omega^{2t+1}IS$} at 1 -1  
\put{$\ss \Omega^{2t}S$} at 4.5 -1 
\put{$\ss \Omega^{2t+1}S$} at -2 3 
\put{$\ss \Omega^{2t+2}IS$} at -3.2 1.8 
\put{$\ss \Omega^{2t}IS$} at 5 2.8 
\endpicture}
$$
	\medskip
{\bf The case $(2t+1,2t+2+u)$ with $u \ge 1$ and $t\ge 1.$ Let $v = 2t+2+u$.}
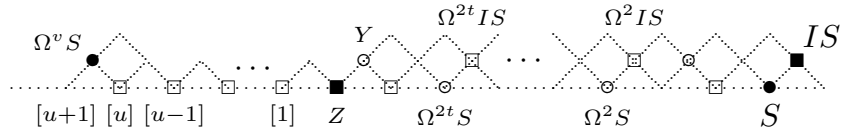
$$
{\beginpicture
    \setcoordinatesystem units <.3cm,.3cm>
\multiput{} at 0 0  30  2  /
\setdots <1mm>
\plot 0 0  30 0 /
\setdots <.5mm>
\plot 2 0  4 2  6 0  7 1  8 0 /
\plot 3 1  4 0  5 1 /
\plot   10 0  11 1 /
\plot 13 1  14 0  16 2  18 0 /
\plot 11 1  12 0  14 2  16 0  18 2 /
\plot 20 0  22 2  24 0  26 2  28 0  29 1 /
\plot 20 2  22 0  24 2  26 0  28 2  30 0 /
\put{$\cdots$} at 9 .7

\multiput{$\sssize \blacksquare$} at 29 1  12 0 /
\multiput{$\bullet$} at 28 0  3 1 /
\multiput{$\sssize \square$} at 12 0  14 0  17 1  26 0  23 1  10 0  8 0  6 0  4 0  /
\put{$S$}$ at 28 -1 
\put{$IS$}$ at 29.9 1.9
\put{$\ss\Omega^v S$} at 1.7 1.7
\put{$\ss Z$} at 12 -1
\put{$\ss [1]$} at 10 -1
\put{$\ss [u-1]$} at 6 -1
\put{$\ss [u]$} at 4 -1
\put{$\ss [u+1]$} at 2 -1
\multiput{$\circ$} at 13 1  16 0  25 1  22 0   /
\put{$\cdots$} at 19 1 
\put{$\ss \Omega^{2}S$} at 22 -1
\put{$\ss \Omega^{2}IS$} at 23 2.7

\put{$\ss \Omega^{2t}S$} at 16 -1
\put{$\ss \Omega^{2t}IS$} at 17 2.7
\put{$\ss Y$} at 13 2
\endpicture}
$$
with $Y = \Omega^{2t+1}S$ and $Z = \Omega^{2t+2}IS$.
On the right, starting with $Z = \Omega^{2t+2}S$ and ending with $\tau^-S$, there
are $3t+5$ simple modules (thus altogether there are $3t+u+6$ simple modules).
Note that $Y = \Omega^{2t+1}S$ is projective, thus $\pd S = 2t+1$. Of course, also
$\Omega^vS = \Omega^{2t+2+u}S = \Omega^u(\Omega^{2t+2}S) = \Omega^2Z$
is projective, thus $\pd IS = v = 2t+2+u.$
	\medskip
{\bf The case $(2t+1,\infty)$ with $t\ge 0.$}
$$
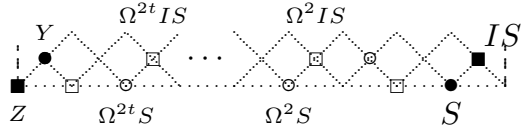
{\beginpicture
    \setcoordinatesystem units <.3cm,.3cm>
\multiput{} at 12 0  30  2  /
\setdots <1mm>
\plot 12 0  30 0 /
\setdots <.5mm>
\plot 13 1  14 0  16 2  18 0 /
\plot 12 0  14 2  16 0  18 2 /
\plot 20 0  22 2  24 0  26 2  28 0  29 1 /
\plot 20 2  22 0  24 2  26 0  28 2  30 0 /

\multiput{$\sssize \blacksquare$} at 29 1  12 0  /
\multiput{$\bullet$} at 28 0  13 1 /
\multiput{$\sssize \square$} at 12 0  14 0  17 1  26 0  23 1  /
\put{$S$}$ at 28 -1 
\put{$IS$}$ at 29.9 1.9
\put{$\ss Z$} at 12 -1
\multiput{$\circ$} at 13 1  16 0  25 1  22 0   /
\put{$\cdots$} at 19 1 
\put{$\ss \Omega^{2}S$} at 22 -1
\put{$\ss \Omega^{2}IS$} at 23 2.7

\put{$\ss \Omega^{2t}S$} at 16 -1
\put{$\ss \Omega^{2t}IS$} at 17 2.7
\put{$\ss Y$} at 13 2

\setdashes <1mm>  
\plot 12 0  12 1.5 /
\plot 30 0  30 1.5 /
\endpicture}
$$
Again, $Y = \Omega^{2t+1}S$ and $Z = \Omega^{2t+2}IS$.
What we have to verify (but this is easy), is that $Z$ is $\Omega$-periodic.
	
$\s$
	\medskip
{\bf Remark.} Let $A$ be a cyclic Nakayama algebra. Then $a(A)$ has been defined to be the
maximum of the numbers $a(S)$, where $S$ is simple and not $\psi$-cyclic. 
	\medskip 
{\it The set $\Cal E(2a(A))$ consists of the $\psi$-cyclic simple modules $S$ with $\pd IS = 2a(A)$.
Similarly, $\Cal E^*(2a(A))$ consists of the $\gamma$-cyclic simple modules $S$ 
with $\id PS = 2a(A)$.}
	\medskip
Proof. Let $a = a(A).$ As we have seen in 3.10, we have $\finpro A \le 2a.$ 
Assume that $S$ is a simple module in $\Cal E(2a).$ Then $\pd S \ge 2a$ and $\pd IS = 2a$.
Since $\finpro A \le 2a,$ it follows from $\pd S \ge 2a$ that $S$ cannot be odd, thus $S$ is 
$\psi$-cyclic. Conversely, assume that $S$ is $\psi$-cyclic and $\pd IS = 2a.$ Then
$\pd S$ is even or infinite. But if $\pd S$ is even, then $\pd S \ge \pd IS$, see 3.2 (b).
We see that always $\pd S \ge 2a$, therefore $S\in \Cal E(2a).$ $\s$
	\bigskip
{\bf 4.8. Ties for small values of $z$.}
We are going to look at small values of $z$, 
and describe the bijection $h_z\:\Cal E(z) \to \Cal E^*(z).$
	\medskip
{\bf The case $z = 0.$}
	\medskip
{\it Let $S, T$ be simple modules. Then $T$ belongs to $\Cal E(0)$ iff $T$ is torsionless, iff $IT$ is projective. Also, $S$ belongs to $\Cal E^*(0)$ iff $S$ is 
divisible, iff $PS$ is projective.} 
	\medskip
{\it The bijection $h_0\:\Cal E(0) \to \Cal E^*(0)$ sends a simple torsionless module $T$
to $\top IT$, the inverse sends a simple divisible module to $\soc PS.$}
	
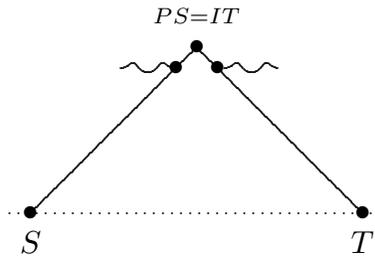
\medskip
$$ 
{\beginpicture
    \setcoordinatesystem units <.23cm,.23cm>
\multiput{} at  0 0  16 8  /
\multiput{$\bullet$} at 0 0  7 7  8 8  9 7  16 0 /
\plot 0 0  8 8  16 0 /
\setdots <1mm>
\plot -1 0  17 0 /
\setsolid 
\multiput{$\snake$} at 5 7  10 7 /
\put{$S$} at 0 -1.4
\put{$T$} at 16 -1.4
\put{$\ss PS = IT$} at 8 9.5
\endpicture}
$$

	\bigskip
{\bf The case $z = 1$.}
	\medskip
{\it Let $S, T$ be simple modules. Then $T$ belongs to $\Cal E(1)$ iff $T$ is not
torsionless, and $\pd T = 1$. And  $S$ belongs to $\Cal E^*(1)$ iff $S$ is not
divisible, and $\id S = 1$.} 
	\medskip
{\it The bijection $h_1\:\Cal E(1) \to \Cal E^*(1)$ is given by the Auslander-Reiten translation
$\tau$.}
	
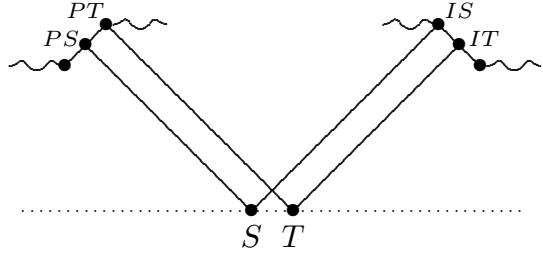
\medskip
$$ 
{\beginpicture
    \setcoordinatesystem units <.23cm,.23cm>
\setdots <1mm>
\plot -12 0  12 0 /
\setsolid
\multiput{} at -11 0   11 11 /
\plot -10 7  -8 9  1 0  9 8 /
\plot -9 8  -1 0  8 9  10  7 / 
\put{$S$} at -1 -1.2 
\put{$T$} at  1 -1.2
\put{$\ss PS$} at -10.2 8.3  
\put{$\ss PT$} at -9 9.7  

\put{$\ss IS$} at 9 9.7  
\put{$\ss IT$} at 10.2 8.3  

\multiput{$\snake$} at -7 9  6 9  -12 7  11 7  /
\multiput{$\bullet$} at -10 7  -9 8  -8 9  10 7  9 8  8 9   -1 0  1 0 /
\endpicture}
$$
	\bigskip
{\bf Remark. $1$-free algebras.}
We say that a Nakayama algebra is {\it $1$-free} provided $\Cal E(1)$ is empty.
Given a Nakayama algebra $A$, we may delete all indecomposable
modules $M$ with $\soc M$ or $\top M$ in $\Cal E(1)$ and obtain the module category of a $1$-free
algebra $A'$ which has similar properties as $A$. 
	\bigskip
{\bf The case $z = 2$.}
	\medskip
{\it Let $S, T$ be simple modules. Then $T$ belongs to $\Cal E(2)$ iff $\pd IT = 2$
and $\pd T \ge 2$. And, $S$ belongs to $\Cal E^*(2)$ iff $\id PS = 2$ and $\id S \ge 2.$}
	\medskip
{\it The bijection $h_2\:\Cal E(2) \to \Cal E^*(2)$ sends $T$ to $\gamma\phi(T) = \tau \soc PIT$.}
(Recall that $\phi(T) = \top IT$ and $\gamma(T') = \tau\soc P(T').$)
	\medskip
We have the following diagram. 
$$ 
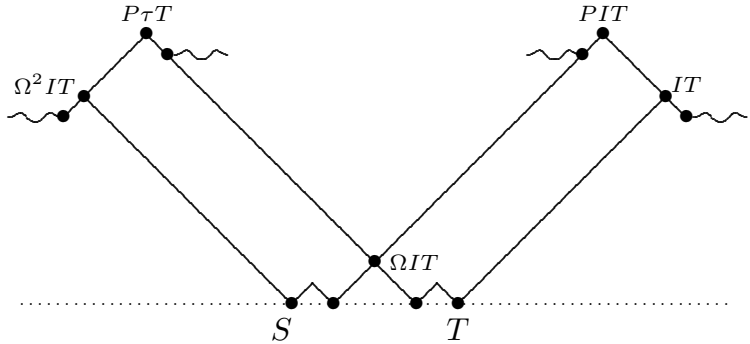
{\beginpicture
    \setcoordinatesystem units <.23cm,.23cm>
\multiput{} at -23 0   15 5 /
\multiput{$\bullet$} at  
    -15 9  -11 13  0 2  11 13  15 9   14 10   -14 10 
    10 12  -10 12  /
\plot 
    -15 9  -11 13  0 2  11 13  15 9 /
\plot  0 2  2 0  3 1 /
\plot  0 2  -2 0  -3 1 /

\plot 14 10  4 0  3 1 /
\plot -14 10  -4 0  -3 1 /
\multiput{$\snake$} at -17 9  -9 12   8 12  16 9 /

\setdots <1mm>
\plot -17 0  17 0 /

\put{$T$} at 4 -1.2
\put{$S$} at -4.5 -1.2

\put{$\ss IT$} at 15 10.6
\put{$\ss \Omega^2IT$} at -15.9 10.6
\put{$\ss PIT$} at 11 14
\put{$\ss P\tau T$} at -11 14
\put{$\ss \Omega IT$} at 1.9 2

\multiput{$\bullet$} at 2 0  4 0  -2 0  -4 0 /

\endpicture}
$$
Note that $S = \tau^{a+1}T$, where $a = |\Omega IT|$.
The modules $\tau^iT$ with $1\le i < a$ belong to $\Cal E(1)$. 
If $A$ is $1$-free,
then $\Omega IT = \tau T,$ so that $a = 1$  and therefore $h(T) = S = \tau^2 T$ and $\tau T$
is both torsionless and divisible. 
	\medskip
{\bf 4.9. Remark.} As we have mentioned, $\Cal E(0)$ consists of the torsionless simple modules,
thus $|\Cal E(0)| = r$, where $r$ is the number of peaks of $A$.
We also have $|\bigcup_{z\ge 2}\Cal E(z)| \le r$. Namely, if $S$ is simple and 
$e(S)\ge 2$, then 
$\tau S$ is divisible, thus $\tau$ maps $\bigcup_{z\ge 2}\Cal E(z)$ into the set
$\Cal E^*(0)$. Since $\tau$ (restricted to $\bigcup_{z\ge 2}\Cal E(z)$) is injective, 
and $|\Cal E^*(0)| = r$, the assertion follows. 
	\bigskip
{\bf 4.10. Remark.} We have seen in 3.7 (4) and (5) that for calculating $\finpro A$,
we may look at the functions $e$ or $f$. Whereas $f$ is easier to handle, 
the importance of $e$ lies in the symmetry assertion formulated in Corollary 4.2.
There is no corresponding assertion for $f$, as we want to show now. Thus, let
$f^*(S) = 0$ in case $\id S$ is odd and $\id PS$ is even, $f^*(S) = \id S$ in case
$\id PS$ is not even, and $f^*(S) = \id PS$ in case $\id S$ is not odd. 
In general, {\it there is no permutation $\pi$ 
of the simple modules with $f^*(\pi S) = f(S)$ for all simple modules $S$,} see the following 
example. 
	\medskip
{\bf Example.} 
$$
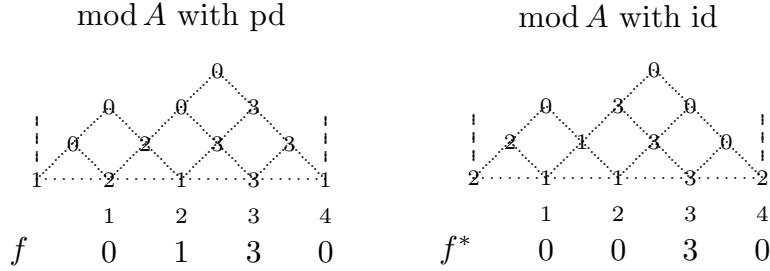
{\beginpicture
    \setcoordinatesystem units <.4cm,.4cm>

\put{\beginpicture
\multiput{} at 0 -1  8 3.5 /
\setdots <1mm>
\plot 0 0  8 0 /
\setdots <.5mm>
\plot  0 0  2 2  4 0  6 2 /
\plot  1 1  2 0  5 3  8 0 /
\plot 4 2  6 0  7 1 /
\multiput{$\ss 0$} at 1 1  2 2  4 2  5 3  /
\multiput{$\ss 1$} at  0 0  4 0  8 0 /
\multiput{$\ss 2$} at 2 0   3 1    /
\multiput{$\ss 3$} at 5 1  6 0  6 2  7 1    /
\setdashes <1mm>
\plot 0 0.5  0 2 /
\plot 8 .5  8 2 /
\put{$\ss 4$} at 8 -1
\put{$\ss 3$} at 6 -1
\put{$\ss 2$} at 4 -1
\put{$\ss 1$} at 2 -1
\put{$\mod A$ with $\pd$} at 4 4.5
\put{$f$} at -.5 -2
\put{$0$} at 8 -2
\put{$3$} at 6 -2
\put{$1$} at 4 -2
\put{$0$} at 2 -2

\endpicture} at 0 0 
\put{\beginpicture
\multiput{} at 0 -1  8 3.5 /
\setdots <1mm>
\plot 0 0  8 0 /
\setdots <.5mm>
\plot  0 0  2 2  4 0  6 2 /
\plot  1 1  2 0  5 3  8 0 /
\plot 4 2  6 0  7 1 /
\multiput{$\ss 0$} at 5 3  6 2  7 1  2 2  /
\multiput{$\ss 1$} at 2 0  3 1  4 0   /
\multiput{$\ss 2$} at 0 0  1 1  8 0 /
\multiput{$\ss 3$} at 4 2  5 1  6 0  /
\setdashes <1mm>
\plot 0 0.5  0 2 /
\plot 8 .5  8 2 /
\put{$\ss 4$} at 8 -1
\put{$\ss 3$} at 6 -1
\put{$\ss 2$} at 4 -1
\put{$\ss 1$} at 2 -1
\put{$\mod A$ with $\id$} at 4 4.5

\put{$f^*$} at -.5 -2
\put{$0$} at 8 -2
\put{$3$} at 6 -2
\put{$0$} at 4 -2
\put{$0$} at 2 -2

\endpicture} at 12 0 
\endpicture}
$$
On the left, we have equipped the Auslander-Reiten quiver with the values $\pd M$;
on the right, with the values $\id M$. We label the simple modules by $1$ up to $4$
as shown below the Auslander-Reiten quivers and provide
the values $f(S)$ or $f^*(S)$, respectively, below the label of $S$. We see that
there is a simple module $S$, namely the one with index 2, such that $f(S) = 1.$
But there is no simple module $S$ with $f^*(S) = 1.$ 
	\medskip

One may also ask, whether there is a permutation $\pi$ of $\Cal E$ such that $\des \pi(S) =
\del S$ for any simple module $S$. Of course, this is true in case $\del S = e(S)$
for all simple modules $S$ (but, as we have mentioned, this is not known, see 3.8). 
	\bigskip\bigskip

{\bf 5. Proof of Theorem 4.1.}
	\medskip
Let $A$ be a connected Nakayama algebra. 
First, we consider the case that $A$ is linear. In case $A$ is cyclic, we will use
covering theory in order to reduce it to the linear case.
	\medskip
{\bf 5.1. The case of $A$ being a linear Nakayama algebra.}
Let $A$ be a linear Nakayama algebra. The proof of the assertions (a) and (b) 
will be quite easy, 
as soon as we have introduced suitable labels for the modules involved.
	\medskip
Dealing with a linear Nakayama algebra $A$, we identify the set of simple modules
with a finite interval of  $\Bbb Z$, say with $[\alpha,\omega] \subset \Bbb Z$, with
$\alpha \le \omega$, 
such that the Auslander-Reiten translation $\tau$ is given by $\tau x = x-1$. 
By abuse of notation, we also write $\tau (\omega+1) = \omega$. In this way,
the indecomposable modules can be written in the form $[x,y[$, where $\alpha \le x < y \le
\omega+1$. We have $\soc [x,y[\ = x$ and $\top [x,y[\ = \tau y$. This convention  
has the following advantage: if $[x,y[$ is non-projective, then
$\Omega[x,y[\ = [\soc P\tau y,x[$. 
Similarly, if  $[x,y[$ is non-injective, then
$\Sigma[x,y[\ = [y, \tau^-\top Ix[$. 
This fact should be considered as
the start of the general setting which we will use.
	\medskip
(1) For a pair $V \subseteq W$ of indecomposable non-projective modules with the same
socle, let $P_0$ and
$Q_0$ be the projective covers and $\Omega V,\ \Omega W$ the first syzygy modules, of
$V$ and $W$, respectively. 
Let $\soc V = \soc W = u_0$ and $\top V = \tau v_0,$ $\top W = \tau w_0$, thus
$V = [u_0,v_0[$ and $W = [u_0,w_0[.$ 
Let $P_0 = [v_1,v_0[$ and $Q_0 = [w_1,w_0[$, then $v_1 \le w_1 < u_0$ (the strict
inequality is due to the fact that we assume that $W$ is not projective),
and $\Omega V = [v_1,u_0[,\ 
\Omega W = [w_1,u_0[.$

For an indecomposable non-projective module $V'$ with a non-zero factor module $W'$, let $B_1$
be a projective cover of $V'$ and $\Omega V' \subseteq \Omega W' \subset B_1$ the syzygy modules
of $V'$ and $W'$ respectively. Let $\soc V' = v_1,$
 $\soc W' = w_1$ and $\top V' = \top W' = \tau u_0.$ The projective cover of $V'$ ad $W'$ is
$B' = [u_1,u_0[$ for some $u_1 < v_1$ and we have 
$\Omega V' = [u_0,v_1[$ and $\Omega W' = [u_1,w_1[.$
	\smallskip
Thus, there are the following two pictures: on the right, we start with the pair
$V\subseteq W$, on the left, we start with $V'$ and its factor module $W'$.
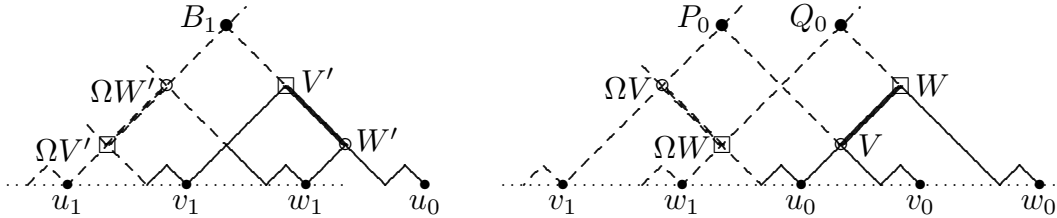
$$
{\beginpicture
    \setcoordinatesystem units <.22cm,.22cm>
\put{\beginpicture
\multiput{} at 0 0  16 10 /
\setdots <1mm>
\plot -1 0  16 0 /
\put{$u_0$} at 20 -1
\put{$w_1$} at 14 -1
\put{$v_1$} at  8 -1
\put{$u_1$} at 2 -1 
\setsolid
\plot 6 0  7 1  8 0  13 5  16 2 /
\plot 10 2  12 0  13 1  14 0  16 2  18 0  19 1 20 0 /

\setdashes <1mm>
\plot 10 8  13 5 /
\plot 0 0  1 1  2 0  11 9 /
\plot 3 3  6 0 /
\plot 6 6  10 2 /

\plot 4.1 2  7.1 5 /
\plot 3.9 2  6.9 5 /

\put{$W'$} at 17.6 2.3
\put{$V'$} at 14.6 5.3
\put{$\Omega W'$} at 4.8 4.8
\put{$\Omega V'$} at 1.8 1.8
\put{$B_1$} at 8.6 8.4
\multiput{$\circ$} at 7 5  16 2 / 
\multiput{$\ss \square$} at 4 2  13 5 /
\put{$\bullet$} at 10 8
\multiput{$\ss \bullet$} at 2 0  8 0   14 0  20 0 /

\setsolid
\plot 13.1 5  16.1 2 /
\plot 12.9 5  15.9 2 /

\endpicture} at 0 0 
\put{\beginpicture

\plot 16.1 2  19.1 5 /
\plot 15.9 2  18.9 5 /

\multiput{} at 0 0  26 10 /
\setdots <1mm>
\plot -1 0  27 0 /
\setsolid
\plot 12 0  13 1  14 0  19 5  24 0  25 1  26 0 /
\plot 16 2  18 0  19 1  20 0 /
\setdashes <1mm>
\plot 16 8  19 5 /
\plot 16 2  10 8 /
\plot 11 9  2 0  1 1  0 0 /
\plot 6 0  7 1  8 0  17 9 /
\plot 6 6  12 0 /

\plot 10.1 2  7.1 5 /
\plot 9.9 2  6.9 5 /

\multiput{$\ss \square$} at 19 5  10 2 /
\multiput{$\circ$} at 16 2  7 5   /
\put{$w_0$} at 26 -1
\put{$v_0$} at 20 -1
\put{$u_0$} at 14 -1
\put{$w_1$} at  8 -1
\put{$v_1$} at 2 -1 
\put{$V$} at 17.5 2
\put{$W$} at 20.6 5
\put{$\Omega V$} at 5.3 4.9
\put{$\Omega W$} at 8 1.9
\put{$P_0$} at 8.6 8.4
\put{$Q_0$} at 14.4 8.4 

\multiput{$\bullet$} at 10 8  16 8 /
\multiput{$\ss \bullet$} at 2 0  8 0   14 0  20 0  26 0  /
\endpicture} at 28 0

\endpicture}
$$
	\medskip
(2) We use successively such labels in order to obtain the first $z$ terms of
minimal projective resolutions (and the corresponding syzygy modules) of say modules 
$V \subseteq W$ with socle $T$, where $\pd V \ge z,\ \pd W\ge z$
(the case, we are interested
in, is $V = T,\ W = IT,$ but first we look at the general case in order to
highlight the pattern).
	\medskip
We assume that $V\subseteq W$ are modules with socle $T$, such that
$\pd V \ge z$ and $\pd W = z.$ Let 
$$ 
 u_0 = T, \quad v_0 = \tau^-\top V,\quad w_0 = \tau^-\top W.
$$
Also, let $S = \top \Omega^z W.$
	\medskip 
Let us denote by $Q_j$ and $P_j$ with $0\le j \le z$  the first projective modules which
occur in a minimal projective resolution of $V$ and $W$, respectively. For $j$ odd, we have
$P_j = Q_j$, thus let $B_j = P_j = Q_j$ for $j$ odd.

	\medskip
In case $z$ is odd, say $z = 2t+1$, the minimal projective
resolutions of $V$ and $W$ are of the form 
$$
\align
 \cdots @>>> B_{2t+1} @>>> P_{2t} @>>> \cdots @>>> B_3 @>>> P_2 @>>> B_1 @>>> &P_0 @>>> V @>>> 0, \cr
 0 @>>> B_{2t+1} @>>> Q_{2t} @>>> \cdots @>>> B_3 @>>> Q_2 @>>> B_1 @>>> &Q_0 @>>> W @>>> 0. \cr
\endalign
$$

In case $z$ is even, say $z=2t$, the minimal projective
resolutions of $V$ and $W$ look as follows
$$
\align
 \cdots @>>> 
 P_{2t} @>>> B_{2t-1} @>>> \cdots @>>> B_3 @>>> P_2 @>>> B_1 @>>> &P_0 @>>> V @>>> 0, \cr
 0 @>>> Q_{2t} @>>> B_{2t-1} @>>> \cdots @>>> B_3 @>>> Q_2 @>>> B_1 @>>> &Q_0 @>>> W @>>> 0. \cr
\endalign
$$

We define for $i\ge 1$
$$
 u_i = \soc B_{2i-1}, \quad v_i = \soc P_{2i-2}, \quad w_i = \soc Q_{2i-2}.
$$
Thus, we obtain for $z = 2t+1$ a sequence of integers 
$$
 v_{t+1} \le  w_{t+1} < u_t < v_t \le w_t < \cdots 
  < u_1 < v_1 \le w_1 < u_0 < v_0 \le w_0.
$$
and for $z = 2t$ the subsequence
$$
  u_t < v_t \le w_t < \cdots 
  < u_1 < v_1 \le w_1 < u_0 < v_0 \le w_0.
$$

Always, we have 
$$
 B_{2i+1} = [u_{2i+1},u_{2i}[, \quad
 P_{2i} = [v_{2i+1},v_{2i}[, \quad
 Q_{2i} = [w_{2i+1},w_{2i}[, \quad
$$
and 
$$
\alignat 2
 \Omega^{2i} V &= [u_{2i},v_{2i}[, &
 \Omega^{2i+1} V &= [v_{2i+1},u_{2i}[,\cr
 \Omega^{2i} W &= [u_{2i},w_{2i}[, & \quad
 \Omega^{2i+1} W &= [v_{2i+1},u_{2i}[.
\endalignat
$$
This is the labelling of the relevant modules which we are going to use. 
	\medskip
For a better vision, let us exhibit the special cases $z = 6$ and $z = 7$. Here,
$t = 3$. First, we show the case $z = 2t = 6$, and second the case $z = 2t+1 = 7$. 
$$
{\beginpicture
    \setcoordinatesystem units <.17cm,.17cm>
\put{\beginpicture
\multiput{} at 10 0  78 10 /
\setdots <1mm>
\plot 9 0  79 0 /
\setsolid 
\plot 12 0  21 9 /
\plot 18 0  27 9 /
\plot 24 0  33 9 /
\plot 30 0  39 9 /
\plot 36 0  45 9 /
\plot 42 0  51 9 /
\plot 48 0  57 9 /
\plot 54 0  63 9 /
\plot 60 0  69 9 /

\plot 66 0  72 6 /
\plot 72 0  74 2  /
\plot 68 8  76 0  77 1  78 0  /

\plot 13 3   16 0  17 1  18 0 /
\plot 17 5  22 0  23 1  24 0 /
\plot 20 8  28 0  29 1  30 0 /
\plot 26 8  34 0  35 1  36 0 /
\plot 32 8  40 0  41 1  42 0 /
\plot 38 8  46 0  47 1  48 0 /
\plot 44 8  52 0  53 1  54 0 /
\plot 50 8  58 0  59 1  60 0 /
\plot 56 8  64 0  65 1  66 0 /
\plot 62 8  70 0  71 1  72 0 /

\multiput{$\ss \bullet$} at  18 0  24 0  30 0  36 0  42 0  48 0  54 0 
   60 0  66 0  72 0  78 0   /
\put{$u_3$} at 12 -1.5 
\put{$v_3$} at 18  -1.5 
\put{$w_3$} at 24 -1.5 
\put{$u_2$} at 30 -1.5 
\put{$v_2$} at 36  -1.5 
\put{$w_2$} at 42 -1.5 
\put{$u_1$} at 48 -1.5 
\put{$v_1$} at 54  -1.5 
\put{$w_1$} at 60 -1.5 
\put{$u_0$} at 66 -1.5 
\put{$v_0$} at 72  -1.5 
\put{$w_0$} at 78 -1.5 

\multiput{$ \circ$} at 14 2  23 5  32 2  41 5  50 2  59 5  68 2   /
\multiput{$\sssize \square$} at  17 5  26 2  35 5  44 2  53 5  62 2  71 5 /
\multiput{$\bullet$} at   20 8  26 8  32 8  38 8  44 8  50 8  56 8  62 8  68 8 /

\put{$\ss\langle 1\rangle$} at 64.1 2 
\put{$\ss\langle 3\rangle$} at 46.1 2
\put{$\ss\langle 5\rangle$} at 28.1 2 

\put{$\ss\langle 6\rangle$} at 19.1 5 
\put{$\ss\langle 4\rangle$} at 37.1 5 
\put{$\ss\langle 2\rangle$} at 55.1 5 

\setshadegrid span <.4mm>
\vshade 26 2 2  <z,z,z,z>  28.1 -.1 4.1  <z,z,z,z> 30.2 2 2  /
\vshade 44 2 2  <z,z,z,z>  46.1 -.1 4.1  <z,z,z,z> 48.2 2 2  /
\vshade 62 2 2  <z,z,z,z>  64.1 -.1 4.1  <z,z,z,z> 66.2 2 2  /

\vshade 17 5 5  <z,z,z,z>  19.1  3 7.1  <z,z,z,z> 21.2 5 5  /
\vshade 35 5 5  <z,z,z,z>  37.1  3 7.1  <z,z,z,z> 39.2 5 5  /
\vshade 53 5 5  <z,z,z,z>  55.1  3 7.1  <z,z,z,z> 57.2 5 5  /

\put{$\ss Q_0$} at 67.5 9.5
\put{$\ss P_0$} at 61.5 9.5
\put{$\ss B_1$} at 55.5 9.5
\put{$\ss Q_2$} at 49.5 9.5
\put{$\ss P_2$} at 43.5 9.5
\put{$\ss B_3$} at 37.5 9.5
\put{$\ss Q_4$} at 31.5 9.5
\put{$\ss P_4$} at 25.5 9.5
\put{$\ss B_5$} at 19.5 9.5

\put{$V$} at 69.8 2
\put{$W$} at 73.4 5.3
\put{$\ss \Omega^6 V$} at 11.4 2
\put{$\ss \Omega^6 W$} at 14.4 5.7  
\put{$S$} at 22 -3.5
\put{$T$} at 66 -3.5
\arr{22 -2.5}{22 -.8}

\plot 14.1 1.9  17.1 4.9 /
\plot 13.9 2.1  16.9 5.1 /
\plot 32.1 1.9  35.1 4.9 /
\plot 31.9 2.1  34.9 5.1 /
\plot 50.1 1.9  53.1 4.9 /
\plot 49.9 2.1  52.9 5.1 /
\plot 68.1 1.9  71.1 4.9 /
\plot 67.9 2.1  70.9 5.1 /

\plot 25.9 1.9  22.9 4.9 /
\plot 26.1 2.1  23.1 5.1 /
\plot 43.9 1.9  40.9 4.9 /
\plot 44.1 2.1  41.1 5.1 /
\plot 61.9 1.9  58.9 4.9 /
\plot 62.1 2.1  59.1 5.1 /

\endpicture} at 0 0 

\put{\beginpicture
\multiput{} at 0 0  78 10 /
\setdots <1mm>
\plot -1 0  79 0 /
\setsolid 
\plot 0 0  9 9 /
\plot 6 0  15 9 /
\plot 12 0  21 9 /
\plot 18 0  27 9 /
\plot 24 0  33 9 /
\plot 30 0  39 9 /
\plot 36 0  45 9 /
\plot 42 0  51 9 /
\plot 48 0  57 9 /
\plot 54 0  63 9 /
\plot 60 0  69 9 /

\plot 66 0  72 6 /
\plot 72 0  74 2  /
\plot 68 8  76 0  77 1  78 0  /
\plot 1 3  4 0  5 1  6 0 /

\plot 5 5  10 0  11 1  12 0 /
\plot 8 8  16 0  17 1  18 0 /
\plot 14 8  22 0  23 1  24 0 /
\plot 20 8  28 0  29 1  30 0 /
\plot 26 8  34 0  35 1  36 0 /
\plot 32 8  40 0  41 1  42 0 /
\plot 38 8  46 0  47 1  48 0 /
\plot 44 8  52 0  53 1  54 0 /
\plot 50 8  58 0  59 1  60 0 /
\plot 56 8  64 0  65 1  66 0 /
\plot 62 8  70 0  71 1  72 0 /

\multiput{$\ss \bullet$} at 0 0  6 0  12 0  18 0  24 0  30 0  36 0  42 0  48 0  54 0 
   60 0  66 0  72 0  78 0   /
\put{$v_4$} at 0 -1.5 
\put{$w_4$} at 6 -1.5 
\put{$u_3$} at 12 -1.5 
\put{$v_3$} at 18  -1.5 
\put{$w_3$} at 24 -1.5 
\put{$u_2$} at 30 -1.5 
\put{$v_2$} at 36  -1.5 
\put{$w_2$} at 42 -1.5 
\put{$u_1$} at 48 -1.5 
\put{$v_1$} at 54  -1.5 
\put{$w_1$} at 60 -1.5 
\put{$u_0$} at 66 -1.5 
\put{$v_0$} at 72  -1.5 
\put{$w_0$} at 78 -1.5 

\multiput{$ \circ$} at 5 5  14 2  23 5  32 2  41 5  50 2  59 5  68 2   /
\multiput{$\sssize \square$} at 8 2  17 5  26 2  35 5  44 2  53 5  62 2  71 5 /
\multiput{$\bullet$} at 8 8  14 8  20 8  26 8  32 8  38 8  44 8  50 8  56 8  62 8  68 8 /

\put{$\ss\langle 1\rangle$} at 64.1 2 
\put{$\ss\langle 3\rangle$} at 46.1 2
\put{$\ss\langle 5\rangle$} at 28.1 2 
\put{$\ss\langle 7\rangle$} at 10.1 2 

\put{$\ss\langle 6\rangle$} at 19.1 5 
\put{$\ss\langle 4\rangle$} at 37.1 5 
\put{$\ss\langle 2\rangle$} at 55.1 5 

\setshadegrid span <.4mm>
\vshade 8 2 2  <z,z,z,z>  10.1 -.1 4.1  <z,z,z,z> 12.2 2 2  /
\vshade 26 2 2  <z,z,z,z>  28.1 -.1 4.1  <z,z,z,z> 30.2 2 2  /
\vshade 44 2 2  <z,z,z,z>  46.1 -.1 4.1  <z,z,z,z> 48.2 2 2  /
\vshade 62 2 2  <z,z,z,z>  64.1 -.1 4.1  <z,z,z,z> 66.2 2 2  /

\vshade 17 5 5  <z,z,z,z>  19.1  3 7.1  <z,z,z,z> 21.2 5 5  /
\vshade 35 5 5  <z,z,z,z>  37.1  3 7.1  <z,z,z,z> 39.2 5 5  /
\vshade 53 5 5  <z,z,z,z>  55.1  3 7.1  <z,z,z,z> 57.2 5 5  /

\put{$\ss Q_0$} at 67.5 9.5
\put{$\ss P_0$} at 61.5 9.5
\put{$\ss B_1$} at 55.5 9.5
\put{$\ss Q_2$} at 49.5 9.5
\put{$\ss P_2$} at 43.5 9.5
\put{$\ss B_3$} at 37.5 9.5
\put{$\ss Q_4$} at 31.5 9.5
\put{$\ss P_4$} at 25.5 9.5
\put{$\ss B_5$} at 19.5 9.5
\put{$\ss Q_6$} at 13.5 9.5
\put{$\ss P_6$} at 7.5 9.5

\put{$V$} at 69.8 2
\put{$W$} at 73.4 5.3
\put{$\ss \Omega^7 V$} at 2.8 5.7
\put{$\ss \Omega^7 W$} at 5.2 2 
\put{$S$} at 10 -3.5
\put{$T$} at 66 -3.5

\arr{10 -2.5}{10 -.8}

\plot 14.1 1.9  17.1 4.9 /
\plot 13.9 2.1  16.9 5.1 /
\plot 32.1 1.9  35.1 4.9 /
\plot 31.9 2.1  34.9 5.1 /
\plot 50.1 1.9  53.1 4.9 /
\plot 49.9 2.1  52.9 5.1 /
\plot 68.1 1.9  71.1 4.9 /
\plot 67.9 2.1  70.9 5.1 /

\plot 7.9 1.9  4.9 4.9 /
\plot 8.1 2.1  5.1 5.1 /
\plot 25.9 1.9  22.9 4.9 /
\plot 26.1 2.1  23.1 5.1 /
\plot 43.9 1.9  40.9 4.9 /
\plot 44.1 2.1  41.1 5.1 /
\plot 61.9 1.9  58.9 4.9 /
\plot 62.1 2.1  59.1 5.1 /

\endpicture} at 0 -20
\put{$t = 3 \quad z = 2t = 6$} [l] at -40 10
\put{$t = 3 \quad z = 2t+1 = 7$} [l] at -40 -10
\endpicture}
$$
The shaded areas with labels of the form $\langle i \rangle$ will be defined in (3):
	\medskip
(3) We define $\langle j \rangle$ for $j\ge 1$ as follows: $\langle 2i \rangle$ is the
set of modules $[x,y[$ with 
$u_{i}\le x < v_{i}$ and $w_{i}\le y < u_{i-1}$
whereas $\langle 2i-1\rangle$ is the set of modules $[x,y[$ with 
$w_{i} \le x < u_{i-1}$ and $u_{i-1} \le y < v_{i-1}$. 
The following assertion is essential:
	\medskip
{\it If $M$ belongs to $\langle j\rangle$ with $j\ge 2$, then $\Sigma M\in \langle j-1\rangle.$}
	\medskip
Proof. First, let $j = 2i$ be even (thus $i\ge 1$).
Let $M = [x,y[\ \in \langle 2i\rangle,$ thus
$u_{i}\le x < v_{i}$ and $w_{i}\le y < u_{i-1}$. 
Then $IM = Ix$. 
Since $B_{2i-1} = [u_{i},u_{i-1}[$
has $x$ and $\tau u_{i-1}$ as composition factors, we see that $\tau u_{i-1}$ is a
composition factor of $Ix$, thus $\tau u_{i-1} \le \top IM$. 
On the other hand, $P_{2i-2} = [v_{i},v_{i-1}[$, thus 
$\top P_{2i-2} = \tau v_{i-1}$. There is no indecomposable module with both $x$ and $\tau v_{i-1}$
as compostion factors, since such a module would have $P_{2i-2}$ as a proper factor module. 
It follows that $\tau v_{i-1}$
is not a composition factor of $IM$. Altogether we see that 
$\tau u_{i-1} \le \top IM < \tau v_{i-1}$, thus $u_{i-1} \le \tau^-\top IM < v_{i-1}$.
Therefore $\Sigma M = [y,\tau^-\top IM[$ belongs to $\langle 2i-1\rangle.$

Next, let $j = 2i-1$ be odd (thus $j\ge 2$).  Let 
$M = [x,y[\ \in \langle 2i-1\rangle,$ thus
$w_{i} \le x < u_{i-1}$ and $u_{i-1} \le y < v_{i-1}$.
Again, $IM = Ix$. The existence of $Q_{2i-2} = [w_{i},w_{i-1}[$ shows that 
$\tau w_{i-1}$ is a composition factor of $IM$, thus $\tau w_{i-1} \le \top IM$.
The fact that $B_{2i-3} = [u_{i-1},u_{i-2}[$ is projective implies that
$\top IM < \tau u_{i-2}$. Altogether we see that $\tau w_{i-1} \le \top IM < \tau u_{i-2}$,
thus $w_{i-1} \le \tau^-\top IM < u_{i-2}$. 
Therefore $\Sigma M = [y,\tau^-\top IM[$ belongs to $\langle 2i-1\rangle.$

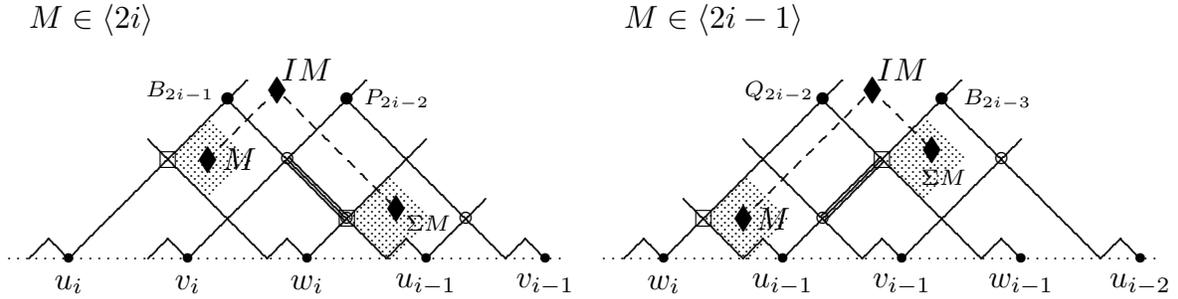
$$
{\beginpicture
    \setcoordinatesystem units <.22cm,.22cm>
\put{\beginpicture
\multiput{} at 0 0  26 10 /
\put{$M \in \langle 2i \rangle$} [l] at 0 12
\setdots <1mm>
\plot -1 0  27 0 /
\setsolid
\plot 12 0  13 1  14 0  20 6 / 
\plot 19 5  24 0  25 1  26 0 /
\plot 16 2  18 0  19 1  20 0  23 3 /
\plot 16 8  19 5 /
\plot 16 2  10 8 /
\plot 11 9  2 0  1 1  0 0 /
\plot 6 0  7 1  8 0  17 9 /
\plot 6 6  12 0 /
\multiput{$\ss \square$} at 16 2   7 5  /
\multiput{$\circ$} at 13 5  16 2  22 2  /
\put{$v_{i-1}$} at 26 -1.3
\put{$u_{i-1}$} at 20 -1.3
\put{$w_{i}$} at 14 -1.3
\put{$v_{i}$} at  8 -1.3
\put{$u_{i}$} at 2 -1.3
\put{$\ss B_{2i-1}$} at 7.6 8.4
\put{$\ss P_{2i-2}$} at 18.5 8 

\multiput{$\bullet$} at 10 8  16 8 /
\multiput{$\ss \bullet$} at 2 0  8 0   14 0  20 0  26 0  /

\setshadegrid span <.4mm>
\vshade 7 5 5  <z,z,z,z>  9.1  3 7.1  <z,z,z,z> 11.2 5 5  /
\vshade 16 2 2  <z,z,z,z>  18.1 -.1 4.1  <z,z,z,z> 20.2 2 2  /

\multiput{$\blacklozenge$} at 9 5  12.5 8.5  18.5  2.5  /
\setdashes <1mm>
\plot 9 5  12.5 8.5  18.5  2.5 /

\put{$M$} at 10.6 5 
\put{$IM$} at 14  9.5 
\put{$\ss \Sigma M$} at 20.1 1.8 

\setsolid
\plot 13.1 5.1  16.1  2.1 /
\plot 12.9 4.9  15.9  1.9 /

\endpicture} at 0 0 
\put{\beginpicture
\multiput{} at 0 0  26 10 /
\put{$M \in \langle 2i-1 \rangle$} [l] at 0 12
\setdots <1mm>
\plot -1 0  27 0 /
\setsolid
\plot 12 0  13 1  14 0  19 5  24 0  25 1  26 0 /
\plot 16 2  18 0  19 1  20 0 /
\plot 16 8  19 5  20 6 /
\plot 16 2  10 8 /
\plot 11 9  2 0  1 1  0 0 /
\plot 3 3  6 0  7 1  8 0  17 9 /
\plot 6 6  12 0 /
\multiput{$\ss \square$} at 13 5  4 2 /
\multiput{$\circ$} at 10 2  19 5   /
\put{$u_{i-2}$} at 26 -1.3
\put{$w_{i-1}$} at 20 -1.3
\put{$v_{i-1}$} at 14 -1.3
\put{$u_{i-1}$} at  8 -1.3
\put{$w_{i}$} at 2 -1.3 
\put{$\ss Q_{2i-2}$} at 7.8 8.4
\put{$\ss B_{2i-3}$} at 18.8 8 

\multiput{$\bullet$} at 10 8  16 8 /
\multiput{$\ss \bullet$} at 2 0  8 0   14 0  20 0  26 0  /

\setshadegrid span <.4mm>
\vshade 4 2 2  <z,z,z,z>  6.1 -.1 4.1  <z,z,z,z> 8.2 2 2  /
\vshade 13 5 5  <z,z,z,z>  15.1  3 7.1  <z,z,z,z> 17.2 5 5  /

\multiput{$\blacklozenge$} at 6 2  12.5 8.5  15.5  5.5  /
\setdashes <1mm>
\plot 6 2  12.5 8.5  15.5  5.5 /
\put{$M$} at 7.5 2 
\put{$IM$} at 14  9.5 
\put{$\ss \Sigma M$} at 16.1 4.1 

\setsolid
\plot 13.1 4.9  10.1  1.9 /
\plot 12.9 5.1  9.9  2.1 /

\endpicture} at 30 0

\endpicture}
$$

	\medskip
(4) {\it If $V = T$, then $\langle 1 \rangle$ is the set of modules
$[x,u_0[$ with $w_1 \le x < u_0$. If $W = IT$, and $w_1 \le x < u_0$, then $\Sigma [x,u_0] = IT.$}
	\medskip

Proof. First, let $V = T$. Now $\langle 1 \rangle$ is the set of modules $[x,y]$ with
$w_1 \le x < u_0$ and $u_0 \le y < v_0.$ Let us assume that $V = T$. Then $V = [u_0,v_0[$ is
simple, thus $\tau v_0 = u_0.$ The condition $u_0 \le y < v_0$ is only satisfied for $y = v_0.$
This is the first assertion.

Next, Let $M = [x,u_0[$ with $w_1 \le x < u_0$. Then $IM = Ix$. The module $Q_0 = [w_1,w_0[$
has $[x,w_0[$ as a factor module, thus this is a submodule of $Ix.$ But $W = [u_0,w_0[$ is
a factor module of $[x,w_0[$. We assume here that $W = IT = Iu_0$. Since the injective module $W$
is a factor module of $[x,w_0[$, also $[x,w_0[$ has to be injective, thus equal to $Ix$.
It follows from $Ix = [x,w_0[$, that $\Sigma [x,u_0[\ = [u_0,w_0[\ = W = IT.$
$\s$
	\medskip
Remark. If $V = T$ and $W = IT$, let us define $\langle 0 \rangle = \{IT\}$. Then 
we have: If $M\in \langle 1\rangle$, then $\Sigma M \in \langle 0 \rangle$, thus
{\it If $M$ belongs to $\langle j \rangle$ with $j\ge 1$, then $\Sigma M \in \langle j-1\rangle.$}
	\medskip
In the case $V = T,\ W = IT$ (and this is the case we are 
interested in), the pictures which we have used are are not quite satisfactory. Here is
a better one of the relevant part:
$$
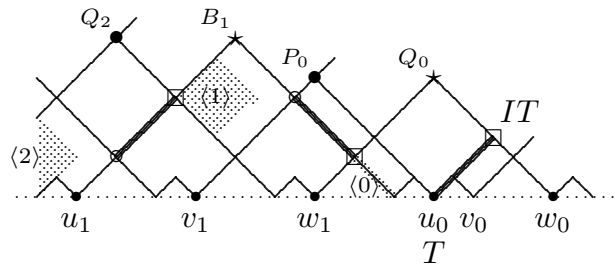
{\beginpicture
    \setcoordinatesystem units <.22cm,.22cm>
\multiput{} at 0 0  26 10 /
\setdots <1mm>
\plot -1 0  29 0 /
\setsolid
\plot 0 0  1 1  2 0  10 8  18 0  19 1 /
\plot 0 6  6 0  7 1  8 0  15 7 /
\plot 14 6  20 0  23 3 /
\plot 0 4  5 9 /
\plot 4 8  12 0  13 1  14 0  20 6  26 0 27 1  28 0  /
\plot 21 1  22 0  25 3 /
\multiput{$\star$} at 10 8  20 6  /
\multiput{$\ss \square$} at 7 5  16 2  23 3 /
\multiput{$\circ$} at 4 2  13 5     /
\put{$w_{0}$} at 26 -1.3
\put{$v_{0}$} at 22 -1.3
\put{$u_{0}$} at 20 -1.3
\put{$w_{1}$} at 14 -1.3
\put{$v_{1}$} at  8 -1.3
\put{$u_{1}$} at 2 -1.3
\put{$\ss Q_{2}$} at 3 9 
\put{$\ss B_{1}$} at 9 9 
\put{$\ss P_0$} at 13 7 
\put{$\ss Q_0$} at 19 7

\put{$T$} at 20 -2.8
\put{$IT$} at 24.3 4.3 
\multiput{$\bullet$} at 4 8  14 6 /
\multiput{$\ss \bullet$} at 2 0  8 0   14 0  20 0  26 0  /
\put{$\ss \langle 2\rangle$} at -.6 2
\put{$\ss \langle 1\rangle$} at 9 5
\put{$\ss \langle 0\rangle$} at 16.5 .5 

\setshadegrid span <.4mm>
\vshade   0.1 -.1 4.1  <z,z,z,z> 2.2 2 2  /
\vshade 7 5 5  <z,z,z,z>  9.1  3 7.1  <z,z,z,z> 11.2 5 5  /
\vshade 15.8 1.8 1.8  <z,z,z,z>  16.2 1.4 2.2 <z,z,z,z> 17.8 -.2 .4 <z,z,z,z> 18.2 0.2 0.2  /

\setsolid
\plot 4.07 1.93  7.07 4.93 /
\plot 3.93 2.07  6.93 5.07 /

\plot 16.07 2.07  13.07 5.07 /
\plot 15.93 1.93  12.93 4.93 /

\plot 20.07 -.07  23.07 2.93 /
\plot 19.93 .07  22.93 3.07 /

\endpicture}
$$
	\medskip
(5) {\it If $z$ is even, then $\Sigma^z PS = IT$. If
$z$ is odd, then $\Sigma^z S = IT$.}
	\medskip
Proof. If $z$ is even, then $PS$ belongs to $\langle z\rangle.$
If $z$ is odd, then $S$ belongs to $\langle z\rangle.$
$\s$
	\medskip
(6) In order to complete the proof of (a), we still have to show:
If $z$ is even, $\Omega^zIT = PS$ and $\pd T \ge z$, then
$\id S \ge z.$ We can assume that $z = 2t\ge 2$. 
WAs we have seen in (2),
$\Omega^{z-1}IT = [w_t,u_{t-1}[$ and $S = \tau w_t = \tau(\soc \Omega^{z-1})$. Since 
$\Omega^{z-1}IT$ is not projective, the module $[\tau w_t,u_{t-1}[$ exists: it
is the middle term of an element in $\Ext^1(\Omega^{z-1}IT,S)$. This shows that 
$\Ext^z(IT,S) = \Omega^{z-1}IT,S) \neq 0.$

For assertion (b), it remains to show: If $z$ is odd, $\Omega^z T = PS$ and $\pd IT \ge z$, then
$\id PS \ge z.$ But this is trivial. Namely, $\Omega^z T = PS$ implies that
$\Ext^z(T,PS) \neq 0$, thus $\id PS \ge z.$ 

This completes the proof of Theorem 4.1 in case $A$ is a linear Nakayama algebra.
$\s$
	\bigskip
{\bf 5.2. Remark.} We use the previous proof in order to mention the following observation.
Recall that for a representation-finite algebra $A$, an indecomposable module $X$ is
a {\it predecessor} of an indecomposable module $Y$, provided there is a path
in the Auslander-Reiten quiver starting in $X$ and ending in $Y$. 
	\medskip
{\it Let $A$ be a linear Nakayama algebra. Let $V \subseteq W$ be indecomposable modules
and let $P_\bullet, Q_\bullet$ be minimal projective 
resolutions of $V$ and $W$, respectively. If both $\Omega^{2t}V$ and $\Omega^{2t}W$
are non-zero, then $\Omega^{2t}W $ is a predecessor of
all the modules $P_i, Q_i$ with $0 \le i < 2t.$ If both $\Omega^{2t+1}V$ and $\Omega^{2t+1}W$
are non-zero, then $\Omega^{2t+1}V $ is a predecessor of
all the modules $P_i, Q_i$ with $0 \le i \le 2t.$}
	\medskip
Proof. See (2) in section 5.2. 
$\s$
	\bigskip
{\bf 5.3. Representations of the quiver $\Bbb Z$.} 
Before we look at cyclic Nakama algebras,
let us consider $\Bbb Z$ as a quiver (its vertex set is the set $\Bbb Z,$ and for any integer $z$, there is an arrow $z \to z\!-\!1$). Let $\rho$ be a set of relations (that means: paths of finite
length $\ge 2$). 
The category $\Cal M = \mod (\Bbb Z,\rho)$ consists of the finite-dimensional representations of
the quiver $\Bbb Z$ which satisfy the relations $\rho$. We assume that $\mod (\Bbb Z,\rho)$
has enough projective and enough injective objects, thus we have projective resolutions and
injective coresolutions in $\Cal M,$ as well as Auslander-Reiten sequences. Thus, as in the
case of Nakayama algebras, we usually will specify $(\Bbb Z,\rho)$ by indicating the
shape of the Auslander-Reiten quiver.

The category $\Cal M$ is quite similar to the module category of a linear Nakayama algebra,
in particular all indecomposable objects $M$ are serial; if $M$ has support $\{x, 
x\!+\!1,\dots, y\}$, we again may write $M = [x,y\!+\!1[$. Of course, an object in $\Cal M$
is simple iff it is of the form $[x,x\!+\!1[$ for some vertex $x$, and
by abuse of language, we will write just $x$ instead of $[x,x\!+\!1[$,
provided no confusion is possible. 

Let us stress that in contrast to  
the module category of a Nakayama algebra, given a vertex $x$, both representations
$x$ and $Ix$ may have infinite projective dimension. Here is a typical example:
$$
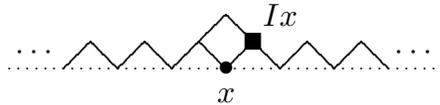
{\beginpicture
    \setcoordinatesystem units <.3cm,.3cm>
\multiput{} at -8 0  8 2 /
\setdots <1mm>
\plot -8 0  8 0 /
\setsolid 
\plot -6 0  -5 1  -4 0  -3 1 -2 0  0 2  2 0  3 1  4 0  5 1  6 0 /
\plot -1 1  0 0  1 1 /
\put{$x$} at 0 -1 
\put{$Ix$} at 2 2
\put{$\bullet$} at 0 0
\put{$\ss \blacksquare$} at 1 1
\multiput{$\cdots$} at -7 0.6  7 .6 /
\endpicture}
$$
However, the analogue of Theorem 4.1 holds true in $\Cal M$.
	\medskip
{\bf Proposition.} {\it Let $S, T$ be simple objects in $\Cal M$, and let
$z\in \Bbb N_0$.} 

(a) {\it If $z$ is even, $\Omega^z IT = PS$, and $\pd T \ge z$, then $\Sigma^z PS = IT$
and $\id S \ge z.$}

(b) {\it If $z$ is odd, $\Omega^z T = PS$, and $\pd IT \ge z$, then $\Sigma^z S = IT$
and $\id PS \ge z.$}
	\medskip
Proof. We show (a); the proof of (b) is similar. 
For any natural number $m$, let
$A_m$ be the path algebra of the full subquiver $[-m,m]$ satisfying the relations in $\rho$
with support in $[-m,m],$ thus $\Cal M = \bigcup_{m\in\Bbb N} \mod A_m$. 
	\medskip
Take  minimal projective resolutions of $T$ and $IT$, say
$$
\align
  \cdots \to   P_z @>>> 
    P_{z-1} @>>> \cdots @>>> 
     P_1 @>>>   P_0 @>>> 
  &  T  @>>> 0, \cr
  \cdots \to   Q_z @>>> 
    Q_{z-1} @>>> \cdots @>>> 
     Q_1 @>>>   Q_0 @>>> 
  &I  T  @>>> 0.
\endalign
$$ 
Since $\Omega^z IT = PS$, we have $Q_z = PS$. 
Choose a natural number $m \ge 1$ such that all the modules $P_i, Q_i$ with $0 \le i \le z$
have support in $[-m,m]$, thus they are $A_m$-modules, and of course projective also
as $A_m$-modules. We have shown in 5.1 that Theorem 4.1 holds for linear Nakayama
algebras, thus for $A = A_m$. In $\mod A$, we have 
$(\Sigma_A)^z PS = IT$
and $\id_A S \ge z$. Let
$$
\align
  0 \to   S @>>> 
    I_0 @>>>  I_1 @>>> \cdots @>>> 
     I_{z-1} @>>>   &(\Sigma_A)^zS 
   @>>> 0, \cr
  0 \to   PS @>>> 
    J_0 @>>>  J_1 @>>> \cdots @>>> 
     J_{z-1} @>>>   &(\Sigma_A)^zPS
   @>>> 0.
\endalign
$$ 
be exact sequences in $\mod A$ such that all the modules $I_i, J_i$ with $0\le i < z$ are
indecomposable injective in $\mod A$.
It remains to be seen that the representations $I_i, J_i$ with $0\le i < z$ are injective
in $\Cal M.$ According to the dual of 5.2, all these modules are predecessors of 
$IT$, and $IT$ is by definition injective in $\Cal M$. The following quite obvious
assertion finishes the proof of (a). 
{\it Let $I$ be an indecomposable injective object of $\Cal M$ which is an $A_m$-module.
Let $M$ be an indecomposable $A_m$-module which is a predecessor of $I$. If $M$ is
injective as an $A_m$-module, then $M$ is injective in $\Cal M$.} 
$\s$
	\bigskip
{\bf 5.4. The case of $A$ being a cyclic Nakayama algebra.}
Let $A$ be a cyclic Nakayama algebra. In order to show (a), let $S, T$ be simple modules
with $\Omega^z IT = PS$ and $\pd T \ge z,$ where $z$ is even. 

Let $\widetilde A$ be a universal cover of $A$
with push-down functor $\pi\:\mod \widetilde A \to \mod A.$ 
Note that $\mod \widetilde A = \mod (\Bbb Z,\rho)$ for a suitable set $\rho$ of relations.
There is a natural number $b$ such that 
all indecomposable projective (and thus also all indecomposable injective) $A$-modules have 
length at most $b$. Thus $\mod \widetilde A$ has enough projective and enough injective
modules (and all indecomposable $\widetilde A$-modules which are projective or injective have
length at most $b$). This shows that Proposition 5.3 holds for $\mod \widetilde A$.

Fix an $\widetilde A$-module $\widetilde T$ with $\pi(\widetilde T) = T.$
Let $I \widetilde T$ be an injective envelope of $\widetilde T$ (as an $\widetilde A$-module).
Then $\pi(I \widetilde T) = IT$. 
Let 
$$
\align
  \cdots \to \widetilde P_z @>>> 
  \widetilde P_{z-1} @>>> \cdots @>>> 
   \widetilde P_1 @>>> \widetilde P_0 @>>> 
  &\widetilde T \to 0 \cr
  \cdots \to \widetilde Q_z @>>> 
  \widetilde Q_{z-1} @>>> \cdots @>>> 
   \widetilde Q_1 @>>> \widetilde Q_0 @>>> 
  &I\widetilde T \to 0 
\endalign
$$ 
be minimal projective resolutions of $\widetilde T$ and $I\widetilde T$
(as $\widetilde A$-modules). 
Applying $\pi$ to these resolutions
we obtain minimal projective resolutions of $T$ and $IT$
(as $A$-modules), respectively. Since $\pd IT = z$, we see that the
map $\widetilde Q_z \to \widetilde Q_{z-1}$ is injective and that
$\pi(\widetilde Q_z) = \Omega^zIT = PS$. In particular, $\widetilde Q_z$
is indecomposable projective, say equal to $P\widetilde S$ with $\widetilde S$
a simple $\widetilde A$-module, and we have $\pi(\widetilde S) = S.$

Altogether, we see: $\widetilde S$ and $\widetilde T$ are simple
$\widetilde A$-modules with $\Omega^zI\widetilde T = P\widetilde S$ and
$\pd \widetilde T \ge z$ (and $z$ is even). Proposition 5.3 (a) asserts that
$\Sigma^zP\widetilde S = I\widetilde T$ and that $\id \widetilde S \ge z.$
Now using the push-down functor $\pi$, it follows that 
$\Sigma^zPS = IT$ and that $\id S \ge z.$ 

This completes the proof of (a). The proof of (b) is, of course, similar.
$\s$
	\medskip
This completes the proof of Theorem 4.1.$\s$
	\bigskip
{\bf 5.5. An example.} Let $A$ be a linear Nakayama algebra, and $S,T$ are simple modules, with
$PS = \Omega^zIT$ and $\Sigma^zPS$, where $z = 2t$ is even. Let $P_\bullet$ be a minimal
projective resolution of $IT$, and $I_\bullet$   a minimal injective coresolution of $PS$,
then these two sequences yield non-zero elements of $\Ext^z(IT,PS)$. 

We present one example in detail, here is $z = 4$, thus $t=2.$ There is a commutative diagram
of the following form
$$
\CD
0 @>>> \Omega^4IT @>>> P_3 @>>> P_0 @>>> P_1 @>>> P_0 @>>> IT @>>> 0 \cr
@.       @|             @VVV     @VVV      @VVV     @VVV    @|     @. \cr
0 @>>>    PS      @>>> I_0 @>>> I_1 @>>> I_2 @>>> I_3 @>>> \Sigma^4 PS  @>>> 0 \cr
\endCD
$$
which yields canonical maps $\Omega^iIT \to \Sigma^{z-i}PS$. 

For the example which we have selected, the map $\Omega^2IT \to \Sigma^{2}PS$
is neither a monomorphism nor an epimorphism: 
$$
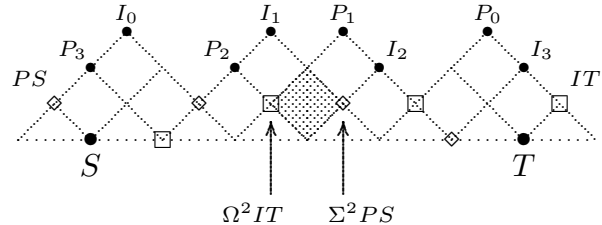
{\beginpicture
    \setcoordinatesystem units <.4cm,.4cm>
\multiput{} at 0 0  16 3 /
\setdots <1mm>
\plot 0 0  16 0 /
\setdots <.5mm>
\multiput{$\ss\bullet$} at 2 2  3 3  6 2  7 3  9 3  10 2  13 3  14 2 /
\multiput{$\bullet$} at 2 0  14 0 /
\plot 0 0  3 3  6 0  9 3  12 0  14 2 /
\plot 1 1  2 0  4 2 /
\plot 2 2  4 0  7 3  10 0  13 3  16 0  /
\plot 6 2  8 0  10 2 /
\plot 12 2  14 0  15 1 /
\multiput{$\ss \square$} at  4 0  7 1  11 1  15 1 /
\multiput{$\diamond$} at 1 1  5 1  9 1  12 0 /
\put{$\ss IT$} at  15.7 1.7 
\put{$T$} at 14 -.7 
\put{$S$} at 2 -.7 
\put{$\ss PS$} at  0.3 1.7 
\setshadegrid span <.4mm>
\vshade 7 1 1  <z,z,z,z> 8 0 2  <z,z,z,z> 9 1 1  /
\put{$\ss P_0$} at 13 3.5  
\put{$\ss P_1$} at 9 3.5  
\put{$\ss P_2$} at 5.5 2.5  
\put{$\ss P_3$} at 1.5 2.5  

\put{$\ss I_0$} at 3 3.5  
\put{$\ss I_1$} at 7 3.5  
\put{$\ss I_2$} at 10.5 2.5  
\put{$\ss I_3$} at 14.5 2.5  
\put{$\ss \Omega^2IT$} at 6.5 -2
\put{$\ss \Sigma^2PS$} at 9.5 -2
\setsolid
\arr{7 -1.5}{7 0.5}
\arr{9 -1.5}{9 0.5}
\endpicture}
$$
	\bigskip\bigskip

{\bf Appendix A.  The functions $\psi$ and $\gamma$.}
	\medskip 
We assume that $A$ is a cyclic Nakayama algebra. 
	\medskip
{\bf A.1. Proposition.} {\it Let $A$ be a cyclic Nakayama algebra.
For all $t\ge 0$, we have }
$$
  \psi^t\gamma^t\psi^t = \psi^t \quad\text{and}\quad
 \gamma^t\psi^t\gamma^t = \gamma^t.
$$
	\medskip
The proof will be given in A.4, using some basic
properties of monotone endofunctions of $\Bbb Z,$ see A.3. 
	\medskip
{\bf Corollary 1.}  {\it Let $A$ be a cyclic Nakayama algebra. For all $t \ge 0,$ 
the map 
$$
 \gamma^t\:\Im\psi^t \to \Im\gamma^t
$$
is a bijection with inverse $\psi^t.$ 
In particular, we have $|\Im \psi^t| = |\Im \gamma^t|$.}
$\s$
	\medskip
Recall that $a(A)$ is the smallest natural number $a$ with $\Im \psi^a = \Im \psi^{a+1}$.
Using duality, the smallest natural number $b$ with $\Im \gamma^b = \Im \gamma^{b+1}$ is 
$a(A^{\op})$.

We denote by $c(A)$ the number of $\psi$-cyclic modules.  
Thus, using duality, $c(A^{\op})$ is the number of $\gamma$-cyclic $A$-modules.
	\medskip
{\bf Corollary 2.} {\it Let $A$ be a cyclic Nakayama algebra. Then}
$$
 a(A) = a(A^{\op}) \qquad\text{\it and}\qquad c(A) = c(A^{\op}).
$$
	\medskip
Proof. 
We have $\Im \psi^0 \supseteq \Im \psi^1 \supseteq \Im \psi^2 \supseteq 
\cdots$. If $S$ is a simple module, then, by definition,  
$a(S) = t$ provided $S$ belongs to 
$\Im \psi^{t-1} \setminus\Im\psi^t,$ and 
$a(S) = \infty$ provided $S$ belongs to $\bigcap_t \Im \psi^t.$

Similarly, $\Im \gamma^0 \supseteq \Im \gamma^1 \supseteq \Im \gamma^2 \supseteq 
\cdots$. If $S$ is a simple module, let $a'(S) = t$ provided $S$ belongs to 
$\Im \gamma^{t-1} \setminus\Im\gamma^t,$ 
and let $a'(S) = \infty$ provided $S$ belongs to $\bigcap_t \Im \gamma^t.$
The simple module $S$ is $\gamma$-cyclic, provided $S$ belongs to $\Im \gamma^t$
for all $t$, thus provided $a'(S) = \infty$.
Using duality, we see that $a(A^{\op})$ is the 
maximum of the numbers $a'(S)$, where $S$ is simple and 
not $\gamma$-cyclic. 

Since $\gamma^t$ provides a bijection between $\Im \psi^t$ and $\Im\gamma^t$,
it provides a bijection from the set of isomorphism classes $[S]$ of
simple modules with $a(S) \ge t$ onto the set of isomorphism classes $[S]$ of
simple modules with $a'(S) \ge t$. It follows that 
$a(A) = a(A^{\op}).$
	\smallskip

For $t\ge a(A)$, 
the map $\gamma^t$ provides a bijection from the
set of $\psi$-cyclic simple modules onto the set of $\gamma$-cyclic modules.
This shows that $c(A) = c(A^{\op})$.
$\s$
	\bigskip
{\bf A.2. Proposition.} {\it Let $A$ be a cyclic Nakayama algebra.
For all $t\ge 0$, we have }
$$
  \psi^t\tau^-\gamma^t\tau\psi^t = \psi^t \quad\text{and}\quad
 \gamma^t\tau\psi^t\tau^-\gamma^t = \gamma^t.
$$
	\medskip
The proof will be similar to the proof of Proposition A.1, see A.5.
As in A.1, there is the following Corollary.
	\medskip
{\bf Corollary.}  {\it Let $A$ be a cyclic Nakayama algebra. For all $t \ge 0,$ 
the map 
$$
 \gamma^t\tau\:\Im\psi^t \to \Im\gamma^t
$$
is a bijection with inverse $\psi^t\tau^-.$ }
$\s$
	\medskip
Proof. The image of $\gamma^t\tau$ is contained in $\Im \gamma^t,$ thus $\gamma^t\tau$
maps $\Im\psi^t$ into $\Im\gamma^t.$ Similarly, $\psi^t\tau^-$ maps $\Im \gamma^t$ into
$\Im\psi^t.$ The equalities in A.2 yield the assertion.
$\s$
	\bigskip

{\bf A.3.} A (set-theoretical) function $f\:\Bbb Z \to \Bbb Z$ is said to be {\it monotone}
provided $i < j$ in $\Bbb Z$ implies $f(i) \le f(j).$
	\medskip
{\bf Lemma.} {\it Let $f,g\:\Bbb Z \to \Bbb Z$ be monotone functions with
$fg(i) \ge i \ge gf(i)$ for all $i\in \Bbb Z$. Then, for all $t\ge 0$, we have
$f^tg^t(i) \ge i \ge g^tf^t(i)$, and $f^tg^t f^t = f^t$ and $g^tf^tg^t = g^t.$}
	\medskip 
Proof. First, we show by induction on $t$ that $f^tg^t(i) \ge i$ for all $i\in \Bbb Z$.
The case $t=1$ is one of the assumptions. Thus, assume that we know already for some $t$
that $f^tg^t(i) \ge i$ for all $i$. Replacing $i$ by $g(i)$, we obtain
$f^tg^{t+1}(i) \ge g(i).$ Since $f$ is monotone, we get $f^{t+1}g^{t+1}(i) \ge fg(i).$
Altogether we see that $f^{t+1}g^{t+1}(i) \ge fg(i)\ge i$.
Similarly, we see that $i \ge g^tf^t(i)$ for all $i$, and all $t\ge 0.$

Applying  $f^t$ to $i \ge g^tf^t(i)$, the monotony gives $f^t(i) \ge f^tg^tf^t(i)$.
On the other hand, we take $f^tg^t(i) \ge i$ and replace $i$ by $f^t(i)$, this
yields $f^tg^tf^t(i) \ge f^t(i).$ Altogether we have 
$f^tg^tf^t(i) \ge f^t(i) \ge f^tg^tf^t(i)$, thus $f^tg^tf^t = f^t.$
Similarly, we see that $g^tf^tg^t = g^t.$
$\s$
	\bigskip
{\bf A.4.} Now, let $A$ be a cyclic Nakayama algebra.
In order to use A.3, we need to work with
the universal cover $\widetilde A$ of the algebra $A$ and a covering functor 
$\pi\:\mod \widetilde A \to \mod A$ (see [Ga] and related papers by 
Bongartz-Gabriel and Gordon-Green). The quiver of $\widetilde A$ is $\Bbb Z$
(we consider the integers $\Bbb Z$ as a quiver with vertex set $\Bbb Z$
and with arrows $i\to i\!-\!1$ for all $i\in \Bbb Z$). We have to fix a simple module
$S = S(0)$, define $S(i) = \tau^{-i} S$ and use as covering map $\pi\:i \mapsto S(i).$ 
Note that the Auslander-Reiten quiver of $\widetilde A$ is a full subquiver of the
translation  quiver $\Bbb Z\Bbb A_\infty$, its lower boundary is formed by the simple
modules. We should stress that the integers which index the simple modules 
{\bf increase} when going from right to left.

Let $\widetilde \psi(i) = i+|I(i)|$ and $\widetilde\gamma(i) = i-|P(i)|.$ Then
$\widetilde\psi$ is a covering of $\psi$, and $\widetilde\gamma$ is a
covering of $\gamma$ (this means: $\pi\widetilde\psi = \psi\pi$ and 
$\pi\widetilde\gamma = \gamma\pi$). Similarly, let 
$\widetilde\tau(i) = i-1$; thus $\widetilde\tau$ is a covering of $\tau.$ 
	\medskip
{\bf (1)} {\it The functions $\widetilde \psi$ and $\widetilde \gamma$ are monotone.}
	\medskip 
Proof. Let $i < j$ in $\Bbb Z$. If $j-i \ge |I(i)|$, then 
$\widetilde\psi(i) = i+|I(i)| \le j < j + |I(j)| = \widetilde\psi(j)$.
Thus, we can assume that $j-i < |I(i)|.$ Since $0 < j-i < |I(i)|$, there is a submodule $U$
of $I(i)$ of length $j-i$ and $I(i)/U$ has socle $j$. 
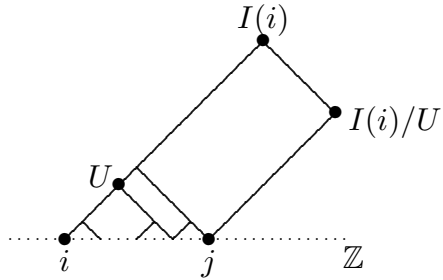
$$
{\beginpicture
    \setcoordinatesystem units <.2cm,.2cm>
\multiput{} at 0 0  15 11 /
\plot 0 0  11 11  /
\plot 1 1  2 0 /
\plot 4 0  5 1 /
\plot 3 3  6 0  7 1 /
\plot 4 4  8 0  15 7  11 11 /
\put{$i$\strut} at 0 -1.3
\put{$j$\strut} at 8 -1.3
\put{$I(i)$} at 11 12.2
\put{$U$} at 2 3.5
\put{$I(i)/U$} at 18.3 6.5 
\setdots <1mm>
\plot -3 0  16 0  /
\setsolid
\multiput{$\bullet$} at 0 0  8 0  3 3  11 11  15 7 /
\put{$\Bbb Z$} at 16 -1 
\endpicture}
$$
Since $I(i)/U$ is a module with socle $j$,
it is a submodule of $I(j)$, thus $|I(j)| \ge |I(i)/U| = |I(i)|-|U| = |I(i)|-j+i.$
This shows that  $\widetilde\psi(i) = i+|I(i)| \le j + |I(j)| = \psi(j)$. This shows that
$\widetilde\psi$ is monotone. 

In the same way, or using duality, one sees that $\widetilde\gamma$ is monotone.
$\s$
	\medskip

{\bf (2)} {\it If $i\in \Bbb Z$, then $\widetilde\psi\widetilde\gamma i\ \le \ 
i \le \widetilde\gamma\widetilde\psi i$.}
	\smallskip
Proof. We show that $|I(i)| \ge |P(\widetilde\psi i)|.$ 
Assume, for the contrary,  that $|I(i)| < |P(\widetilde\psi i)|$. Then 
$P(\widetilde\psi i)$ has a submodule $V$ of length $|I(i)|+1$. Note that the
socle of $V$ has to be $i$. 
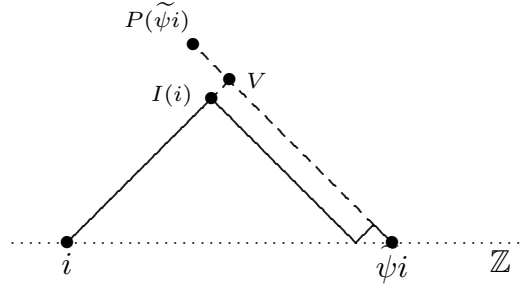
$$
{\beginpicture
    \setcoordinatesystem units <.2cm,.2cm>
\multiput{} at 0 0  15 11 /
\plot 0 0  8 8  16 0  17 1  18 0 /
\setdashes <1mm>
\plot 17 1  7 11 /
\plot  9 9  8 8 /
\put{$i$\strut} at 0 -1.3
\put{$\widetilde \psi i$\strut} at 18 -1.3
\put{$\ss V$} at 10.5 9
\put{$\ss I(i)$} at 5.8 8.2
\put{$\ss P(\widetilde\psi i)$} at 5 12.5 
\setdots <1mm>
\plot -3 0  26 0  /
\setsolid
\multiput{$\bullet$} at 0 0  8 8  9 9  7 11  18 0  /
\put{$\Bbb Z$} at 24 -1 
\endpicture}
$$
This implies that $I(i)$ is a submodule of $V$, and,
of course, a proper submodule. This is impossible, since $I(i)$ is injective.

Since $|I(i)| \ge |P(\widetilde\psi i)|,$ we see that  
$i \le i + |I(i)| - |P(\widetilde\psi i)| = \widetilde\psi (i) - |P(\widetilde\psi i)| =
\widetilde \gamma\widetilde\psi(i).$ 
In the same way, or using duality, one sees that $\widetilde\psi\widetilde\gamma i \le i$.
$\s$
	\medskip
Proof of Proposition A.1. The assertions (1) and (2) show that we can apply Lemma A.3 to 
the functions $f = \widetilde \psi$ and $g = \widetilde\gamma$. We get 
$\widetilde \psi^t\widetilde \gamma^t \widetilde \psi^t = \widetilde \psi^t$ as well as
$\widetilde \gamma^t\widetilde \psi^t \widetilde \gamma^t = \widetilde \gamma^t$, for all $t\ge 0.$ 
But $\widetilde \psi^t\widetilde \gamma^t \widetilde \psi^t = \widetilde \psi^t$ implies
that $\psi^t\gamma^t\psi^t = \psi^t$.
Similarly, $\widetilde \gamma^t\widetilde \psi^t \widetilde \gamma^t = \widetilde \gamma^t$ implies
that $ \gamma^t\psi^t\gamma^t = \gamma^t.$
$\s$
	\bigskip
{\bf A.5.} As a second application of A.3, we show Proposition A.2.
	\medskip 
{\bf (1)}
{\it The function $\widetilde\tau^-\widetilde\gamma\widetilde\tau$ 
is monotone.}
	\medskip
Proof. This follows directly from A.2, Lemma 1 and the fact that $\widetilde\tau$ and
$\widetilde\tau^-$ are, of course, monotone. $\s$
	\medskip
{\bf (2)} {\it If $i\in \Bbb Z$, then 
$\widetilde\tau^-\widetilde\gamma\widetilde\tau\widetilde\psi(i) \ \le\  i \ \le \  \widetilde\psi\widetilde\tau^-
\widetilde\gamma\widetilde\tau(i)$.}
	\medskip
Proof. First, let us show that $\widetilde\tau^-\widetilde\gamma\widetilde\tau\widetilde\psi(i)
\le i.$ Let $j = \widetilde\tau\widetilde\psi(i).$
Always, $I(i)$ is a factor module of $P(j)$ and $x = \widetilde\tau^-\gamma(j)$ is the socle
of $P(j).$ Thus $x \le i.$ 
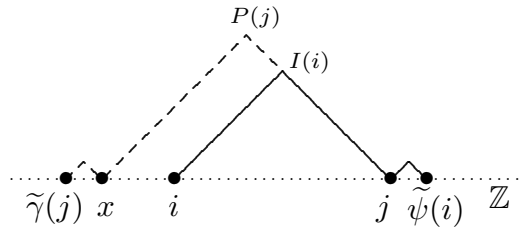
$$
{\beginpicture
    \setcoordinatesystem units <.2cm,.2cm>
\multiput{} at 0 0  15 11 /
\plot 6 0  12 6  18 0  19 1  20 0 /
\setdashes <1mm>
\plot 0 0  1 1  2 0  10 8  12 6  /
\put{$x$\strut} at 2.2 -1.6
\put{$i$\strut} at 6 -1.6
\put{$j$\strut} at 17.6 -1.6
\put{$\widetilde\psi(i)$\strut} at 20.5 -1.6
\put{$\widetilde\gamma(j)$\strut} at -.6 -1.6
\put{$\ss I(i)$} at 13.5 6.5
\put{$\ss P(j)$} at 10.5 9
\setdots <1mm>
\plot -3 0  26 0  /
\setsolid
\multiput{$\bullet$} at 0 0  2 0  6 0  18 0  20 0  /
\put{$\Bbb Z$} at 24 -1 
\endpicture}
$$

Second, we want to show that $i \ \le \  \widetilde\psi\widetilde\tau^-
\widetilde\gamma\widetilde\tau(i)$. We start with $i$ and $\widetilde\tau(i) = i\!-\!1$
and let $j = \widetilde\tau^-\widetilde\gamma(i\!-\!1)$. Note that $j = \soc P(i\!-\!1).$
It follows that $P(i\!-\!1) \subseteq I(j)$. If $P(i\!-\!1) = I(j)$, then
$i = \widetilde\psi(j).$ Otherwise, $P(i\!-\!1)$ is a proper submodule of $I(j)$
and then $i <  \widetilde\psi(j)$.
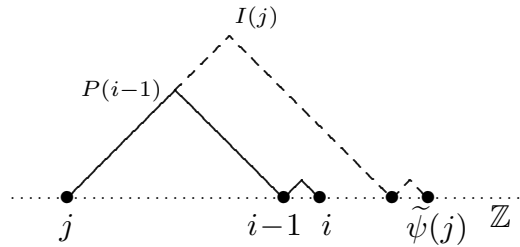
$$
{\beginpicture
    \setcoordinatesystem units <.2cm,.2cm>
\multiput{} at 0 0  15 11 /
\plot 0 0  6 6  12 0  13 1  14 0 /
\setdashes <1mm>
\plot 6 6  9 9  18 0  19 1  20 0  /
\put{$j$\strut} at 0 -1.6
\put{$i\!-\!1$\strut} at 11.5 -1.6
\put{$i$\strut} at 14.4 -1.6
\put{$\widetilde\psi(j)$\strut} at 20.5 -1.6
\put{$\ss P(i-1)$} at 3 6
\put{$\ss I(j)$} at 10.5 10
\setdots <1mm>
\plot -3 0  26 0  /
\setsolid
\multiput{$\bullet$} at 0 0  12 0  14 0  18 0  20 0  /
\put{$\Bbb Z$} at 24 -1 
\endpicture}
$$
$\s$
	\medskip
Proof of Proposition A.2. The functions $f = \widetilde\tau^-\widetilde\gamma\widetilde\tau$ and
$g = \widetilde\psi$ are monotone, according to A.5 (1) and A.4 (1), respectively.
According to A.5 (2), we see that the functions $f$ and $g$ satisfy the conditions of
Lemma A.3. Thus, for all $t\ge 0,$ we have  
$\psi^t\tau^-\gamma^t\tau\psi^t = \psi^t(\tau^-\gamma\tau)^t\psi^t = \psi^t$ and 
$\tau^-\gamma^t\tau\psi^t\tau^-\gamma^t\tau =
(\tau^-\gamma\tau)^t\psi^t(\tau^-\gamma\tau)^t = (\tau^-\gamma\tau)^t = \tau^-\gamma^t\tau$.
Multiplying the latter equality from the left by $\tau,$ from the right by $\tau^-$,
we get $\gamma^t\tau\psi^t\tau^-\gamma^t =
\gamma^t$, as required.
$\s$

	\bigskip
{\bf A.6. The case $t=1$.}
Note that $\Im \psi$ are just the simple modules
with projective dimension greater than $1$. Similarly, 
$\Im \gamma$ are just the simple modules
with injective dimension greater than $1$.

There is the bijection 
$$
 \gamma\: \Im \psi =\{T\mid \pd T \ge 2\} \longrightarrow \Im \gamma =\{S\mid \id S \ge 2\}
$$ 
with inverse $\psi$ :
$$
{\beginpicture
    \setcoordinatesystem units <.3cm,.3cm>
\multiput{} at 0 0  20 0  /
\plot 2 2  9 9  10 8  11 9  18 2 /
\setdots <1mm>
\plot 1 2  19 2 /
\multiput{$\bullet$} at 9 9  10 8  11 9  2 2  18 2 /
\setsolid 
\multiput{$\snake$} at 7.5 9  11.5 9 /
\put{$\ss IS$} at 9 10
\put{$\ss PT$} at 11 10 
\put{$S$} at 2 1 
\put{$T$} at 18 1 
\put{$\ss=\gamma T$} at 2 0 
\put{$\ss=\psi S$} at 18 0 
\setdots <.5mm> 
\plot 3 3  4 2  10 8  /
\setdashes <1mm> 
\plot  10 8  16 2  17 3 /
\endpicture}
$$ 
This bijection between $\Im \psi =\{T\mid \pd T \neq 1\}$ 
and $\Im \gamma =\{S\mid \id S \neq 1\}$  can also be seen also as bijections of 
$\Im\psi$ and $\Im \gamma$ with the set of ``valleys'' of the ``roof'' of $A$ (see [R2]);
by definition, a module $V$ is called a {\it valley} provided $V$ 
is both the radical of a projective
module as well of the form $IS/S$, where $S$ is simple. Corresponding to the valleys in the
roof, there are the peaks (the indecomposable modules which are both
projective and injective).

There is also the bijection 
$$
 \gamma\tau\: \Im \psi =\{T\mid \pd T \ge 2\} \longrightarrow \Im \gamma =\{S\mid \id S \ge 2\}
$$ 
with inverse $\psi\tau^-$ 
$$
{\beginpicture
    \setcoordinatesystem units <.3cm,.3cm>
\multiput{} at 0 0  20 0  /
\plot 2 2  3 3  4 2  10 8  16 2  17 3  18 2  /
\setdots <1mm>
\plot 1 2  19 2 /
\put{$\star$} at 10 8 
\multiput{$\bullet$} at     2 2  18 2 /
\setsolid 
\put{$S$} at 2 1 
\put{$T$} at 18 1 
\setdots <.5mm> 
\plot 3 3  4 2  10 8  /
\setdashes <1mm> 
\plot  10 8  16 2  17 3 /
\put{$\ss=\gamma\tau T$} at 2 0 
\put{$\ss=\psi\tau^- S$} at 18 0 
\put{$\ss \tau T$} at 16 1.2 
\put{$\ss \tau^- S$} at 4 1.2 
\put{$\ss I\tau^- S = P\tau T$} at 10 9 
\setdashes <1mm>
\plot 3 3  7 7  8 6 /
\plot 17 3  13 7  12 6 /
\endpicture}
$$ 
Here, the middle module $I\tau^-S = P\tau T$ is a peak and 
$\gamma\tau$ provides bijections between $\{T\mid \pd T \ge 2\}$ and the set of peaks
as well as between the set of peaks and $\{S\mid \id S \ge 2\}$.

Note that all these sets 
$\{T\mid \pd T \ge 2\}$, $\{S\mid \id S \ge 2\}$, the set of valleys and 
the set of peaks have the same cardinality as a minimal set of
admissible relations which are needed to define $A$, see Appendix C.
	\medskip
{\bf A.7. Historical remark.} Some of the effects of the interrelation 
between $\psi$ and $\gamma$ were already discussed by Shen [Sh3]. 
The assertion that there are as many $\psi$-cyclic simple modules as there are 
$\gamma$-cyclic simple modules was shown by Shen in [Sh2]. 
	\bigskip
{\bf A.8. The $\psi$-quiver and the $\gamma$-quiver of a cyclic Nakayama algebra $A$.}
	\medskip 
Using $\gamma$, one can define the {\it resolution quiver} (or {\it $\gamma$-quiver})
of $A$, see [R2]: its vertices
are the isomorphism classes $[S]$ of the simple modules $S$, and there is an arrow
$[S] \to [T]$ provided $T = \gamma(S)$. Of course, the $\gamma$-paths are just the paths 
in the resolution quiver. 

Dually, using $\psi$, one can define the {\it coresolution quiver} 
(or {\it $\psi$-quiver}) of $A$: its vertices
are the isomorphism classes $[S]$ of the simple modules $S$, and there is an arrow
$[S] \to [T]$ provided $T = \psi(S)$. The $\psi$-paths are just the paths 
in the coresolution quiver. 

Note that the coresolution quiver
of $A$ is just the resolution quiver of $A^{\op}$. 
	\bigskip
{\bf Warning.} The $\psi$-quiver and the $\gamma$-quiver of a Nakayama algebra $A$
have many properties in common. Several such properties are mentioned above. In addition,
Shen [Sh1] has shown, that also the number of components are equal.
However, the cardinalities of the components of the $\gamma$-quiver 
may be different from the cardinalities of the components of the $\psi$-quiver 
as the following example shows:
$$
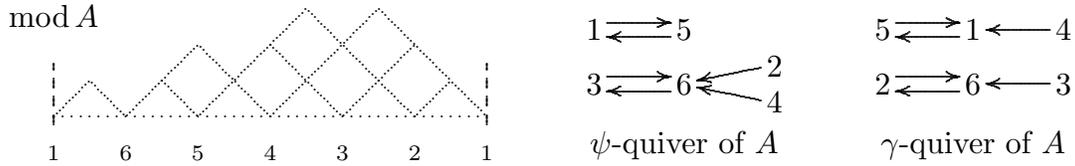
{\beginpicture
    \setcoordinatesystem units <.4cm,.4cm>
\put{\beginpicture
\put{$\mod A$} at 0 2.8
\multiput{} at 0 0  12 3    /
\setdots <1mm>
\plot 0 0  12 0 /
\setdots <.5mm>
\plot  0 0  1 1  2 0  4 2  6 0  9 3  12 0  /
\plot  3 1  4 0  7 3  10 0  11 1 /
\plot  6 2  8 0  10 2 /
\multiput{$\ss 1$} at 0 -1  12 -1 /
\put{$\ss 2$} at 10 -1
\put{$\ss 3$} at 8 -1
\put{$\ss 4$} at 6 -1
\put{$\ss 5$} at 4 -1
\put{$\ss 6$} at 2 -1
\setdashes <1mm>
\plot 0 -.2  0 1.5 /
\plot 12 -.2  12 1.5 /
\endpicture} at 0 0
\put{\beginpicture
    \setcoordinatesystem units <.5cm,.4cm>
\put{$1$} at 0 1.5
\put{$2$} at 4 0.5
\put{$3$} at 0 0
\put{$4$} at 4 -.5
\put{$5$} at 2 1.5
\put{$6$} at 2 0
\arr{0.3 1.7}{1.7 1.7}
\arr{1.7 1.3}{0.3 1.3}
\arr{0.3 0.2}{1.7 0.2}
\arr{1.7 -.2}{0.3 -.2}
\arr{3.7 0.45}{2.3 0.1}
\arr{3.7 -.45}{2.3 -.1}
\put{$\psi$-quiver of $A$} at 2 -1.7
\endpicture} at 12 -.2
\put{\beginpicture
    \setcoordinatesystem units <.5cm,.4cm>
\put{$5$} at 0 1.5
\put{$2$} at 0 0
\put{$3$} at 4 0
\put{$4$} at 4 1.5
\put{$1$} at 2 1.5
\put{$6$} at 2 0
\arr{0.3 1.7}{1.7 1.7}
\arr{1.7 1.3}{0.3 1.3}
\arr{0.3 0.2}{1.7 0.2}
\arr{1.7 -.2}{0.3 -.2}
\arr{3.7 0}{2.3 0}
\arr{3.7 1.5}{2.3 1.5}
\put{$\gamma$-quiver of $A$} at 2 -1.7
\endpicture} at 20 -.2
\endpicture}
$$
	\medskip
{\bf Remark.} The $\psi$-quiver (and similarly, 
the $\gamma$-quiver) of a cyclic Nakayama algebra $A$ do not 
determine the finitistic dimension of $A$. Namely, the following 
two algebras $A, A'$ have the same $\psi$-quiver, but $\finpro A = 2,$ whereas
$\finpro A' = 1.$ 
$$
{\beginpicture
    \setcoordinatesystem units <.4cm,.4cm>
\put{\beginpicture
\put{$\mod A$} at 0 3.6
\multiput{} at 0 0  4 2 /
\setdots <1mm>
\plot 0 0  4 0 /
\setdots <.5mm>
\plot  0 0  2 2  4 0 /
\plot 1 1  2 0  3 1 /
\multiput{$0$} at 1 1  2 2  /
\multiput{$1$} at 0 0  4 0 /
\multiput{$2$} at 2 0  3 1  /
\put{$\ss 1$} at 4 -1
\put{$\ss 2$} at 2 -1
\setdashes <1mm>
\plot 0 -.2  0 1.5 /
\plot 4 -.2  4 1.5 /
\endpicture} at 0 0.3 
\put{\beginpicture
\put{$\mod A'$} at 0 3.2
\multiput{} at 0 0  4 2 /
\setdots <1mm>
\plot 0 0  4 0 /
\setdots <.5mm>
\plot  0 0  2 2  4 0 /
\plot 1 1  2 0  3 1 /
\plot 0 0  1 -1  2 0  3 -1  4 0  5 -1 /
\multiput{$0$} at 1 1  2 2  /
\multiput{$1$} at 1 -1  5 -1 /
\multiput{$\infty$} at 2 0  3 1  0 0  3 -1 4 0  /
\put{$\ss 1$} at 5 -2
\put{$\ss 2$} at 3 -2
\setdashes <1mm>
\plot 0 -1.2  0 1.5 /
\plot 4 -1.2  4 1.5 /
\endpicture} at 8 0 
\put{\beginpicture
\put{$\psi$-quiver} at  0 3
\put{$1$} at 0 0
\put{$2$} at 2 0
\arr{0.3 0}{1.7 0}
\circulararc 300 degrees from 2.0 -.5  center at 3 0
\arr{2.07 .6}{2 .5}
\endpicture}  at 15 1
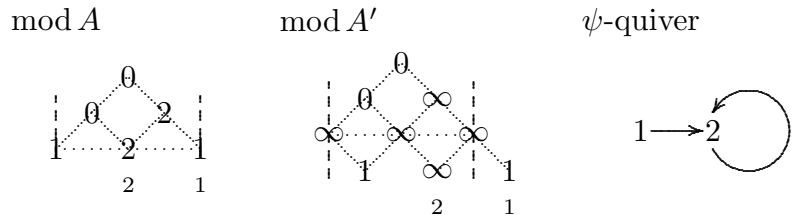
\endpicture}
$$
	\bigskip
{\bf Appendix B. Sen's $\epsilon$-construction.}
	\medskip
It seems to be worthwhile to add a short outline of a reduction method which was 
introduced and very skillfully used by Sen [S1, S2, S3],
since it relates to the functions $\gamma$ and $\psi$ (discussed in sections 3 and 4, as well 
as in Appendix A). 

Let $A$ be a cyclic Nakayama algebra.
Sen considers the subcategory $\Cal F$ of modules $M$ with what he calls ``syzygy filtrations'':
this is an abbreviation of saying that $M$ has a filtration whose factors are {\bf second} syzygies
of simple modules (denoted below by $\Delta T$).   
	\medskip 

{\bf B.1.} Three sets of simple modules will play a role: $\Cal S$, $\Cal T$ and $\Cal U$
(see also Appendix C). 
Here, $\Cal S$ is the set of simple torsionless modules. Note that an indecomposable module $M$
is torsionless if and only if $\soc M$ belongs to $\Cal S$.
Then, $\Cal T$ is the set of simple modules $T$ with $\id T \ge 2;$ these are the simple
modules in $\Im \gamma$. Finally, $\Cal U$ is the set of simple modules $U$ with $\id U \ge 2;$ 
these are the simple modules in $\Im \psi$. 
{\it A simple module $S$ belongs to $\Cal S$ if and only if $\tau S$ belongs to $\Cal T.$}
	\medskip
Let $\Cal F$ be the full subcategory of all modules with top in $\add \Cal T$ and socle in
$\add \Cal S$. Thus, an indecomposable module $M$ belongs to $\Cal F$ if and only if $\top M$
belongs to $\Cal T$ and $M$ is torsionless.
	\medskip
{\bf B.2. The modules $\Delta T$.} 
For $T\in \Cal T$, let $\Delta T$ be the smallest non-zero factor module of $PT$ which is torsionless
(since $PT$ is torsionless, $\Delta T$ exists; also, 
$\Delta T$ is, of course, uniquely determined by $T$).
	\medskip
(a) {\it The module $\Delta T$ has a unique composition factor which belongs to $\Cal T$,
namely $T = \top\Delta T$. Similarly, $\Delta T$ has a unique composition factor which belongs to 
$\Cal S$, namely $\tau^{-}T = \soc\Delta T$.}
	\smallskip
(b) {\it If $C$ is any simple module, there is a unique $T\in \Cal T$ such that $C$ is a composition factor of $\Delta T$ and the multiplicity of $C$ in $\Delta T$ is $1$.} (Thus, every simple module occurs with
multiplicity 1 as a composition factor of $\bigoplus_{T\in \Cal T} \Delta T$.)

	\smallskip
(c)  {\it If $T \in \Cal T,$ then 
$$
 \Delta T = \Omega^2\psi T,
$$
thus $\Delta T = \Omega^2U$, where $U = \psi T\in \Cal U$. 
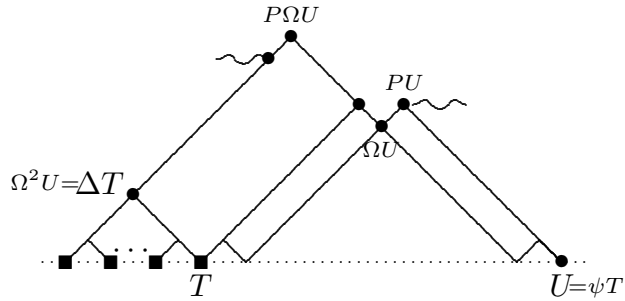
$$
{\beginpicture
    \setcoordinatesystem units <.25cm,.25cm>
\multiput{} at 0 0  22 10 /
\plot 0 0  10 10  20 0  21 1  22 0 /
\plot 3 3  6 0  13 7 /
\plot 1 1  2 0  /
\plot 4 0  5 1 /
\plot 7 1  8 0  15 7  22 0  /
\put{$\cdots$} at 3 .5
\multiput{$\snake$} at 7.3 9  16 7 /
\setdots <1mm>
\plot -1 0  23 0 /
\multiput{$\bullet$} at 3 3  9 9  10 10  13 7  14 6  15 7  22 0 /
\put{$\ss PU$} at 15 8 
\put{$\Delta T$} at 1.5 3.5
\put{$\ss\Omega^2 U=$} at -.9 3.65
\put{$T$} at 6 -1
\put{$U$} at 22 -1
\put{$\ss =\psi T$} at 23.5 -1.1
\put{$\ss\Omega U$} at 14 5 
\put{$\ss P\Omega U$} at 10 11
\multiput{$\sssize\blacksquare$} at 0 0  2 0  4 0  6 0 / 
\endpicture}
$$

Conversely, if $U$ is simple, then
either $\pd U \le 1$ and $\Omega^2 U = 0$, or else $U \in \Cal U,$ and then 
$T = \gamma U$ belongs to $\Cal T$ and $\psi T = U.$ 
It follows that the modules $\Delta T$ are just the non-zero modules which are second syzygy modules
of simple modules.}

	\medskip
Note that the relevance of the modules $\Omega^2 U$ mentioned in (c) was stressed
already by Madsen [M].  
	\medskip
Proof. (a) First, we show that 
$\top \Delta T$ is the only composition factor of $\Delta T$ which belongs to $\Cal T$.
Namely, assume, for the contrary, that there are submodules 
$X'\subset X \subset \Delta T$ with $X/X'\in \Cal T$.
Then $\soc (\Delta T/X) \simeq \tau^-\top X = \tau^-(X/X')$ is in $\Cal S$. This implies that 
$\Delta T/X$ is torsionless. But this contradicts the minimality of $\Delta T.$

Next, we show that $\soc \Delta T$ is the only composition factor of $\Delta T$ which 
belongs to $\Cal S$. 
If $0 \subset X' \subset X \subseteq \Delta T$ are submodules such that $X/X'$ is in $\Cal S$,
then $\top X' \simeq \tau \soc (\Delta T/X') = \tau (X/X')$ 
belongs to $\Cal T$. But this contradicts the first assertion. 

(b) Let $S$ be simple. Let $c\ge 1$ be minimal with $\tau^{-c}S$ torsionless. Thus $T = \tau^{-c+1}S =
\tau(\tau^{-c} S)$ belongs to $\Cal T$, whereas the modules $\tau^iT$ with $1\le i < c$ do not belong to
$\Cal T.$ It follows that $\Delta T$ has length at least $c$ and that $S$ is a composition factor of $\Delta T$. Since the top $T$ of $\Delta T$ has Jordan-H\"older multiplicity 1, any composition factor 
of $\Delta T$ has multiplicity 1.  

(c) Let $T \in \Cal T$ and $U = \psi T.$ Then $\pd U \ge 2$ and $\gamma U = T$. 
Since $\pd U \ge 2$, it follows from 3.2 (b) that $\top \Omega^2 U = \gamma U = T.$
Of course, $\Omega^2 U$ is torsionless, thus
we see that $\Delta T$ is a factor module of $\Omega^2 U.$ According to 3.2 (c), only the socle
of $\Omega^2 U$ can be torsionless, thus $\Delta T = \Omega^2 U.$

Conversely, let $U$ be simple. If $\pd U \le 1$, then $\Omega^2 U = 0.$ Thus, let $\pd U \ge 2$ 
and $T = \gamma U$. Thus, $T\in \Cal T$ and $T = \gamma U.$ 
$\s$
	\bigskip
{\bf B.3. The category $\Cal F.$} {\it The modules $\Delta T$ with $T\in \Cal T$ are pairwise orthogonal bricks and $\Cal F$ is the full subcategory of all modules which have a filtration
with factors of the form $\Delta T$ with $T\in \Cal T.$ 

The category $\Cal F$ is an extension-closed exact abelian category with simple objects the modules
$\Delta T$ with $T\in \Cal T.$ The category $\Cal F$ is a serial length category with
$$
  \tau_{\Cal F}\Delta T = \Delta (\tau\soc \Delta T),
$$
provided $\Delta T$ is not projective.} (Here, $\tau_{\Cal F}$ is the Auslander-Reiten
translation inside the category $\Cal F$.)
	\medskip
Proof. Let $T \in \Cal T$. If $f\:\Delta T \to \Delta T$ is a non-zero non-invertible endomorphism, 
then the image of $f$ is a proper submodule of $\Delta T$ with top $T$, a contradiction.
Also, if $T', T$ are non-isomorphic modules in $\Cal T$, and $f\:\Delta T'\to \Delta T$
is a non-zero homomorphism, then the image of $f$ is a submodule of $\Delta T$ with top $T'$,
again a contradiction. Thus, we see that the modules $\Delta T$ with $T\in \Cal T$ are pairwise
orthogonal bricks. 

Next, we show that an indecomposable module $M$ with top in $\Cal T$ and socle in $\Cal S$ has a filtration with factors of the form $\Delta T$ with $T\in \Cal T$.
The proof is by induction on the length of $M$. 
Let $M$ be a module with top in $\Cal T$ and socle in $\Cal S$.
Say let $\top M = T'$. Then $M$ is a factor module of $PT'$. Since $\soc M$ is torsionless, 
also $M$ is torsionless. By definition of $\Delta T'$ we see that $\Delta T'$ is a factor module of
$PT'$, say $M = PT'/X$ for some submodule $X$ of $PT'$. Now $\top X = \tau\soc (PT'/X) = \tau\soc M$,
thus $\top X$ is in $\Cal T.$ By induction, $X$ has a filtration with
factors of the form $\Delta T$. $\s$

It follows that 
$\Cal F$ is the extension closure of the class of modules $\Delta T$ with
$T\in \Cal T$, and that $\Cal F$ is an extension-closed exact abelian subcategory, whose
simple objects are just the objects $\Delta T$ with $T\in \Cal T$, see for example [R1].

It remains to calculate $\tau_{\Cal F}\Delta T.$
Let $\pi\:PT \to \Delta T$ be a projective cover. If $\Delta T$ is not
projective, then $\Ker(f)$ is an indecomposable module, 
which belongs to $\Cal F$ and has $\Delta T'$ with $T' = \tau \soc\Delta T$ 
as a factor object (its top in $\Cal F$). Thus $\tau_{\Cal F}\Delta T = \Delta(T').$
$\s$
	\bigskip
Note: The algebra $\epsilon(A)$ is not necessarily cyclic, nor even connected.
For example, for 
$$
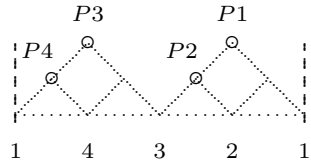
{\beginpicture
    \setcoordinatesystem units <.4cm,.4cm>
\put{\beginpicture
\multiput{} at 0 -1  8 3.5 /
\setdots <1mm>
\plot 0 0  8 0 /
\setdots <.5mm>
\plot  0 0  2 2  4 0   6 2  8 0  /
\plot  1 1  2 0  3 1  /
\plot  5 1  6 0  7 1 /
\multiput{$\ss 1$} at  0 -1  8 -1 /
\multiput{$\ss 2$} at 6 -1  /
\multiput{$\ss 3$} at 4 -1  /
\multiput{$\ss 4$} at 2 -1  /
\multiput{$\circ$} at 1 1  2 2  5 1  6 2 /
\put{$\ss P1$} at 6 2.8
\put{$\ss P2$} at 4.6 1.8 
\put{$\ss P3$} at 2 2.8
\put{$\ss P4$} at 0.6 1.8 
\setdashes <1mm>
\plot 0 -.2  0 2 /
\plot 8 -.2  8 2 /
\endpicture} at 0 0 
\endpicture}
$$
the set $\Cal F$
consists of the projective modules $P2$ and $P4$. Thus, $\Cal F$ is a semisimple category
and its quiver consists of two singletons (and $\epsilon(A) = k\times k$).
	\bigskip
{\it The subcategory $\Cal F$ of $\mod A$ is closed under projective covers.} It follows
that the abelian category $\Cal F$ has enough projective objects. The indecomposable projective
objects in $\Cal F$ are the modules $PT$ with $T\in \Cal T.$ 
	\medskip
Proof that $\Cal F$ is closed under
projective covers: if $M$ is indecomposable and in $\Cal F$, say with top $T$, then
$PM = PT$ also belongs to $\Cal F$. $\s$ 

	\bigskip
{\it The module
$G = \bigoplus_{T\in \Cal T} PT$ is a progenerator of $\Cal F.$ It follows that 
$\epsilon(A) = (\End G)^{\op}$ is a Nakayama algebra and
there is a categorical equivalence}
$$
 \Cal F \simeq \mod \epsilon(A).
$$  
\vglue-1cm $\s$
	\bigskip
{\bf B.4.  The $\gamma$-quiver of $\epsilon(A).$} We assume that none
of the modules $\Delta T$ is projective, thus $\epsilon(A)$ is again a {\bf cyclic}
Nakayama algebra. Using the categorical equivalence 
$\Cal F \simeq \mod \epsilon(A)$, the $\gamma$-quiver of $\epsilon(A)$ may be
considered as a quiver with vertex set the set of modules $\Delta T$ with $T\in \Cal T$
and with arrows $\Delta T \to \gamma_{\Cal F}(\Delta T).$ There is the following observation:
	\medskip
{\bf Lemma.} {\it If $T\in \Cal T$, then }
$$
 \gamma_{\Cal F}(\Delta T) = \Delta(\gamma T).
$$
	\medskip
Proof. Let $\soc_{\Cal F}P_{\Cal F}(\Delta T) = \Delta(T')$ for some $T'$. Then
$\soc \Delta(T') = \soc PT$. Since $\gamma_{\Cal F} = \tau_{\Cal F}\soc_{\Cal F}P_{\Cal F}(-)$
and $\gamma = \tau\soc P(-)$, we see that
$$
\align 
 \gamma_{\Cal F}(\Delta T) &= \tau_{\Cal F}\soc_{\Cal F}P_{\Cal F}(\Delta T) \cr
 &= 
 \tau_{\Cal F}\Delta(T') = 
 \Delta(\tau \soc\Delta(T')) \cr
 &= \Delta(\tau \soc PT) = \Delta(\gamma T).
\endalign
$$
\vglue-1cm
$\s$
	\bigskip
This shows: {\it Let $\Gamma(A)$ be the $\gamma$-quiver of $A$ and $\Gamma(A)|\Im\gamma$
its restriction to $\Im\gamma$. Then 
the map $T \mapsto \Delta(T)$ identifies $\Gamma(A)|\Im\gamma$ 
with the $\gamma_{\Cal F}$-quiver of $\epsilon(A).$ }
	\bigskip

{\bf B.5. The subcategory $\Cal F$ and the class of reflexive modules.} 
	\medskip
{\bf Proposition.} 
{\it The modules in $\Cal F$ are reflexive.}
	\medskip
Proof. Assume that $M$ belongs to $\Cal F$. We want to show that $M$ is reflexive.
We can assume that $M$ is indecomposable and not projective. Now $M$ is torsionless,
thus $IM$ is projective. Let $0 \neq P \subseteq IM$ be minimal projective (this means: projective,
and $\rad P$ is not projective). Since $M$ is a non-projective submodule of $IM$, we see that
$M$ is a proper submodule of $P$. It follows that the inclusion map $u\:M \to P$ is a minimal
left $\add{}_AA$-approximation, therefore $\mho M = \Cok u$. Now $\top M$ belongs to $\Cal T$,
therefore $\soc \Cok u = \tau^-\top M$ belongs to $\Cal S$, thus $\Cok u$ 
is torsionless. 
It is well-known that a module
$M$ is reflexive if and only if both $M$ and $\mho M$ are torsionless (see for example
[AZ], Theorem 1.5). 
 $\s$
	\medskip
Let $\Cal R$ be the full subcategory of all reflexive modules. Let $\Cal R_0$ be the
full subcategory of all reduced reflexive modules (a module $M$ is said to be {\it reduced}
provided $M$ has no non-zero projective direct summand). 
	\medskip
{\bf Corollary.} {\it We have }
$$
 \Cal R_0 \subseteq 
  \{\Omega^2X\mid X\in \mod A\} \subseteq \Cal F \subseteq 
  \{P\oplus \Omega^2X \mid X\in \mod A,\ P\in \add{}_AA \} = \Cal R.
$$

Proof. First, let $M$ be indecomposable, reflexive and non-projective.
If $X = \mho^2 M$, 
then $M = \Omega^2X$. This shows that $\Cal R_0 \subseteq 
  \{\Omega^2X\mid X\in \mod A\}$.

Next, let us show that 
any module of the form $\Omega^2 M$ belongs to $\Cal F.$
On the one hand, $\Omega^2 M$ is torsionless. On the other hand, either $\Omega^2 M = 0$,
or else $\top \Omega^2 M = \gamma \top M$, thus $\top \Omega^2M \in \Cal T.$ This shows that
$M\in \Cal F$. 

Second, let $M \in \Cal F$. We want to show that $M = P\oplus \Omega^2X$ for some modules
$P,X$ with $P$ projective. We can assume that $M$ is indecomposable and not projective.
As we saw, $M$ is reflexive. Thus, $M = \Omega^2 X$ for some module $X$. $\s$
	\medskip
{\bf Examples.} {\it All three inclusions mentioned in the Corollary may be proper.}
	\medskip
$$
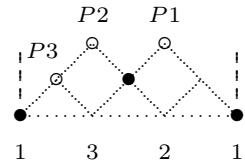
{\beginpicture
    \setcoordinatesystem units <.4cm,.4cm>
\put{\beginpicture
\put{Example 1} at -3 3 
\multiput{} at 0 -1  6 3.5 /
\setdots <1mm>
\plot 0 0  4 0 /
\setdots <.5mm>
\plot  0 0  2 2  4 0   /
\plot  1 1  2 0  3 1   /
\multiput{$\ss 1$} at  0 -1  4 -1 /
\multiput{$\ss 2$} at 2 -1  /
\multiput{$\bullet$} at 1 1  /
\put{$\ss P1$} at 2 2.8
\put{$\ss P2$} at 0.6 1.8 
\setdashes <1mm>
\plot 0 -.2  0 2 /
\plot 4 -.2  4 2 /
\endpicture} at 0 0 
\put{\beginpicture
\put{Example 2} at -3 3 
\multiput{} at 0 -1  6 3.5 /
\setdots <1mm>
\plot 0 0  6 0 /
\setdots <.5mm>
\plot  0 0  2 2  4 0  5 1  5 1  /
\plot  1 1  2 0  4 2  6 0  /
\multiput{$\ss 1$} at  0 -1  6 -1 /
\multiput{$\ss 2$} at 4 -1  /
\multiput{$\ss 3$} at 2 -1 /
\multiput{$\bullet$} at 0 0  3 1  6 0 /
\multiput{$\circ$} at 1 1  2 2  4 2  /
\put{$\ss P1$} at 4 2.8
\put{$\ss P2$} at 2 2.8
\put{$\ss P3$} at 0.6 1.8 
\setdashes <1mm>
\plot 0 0.5  0 2 /
\plot 6 0.5  6 2 /
\endpicture} at 12 0 

\endpicture}
$$
	\medskip
Example 1: Here is $\Cal R_0$ the zero category, whereas there is
an indecomposable module of the form $\Omega^2 X$, namely the projective $P2.$ 

In example 2, 
the indecomposable modules of the form $\Omega^2 X$ are just $1$ and $\rad P1$ (marked
by bullets). The subcategory $\Cal F$ contains also $P1$ and $P2$, but not $P3.$
Note that $\Cal F$ is closed under extensions, whereas the subcategories
$\{\Omega^2X\mid X\in \mod A\}$ and $\{P\oplus \Omega^2X \mid X\in \mod A,\ P\in \add{}_AA \}$
are not closed under extensions (the extension closure of 
$\{\Omega^2X\mid X\in \mod A\}$ is $\Cal F$, and the extension closure of 
$\{P\oplus \Omega^2X \mid X\in \mod A,\ P\in \add{}_AA \}$ contains $3\oplus P2.$ 
	\medskip
{\bf Warning.} Some mathematical papers (for example 
Auslander-Reiten [AR1], [AR2]) denote by $\Omega^t(\mod A)$ the full subcategory of all modules 
of the form $P\oplus \Omega^t X$, with $P$
projective (and not, what one would expect, the full subcategory $\Cal C$ of all modules of the form
$\Omega^tX$, where $X$ is a module, or the additive closure of $\Cal C$)
and consider the question whether this subcategory $\Omega^t(\mod A)$ is
closed under extension or not. The examples above show that subcategories related to 
this subcategory may be closed under extensions whereas $\Omega^t(\mod A)$ itself is not.

	\bigskip
{\bf Appendix C. Some canonical bijections.}
	\medskip
As an addition to appendix B, we assume again that $A$ is a cyclic Nakayama algebra.
We want to focus the attention to several further sets of modules which are in natural bijection
to the subsets $\Cal S, \Cal T, \Cal U$ mentioned above. The category $\Cal F$ is defined by using
as building blocks the modules $\Omega^2 U$ with $U\in \Cal U$, thus one may draw the attention also
to the modules $\Omega U$ with $U \in \Cal U,$ we call them ``valley'' modules.

Recall from A.6 that an indecomposable module $V$ is a valley provided that
$V = \rad P = I/\soc I$, where $P$ is projective and $I$ is injective. Also,
an indecomposable module $Z$ is a peak provided that 
$Z$ is both projective and injective.
(Looking at the roof of the Auslander-Reiten quiver of $A$, we see that valleys and peaks alternate.
Of course.)

An indecomposable projective module $P$ is said to be {\it minimal} provided that
$\rad P$ is not
projective.
An indecomposable injective module $I$ is said to be {\it minimal} provided that 
$I/\soc I$ is not
injective.

$$ 
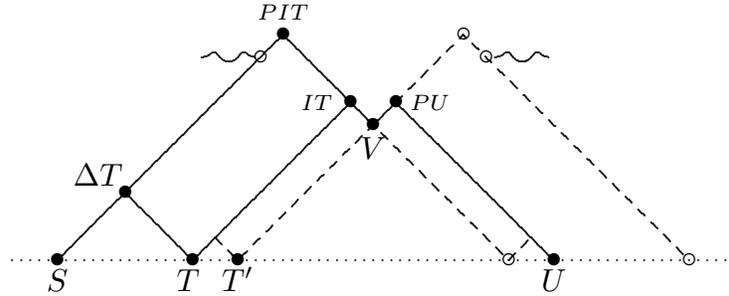
{\beginpicture
    \setcoordinatesystem units <.25cm,.25cm>
\multiput{} at  0 0  19 10 /
\multiput{$\bullet$} at 0 0  3 3  6 0    10 10  13 7  14 6  15 7  22 0  8 0 /
\multiput{$\circ$} at   9 9  20 0  28 0  18 10  19 9  /
\setdots <1mm>
\plot -2 0  30 0 /
\setsolid
\plot 0 0  10 10  14 6  15 7  22 0 /
\plot 3 3  6 0  13 7 /
\multiput{$\snake$} at 7 9  20 9 /
\setdashes <1mm>
\plot 7 1  8 0  14 6  20 0  21 1 /
\plot 15 7  18 10  28 0   /
\put{$T$\strut} at 5.8 -1
\put{$T'$\strut} at 8 -1
\put{$U$\strut} at 22 -1
\put{$\ss IT$\strut} at  11.5 7 
\put{$\ss PIT$\strut} at 10 11 
\put{$\Delta T$} at 1.8 3.8
\put{$V$} at 14  5
\put{$\ss PU$} at 16.5 7
\put{$S$\strut} at 0 -1
\endpicture}
$$
	
{\bf Some relevant bijections.} 
We consider the modules marked above by bullets: these are modules which 
determine each other uniquely. Thus, we look at 
the set $\Cal T$ of simple modules $T$ with $\id T \ge 2;$
the set $\Cal U$ of simple modules $U$ with $\pd U \ge 2;$ the set $\Cal S$
of simple torsionless modules;
the set of modules $\Delta T$ with $T\in \Cal T;$ the set of valley modules, the set of peak
modules, the set of minimal injective modules and finally the set of minimal projective modules
(the sets are arranged roughly in the same way as the modules are located in the
Auslander-Reiten quiver):
{\it All these sets have the same cardinality, and any arrow in 
the following pictures exhibits a canonical
bijection.} 
Actually, the set $\Cal S$ of simple torsionless modules appears twice in the picture,
since two different bijections between $\Cal S$ and $\Cal T$ are relevant in our
setting: if $T\in \Cal T,$ then, on one hand, the socle of $\Delta T$ is torsionless and simple, 
on the other hand, also $\tau^-T$ is torsionless and simple. These bijections are combined 
when we look at the Auslander-Reiten translation $\tau_{\Cal F}$ of $\Cal F$, see (4).

$$ 
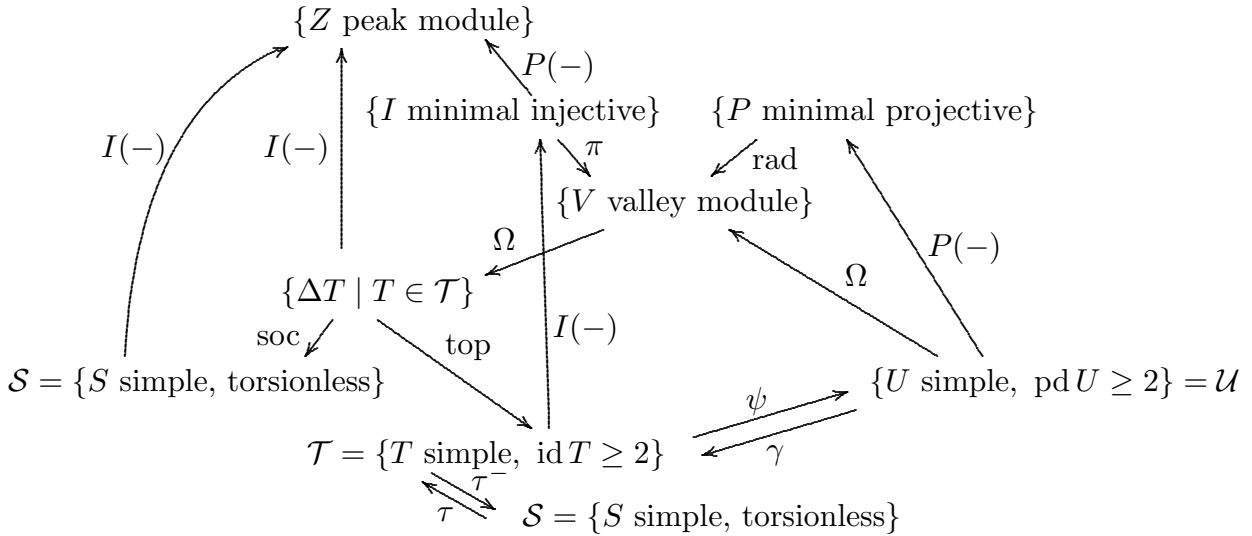
{\beginpicture
    \setcoordinatesystem units <1cm,1cm>
\put{$\left\{ \text{$Z$ peak module} \right \}$} at -2.8 4
\put{$\left\{ \text{$I$ minimal injective} \right \}$} at -1.7 3
\put{$\left\{ \text{$P$ minimal projective} \right \}$} at 2.3 3
\put{$\left\{ \text{$V$ valley module} \right \}$} at 0.2 2
\put{$\left\{ \Delta T\mid T\in \Cal T \right\}$} at -3.2 1
\put{$\Cal T = \left\{ T\ \text{simple,\ } \id T \ge 2 \right\}$} at -2 -.8
\put{$\Cal S = \left\{ S\ \text{simple, torsionless}  \right\}$} at .5 -1.5
\put{$\Cal S = \left\{ S \text{\ simple, torsionless} \right\}$} at -5.2 0
\put{$\left\{ U \text{\ simple,\ } \pd U \ge 2 \right\} = \Cal U$} at 4.3 0
\arr{.3 -.6}{2 -.1}
\arr{2.1 -.3}{.4 -.8}
\put{$\psi$} at 1 -.2
\put{$\gamma$} at 1.2 -.8
\put{$\tau^-$} at -1.95 -1.05
\put{$\tau$} at -2.45 -1.45

\arr{-2 -1.5}{-2.7 -1.1}
\arr{-2.6 -}{-1.9 -1.4}

\arr{-3.7 .7}{-4 .3}
\arr{-3.2 .7}{-1.5 -.5}
\put{$\soc$} at -4.3 .5
\put{$\top$} at -2.2 0.4

\arr{-.7 1.7}{-2 1.2}
\put{$\Omega$} at -1.8 1.6
\arr{3 .3}{.7 1.7}
\put{$\Omega$} at 2.1 1.2

\arr{3.5 .3}{2 2.7}
\put{$P(-)$} at 3.3 1.5
\arr{1 2.7}{.5 2.3}
\put{$\rad$} at 1.2 2.5

\arr{-3.6 1.5}{-3.6 3.7}

\arr{-1.3 -.5}{-1.4 2.7}
\put{$I(-)$} at -.9 0.6
\put{$I(-)$} at -4.1 2.6

\arr{-1.5 3.2}{-2 3.8}
\put{$P(-)$} at -1.2 3.5

\arr{-1.2 2.7}{-.85 2.3}
\put{$\pi$} at -.8 2.6

\setquadratic
\plot -6 .3  -5.5 2.5  -4.3 3.7 /
\arr{-4.3 3.68}{-4.2 3.72}
\put{$I(-)$} at -5.9 2.6
\endpicture}
$$
Here, $\pi(I) = I/\soc I$ for $I$ injective.
	\medskip
The display shows the interrelation between a large number of sets of cardinality $r$,
where $|\Cal T| = r.$
However one should be aware that even this display is not yet complete. For example,
similar to the torsionless simple modules, also the divisible simple modules play a role.

And, parallel to $\Delta T = \Omega V$, where $V$ is a valley, also the modules $\Sigma V$ 
should be taken into account. But we have refrained from doing so in order to keep the display manageable. 

In general, if $\Cal W$ is a set of simple modules and has
cardinality $r$ (such as $\Cal S,\ \Cal T,\ \Cal U$ or the set $\Cal D$ of divisible modules),
also the sets $\{PW\mid W\in \Cal W\}$, $\{\rad PW\mid W\in \Cal W\}$,
$\{IW\mid W\in \Cal W\}$ and $\{IW/\soc \mid W\in \Cal W\}$ are sets of cardinality $r$.
	\bigskip 
{\bf Minimal sets of relations.} Let $r$ be the cardinality of $\Cal T$. 
When looking at the large number of sets of modules
which are in canonical bijection with $\Cal T$, thus of cardinality $r$, 
one should be aware that it is 
another set of cardinality $r$ which really is the decisive one, namely
any minimal set of relations which defines 
$A$. This is Sen's approach and it is definitely the relevant one. It leads immediately
to an intrinsic choice of indices. Unfortunately, it seems that the 
use of these indices turns out to be not easy to digest, at least 
at a first reading. This is the reason that in contrast to Sen, we have started in
a different way, namely with the set $\Cal T$. As we saw, there is a canonical bijection 
between the set $\Cal T$ and the set $\Cal V$ of valleys $V,$ and it is the set $\Cal V$ which
has to be considered as a natural setting for dealing with relations. We recall that 
relations can be, of course, interpreted as
non-zero elements of $\Ext^2$, say of $\Ext^2(-,-)$. Any valley $V$ yields an
exact sequence
$$ 
 0 @>>> T @>>> IT @>f>> PU @>>> U @>>> 0
$$
with $V = \Im f$. Its equivalence class in $\Ext^2(U,T)$ is non-zero. In this way, the set of
valleys correspond to a minimal set of relations which defines $A$.
	\bigskip\medskip

{\bf Appendix D. Further invariants related to the finitistic dimension.}
	\medskip
In order to get hold of $\finpro A$, 
G\'elinas discusses in the paper [G1] not only $\del S$, but also the grade of
simple modules $S$ (thus the depth of $A$), and he looks   
for the numbers $t\ge 0$ such that $\mho^t\Omega^t S$ is torsionless. 
Examples of Nakayama algebras can be used in order to
get a better understanding of these settings. 
	\medskip
{\bf D.1. The grade of a simple module and the depth of the algebra.}
	\medskip
Let $A$ be any artin algebra. 
Let $M$ be a module. The {\it grade $\grade M$} is defined as follows.
We have $\grade M = 0$ iff $\Hom(M,A) \neq 0$, and	
$\grade M = d \ge 1$ iff $\Hom(M,A) = 0$, $\Ext^i(M,A) = 0$ for $1\le i < d$ and
$\Ext^d(M,A) \neq 0.$ 
The {\it depth} of $A$ is the maximum of $\grade S$, where $S$ are the simple modules.
	\medskip
{\bf Proposition (Jans--G\'elinas).} {\it Let $S$ be a simple module with finite
$\grade S = d\ge 1$. There are modules $N$ with 
injective dimension $d$ which satisfy $\Ext^d(S,N) \neq 0,$ for example $N =  \tau \Omega^{d-1}S$.}
	\medskip
Proof.
Let $P_\bullet$ be a minimal projective resolution of $S$ and
consider its truncation
$$
 P_d @>f_d>> P_{d-1} @>>> \cdots @>>> P_1 @>f_1>> P_0 @>>> 0.
$$
Here, $P_0 = P(S)$ is a projective cover of $S$. For all $0\le i < d,$
the cokernel of $f_{i+1}$ is  $\Omega^iS$. 
In particular, $X = \Omega^{d-1}S$ is the cokernel of $f_d.$ 
We form the $A$-dual
$$
 P^*_d @<f^*_d<< P^*_{d-1} @<<< \cdots @<<< P^*_1 @<f^*_1<< P^*_0 @<<< 0.
$$
This is an exact sequence, since $\Ext^i(S,A) = 0$ for $0\le i <d.$
By definition, $\Tr X$ is the cokernel of $f^*_d$. Thus, we have obtained 
a minimal projective resolution of the right module $\Tr X$
$$
 0 @<<< \Tr X @<<< P^*_d @<<< P^*_{d-1} @<<< \cdots @<<< P^*_1 @<<< P^*_0 @<<< 0. 
$$

We apply $D$ to this sequence  and get a minimal injective coresolution of $N =
\tau \Omega^{d-1}S = D\Tr \Omega^{d-1}S = D\Tr X$ 
with $DP^*_0 = DP(S)^* = I(S)$. This shows that $\id N = d$ and $\Ext^d(S,N) \neq 0$.
	\bigskip 
{\bf Corollary 1.} {\it For any simple module $S$, we have $\grade S \le \del S.$}
	\medskip
Proof. This is trivially true if $\grade S = 0.$ Thus we can assume that 
$\grade S = d \ge 1$. The proposition asserts that there is a module $N$
with injective dimension $d$ which satisfies $\Ext^d(S,N) \neq 0.$
According to Proposition 2.2, we get $d \le \del S.$
$\s$
	\medskip 
{\bf Corollary 2.} 
$$ 
  \depth A \le \del A.
$$
\vglue-.5cm
$\s$
	\bigskip
{\bf D.2. The case of Nakayama algebras.} 
	\medskip
{\bf Lemma.} {\it Let $A$ be Nakayama and let $S$ be simple. Then $\del S = 0$ iff
$\grade S = 0$ iff $S$ is torsionless. Next,
$\del S = 1$ iff
$\grade S = 1$ iff $S$ is not torsionless and $\pd S = 1.$
Also, $\del S = 2$ implies that $\grade S = 2$}
(the converse is not rue, as we will see at the end of this section). 

	\medskip
Proof. If $\Omega S$ is non-zero, then $\Ext^1(S,\Omega S) \neq 0$. 
Thus, if $\Omega S$ is non-zero and projective, then $\grade S \le 1.$
Conversely, let $\grade S = 1,$ and assume that $\Omega S$ is not
projective. Then $\Omega S \to P(S) \to S $ is not an $\mho$-sequence,
but this implies that $\rad P(S) = \Omega S \to P(S)$ is not an $\add A$-approximation, 
thus $P(S)$ is not a minimal projective. But this implies that $\rad P(S)$ is projective,
a contradiction. 

Next, let $\del S = 2$. Then $\grade S \ge 2$. But in general, $\grade S \le \del S$.
Thus $\grade S = 2.$
  $\s$
	\bigskip
{\bf Examples.} 
{\it There are cyclic Nakayama algebras of arbitrarily large finite global dimension with
depth equal to $2$.}
	\medskip 
Let $m\ge 2$ and $n = 2m.$ Let $A_n$ be the 
cyclic Nakayama algebra with Kupisch series $(2,3,3,4,4,5,5,\cdots, 
m+1,m+1,m+2).$ We denote the simple modules by $S_i$ with $1\le i \le n$, where $P(S_1) = 2$ and
$S_{i+1} = \tau^-S_i.$ Thus, $|P(S_{2i-1})| = i+1$ and $|P(S_{2i})| = i+2$ for all $i\ge 1.$ In
particular, we see that $\pd S_{2i} = 1$ for all $i\ge 1.$ 
Let us look at the simple modules with odd index.
Since $\Omega S_1 = S_n = S_{2m}$ has projective dimension $1$, we have $\pd S_1 = 2.$ 
Since $\Omega^2 S_3 = S_n = S_{2m}$ has projective dimension $1$, we have $\pd S_3 = 3.$ 
For $i\ge 2,$ we have $\Omega^2 S_{2i+1} = S_{i-1}.$ This shows that $A_n$ has finite global dimension.
Also, 
$$
  \pd S_{2^t-3} = 2t-2,  \quad\text{for}\quad t\ge 2.
$$ 
(The proof is by induction on $t \ge 2$. As we know, $\pd S_1  = 2$. We have $\Omega S_{2^{t+1}-3}
= P(S_{2^{t+1}-4})/\soc$ and $\soc P(S_{2^{t+1}-4}) = S_{2^t-3}$. By induction, we assume that 
$\pd S_{2^t-3} = 2t-2,$ therefore $\pd S_{2^{t+1} -3} = 2+\pd S_{2^t-3} = 2t = 2(t+1)-2,$ as we want to show.)
Thus the global dimension of $A_n$ gets arbitrarily large.

It remains to be seen that $\depth A_n = 2.$ Since the simple modules with even index have projective
dimension $1$, we have $\grade S_{2i} \le 1$ for all $i$. We have 
$\Omega S_{2i+1} = \rad P(S_{2i+1}) = P(S_{2i})/\soc$
and there is a non-split exact sequence
$$ 
 0\to   P(S_{2i-1}) \to P(S_{2i}) \oplus P(S_{2i-1})/\soc \to  P(S_{2i})/\soc \to 0
$$
(an Auslander-Reiten sequence), thus 
$$ 
 \Ext^2(S_{2i+1},P(S_{2i-1})) = \Ext^1(P(S_{2i})/\soc,P(S_{2i-1})) \neq 0.
$$
This shows that $\grade S_{2i+1} \le 2$ (and actually, $\grade S_{n-1} = 2$). 
$\s$

Here is the case $n = 8$. On the left, there is the Auslander-Reiten quiver of $A_8$ with 
the function $\pd$,  
and below the indices of the simple modules. On the right we  highlight the module $S_7$ and 
the Auslander-Reiten
sequence which shows that $\Ext^2(S_7,P(S_5))\neq 0.$
$$ 
{\beginpicture
    \setcoordinatesystem units <.4cm,.4cm>
\put{\beginpicture
\multiput{} at  0 0  16 5 /

\setdots <1mm>
\plot 0 0  16 0 /

\setdots <0.5 mm>
\plot 0 0  2 2  4 0  8 4  12 0  14 2  /
\plot 1 1  2 0  5 3  8 0  12 4 /
\plot 4 2  6 0  11 5  16 0 /
\plot 7 3  10 0 13 3 /
\plot 10 4  14 0  15 1 / 
\setdashes <1mm>
\plot 0 0  0 2 /
\plot 16 0  16 2 /
\multiput{$\ss 0$} at 1 1  2 2  4 2  5 3   7 3  8 4  10 4  11 5 /
\multiput{$\ss 1$} at 0 0  4 0  8 0  12 0  16 0 /
\multiput{$\ss 2$} at 2 0  3 1  8 2  9 3   /
\multiput{$\ss 3$} at 5 1  6 0  6 2  7 1  13 1  14 0  14 2  15 1 /
\multiput{$\ss 4$} at 9 1  10 2  11 3  12 4    10 0  11 1  12 2  13 3  /

\multiput{$\ss 4$} at 8 -1 /
\multiput{$\ss 3$} at 6 -1  /
\multiput{$\ss 2$} at 4 -1  /
\multiput{$\ss 1$} at 2 -1 /
\multiput{$\ss 5$} at 10 -1 /
\multiput{$\ss 6$} at 12 -1 /
\multiput{$\ss 7$} at 14 -1 /
\multiput{$\ss 8$} at 16 -1 /
\endpicture} at 0 0
\put{\beginpicture
\multiput{} at  0 0  16 5.6 /

\setdots <1mm>
\plot 0 0  16 0 /

\setdots <0.5 mm>
\plot 0 0  2 2  4 0  8 4  12 0  14 2  /
\plot 1 1  2 0  5 3  8 0  12 4 /
\plot 4 2  6 0  11 5  16 0 /
\plot 7 3  10 0 13 3 /
\plot 10 4  14 0  15 1 / 
\setdashes <1mm>
\plot 0 0  0 2 /
\plot 16 0  16 2 /

\put{$S_7$} at 14 0 
\put{$\ss\Omega S_7$} at 9.5 3 
\put{$\ss P(S_5)$} at  6.5 3 

\multiput{$\ss 4$} at 8 -1 /
\multiput{$\ss 3$} at 6 -1  /
\multiput{$\ss 2$} at 4 -1  /
\multiput{$\ss 1$} at 2 -1 /
\multiput{$\ss 5$} at 10 -1 /
\multiput{$\ss 6$} at 12 -1 /
\multiput{$\ss 7$} at 14 -1 /
\multiput{$\ss 8$} at 16 -1 /
\setshadegrid span <.5mm>
\vshade 7 3 3  <z,z,,> 8 2 4  <z,z,,> 9 3 3 /

\endpicture} at 18 0
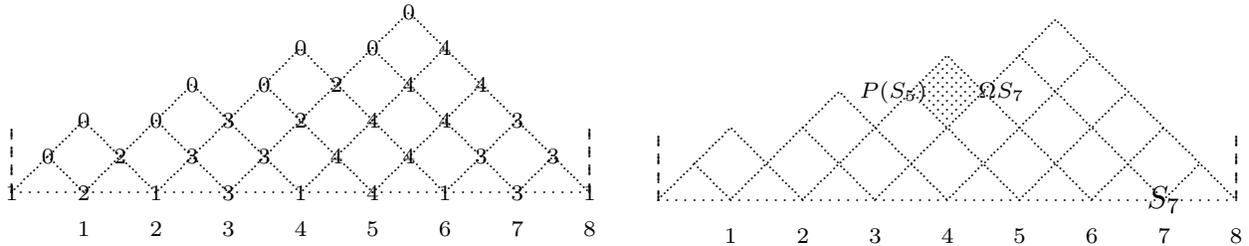
\endpicture}
$$
	\medskip
Already $A_4$ has a simple module $S$ with $\grade S = 2$ and $\del S = 3.$
Here is the algebra with the module $S = S_3$ encircled. 
$$ 
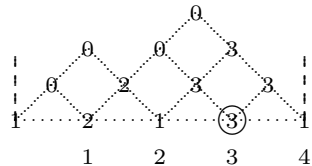
{\beginpicture
    \setcoordinatesystem units <.4cm,.4cm>
\multiput{} at  0 0  8 3  /

\setdots <1mm>
\plot 0 0  8 0 /

\setdots <0.5 mm>
\plot 0 0  2 2  4 0  6 2   /
\plot 1 1  2 0  5 3  8 0   /
\plot 4 2  6 0  7 1 /
\setdashes <1mm>
\plot 0 0  0 2 /
\plot 8 0  8 2 /
\multiput{$\ss 0$} at 1 1  2 2  4 2  5 3  /
\multiput{$\ss 1$} at 0 0  4 0  8 0 /
\multiput{$\ss 2$} at 2 0  3 1  /
\multiput{$\ss 3$} at 5 1  6 0  6 2  7 1 /

\multiput{$\ss 4$} at 8 -1 /
\multiput{$\ss 3$} at 6 -1  /
\multiput{$\ss 2$} at 4 -1  /
\multiput{$\ss 1$} at 2 -1 /
\put{$\bigcirc$} at 6 0 
\endpicture}
$$
We have mentioned already that $\grade S_3 = 2$. Its projective dimension is $3$, thus
$\Omega^2 S_3$ is not projective. On the other hand, since the global dimension of $A_4$
is 3, all modules of the form $\Omega^3 M$ are projective. Thus, we see that 
$\del S_3 \ge 3$ (and therefore $\del S_3 = 3$). 
	\bigskip
\vfill\eject
{\bf D.3. The modules $\mho^t\Omega^t S$.}
	\medskip
If $M$ is a module, let $\mho M = \Tr\Omega\Tr M$ (note that $\mho$ is a stable functor; it was
considered by Auslander and Reiten in [AR], and is left adjoint to $\Omega$;
one can show that $\mho M$ is just the cokernel of
a minimal left $\add{}_AA$-approximation of $M$, see [RZ], 4.4).

If $S$ is a simple module, G\'elinas looks for a non-negative number $t$
such that $\mho^t\Omega^tS$ is torsionless. Namely,  
Theorem 1.13 of G\'elinas [Ge1] asserts:
{\it If \ $\mho^t\Omega^t M$ is torsionless, then $\del M \le t$.}
There is the following immediate consequence
(the ``torsionfreeness criterion'' 1.14): {\it If there is $t\ge 0$ such that 
$\mho^{t}\Omega^{t} S$ is torsionless, then $\del A \le
\sup n_S.$}
	\medskip
Here is an example of a cyclic Nakayama algebra and 
a simple module $S$ such that no module $\mho^t\Omega^t S$ with $t\ge 0$ is torsionless.
$$ 
{\beginpicture
    \setcoordinatesystem units <.4cm,.4cm>
\put{\beginpicture
\multiput{} at 0 -1.8   10 3 /
\setdots <0.5 mm>
\plot 0 0  2 2  4 0  6 2 /
\plot 8 2  10 0 /
\plot 1 1  2 0  5 3  8 0  9 1 /
\plot 4 2  6 0  8 2 /
\plot 2 2  3 3  4 2 /
\setdots <1mm>
\plot 0 0  10 0 /
\setdashes <1mm> 
\plot 0 2  0 -.2 /
\plot 10 2  10 -.2 /
\put{$S$} at 8 -.8
\multiput{$\ss 0$} at 1 1  2 2  3 3  5 3  8 2  /
\multiput{$\ss 1$} at 4 0  5 1  6 0   /
\multiput{$\ss 2$} at 9 1 /
\multiput{$\ss \infty$} at 0 0  2 0  3 1  4 2  6 2  7 1  8 0  10 0  /
\endpicture} at 0 0 
\put{\beginpicture
\multiput{} at 0 0   10 3 /
\setdots <0.5 mm>
\plot 0 0  2 2  4 0  6 2 /
\plot 8 2  10 0 /
\plot 1 1  2 0  5 3  8 0  9 1 /
\plot 4 2  6 0  8 2 /
\plot 2 2  3 3  4 2 /
\setdots <1mm>
\plot 0 0  10 0 /
\setdashes <1mm> 
\plot 0 2  0 -.2 /
\plot 10 2  10 -.2 /
\multiput{$\bullet$} at 8 0  4 2  0 0    3 1  7 1  10 0 /
\put{$S$} at 8 -.8
\put{$\ss \Omega S$} at 4 2.8 
\put{$\ss \Omega T$} at 7 -1.5
\put{$\ss T$\strut } at 10 -.9
\arr{7 -1}{7 .5} 
\put{$\ss \Omega^2 T$} at 3 -1.5
\arr{3 -1}{3 .5} 
\put{$\ss T = \Omega^2 S$\strut} at 0 -.9
\multiput{$\circ$} at 2 0 /
\put{$\ss \mho T$\strut} at 2 -.6 
\put{$\ss X$} at 6.1 2.6
\put{$\circ$} at 8 2
\put{$\ss PT$} at 8.1 2.6
\multiput{$\star$} at 6 2 /
\multiput{$\ss 0$} at   /
\multiput{$\ss 1$} at  /
\multiput{$\ss 2$} at  /
\multiput{$\ss \infty$} at  /
\endpicture} at 14 0 
\endpicture}
$$
{\bf Claim.} 
{\it We have $\mho \Omega S = S$ and $\mho^t\Omega^tS = X$ for $t\ge 2$} (and both modules
$S$ and $X$ are not torsionless). On the right, we have marked by bullets
the modules $\Omega^t S$ with $0 \le t \le 4$ 
and stress that $\Omega^5 S = \Omega^2 S$. We write $T = \Omega^2S$. 
Also, we have inserted $X = \mho^2 T$ as well as the projective cover $PT$ of $T$. 

Proof. One checks immediately that $\mho\Omega S = S$. Thus, let us show that 
$\mho^t\Omega^tS = X$ for $t\ge 2.$ 
It is easy to see that we have: (a) $\Omega^3 T = T$,  
(b) $\mho\Omega T = T,$ and (c) $\mho^2\Omega^2 T = T$. 
We define (d) $X = \mho^2 T.$ Then we have (e)
$\mho X = T$ 
(we stress that a minimal $\add A$-approximation $X \to PT$ is not a monomorphism, 
but this does not matter!).
It follows from (d) and (e) that we have (f) $\mho^3 T = T.$

Now, let $t\ge 2.$ If 
$t = 3s$ with $s\ge 1,$ then 
$$
 \mho^t\Omega^t S = \mho^2\mho^{3(s-1)}(\mho\Omega)\Omega^{3(s-1)}T = X
$$ 
using (a), (b), (f) and (d). Next, for $t = 3s+1$
with $s\ge 1,$ we get
$$
 \mho^t\Omega^t S = \mho^2\mho^{3(s-1)}(\mho^2\Omega^2)\Omega^{3(s-1)}T = X
$$ 
using (a), (c), (f) and (d). Finally, for $t = 3s+2$
with $s\ge 0,$ we get
$$
 \mho^t\Omega^t S = \mho^2\mho^{3s}\Omega^{3s}T = X
$$ 
using (a), (f) and (d). $\s$	
	\bigskip\bigskip
	
{\bf References.}
	\medskip
\item{[AR1]} M. Auslander, I. Reiten.
    k-Gorenstein algebras and syzygy modules. J. Pure Appl. Algebra 92 (1994) 1--27.
\item{[AR2]}  M. Auslander, I. Reiten.
   Syzygy modules for noetherian rings. J. Algebra 183 (1996), 167--185.
\item{[B]} H. Bass. Finitistic dimension and a homological generalization of 
   semi-primary rings. Trans. Amer. Math. Soc. 95 (1960), 466--488 
\item{[Ga]} P. Gabriel. The universal cover of a representation-finite algebra. In: Representations
  of Algebras. Springer Lecture Notes in Mathematics 903 (1981), 68--105. 
\item{[Ge]} V. G\'elinas. The depth, the delooping level and the finitistic dimension.
   \newline
   arXiv:2004.04828 (2020).
\item{[Gu]} W. H. Gustafson. Global dimension in serial rings. Journal of Algebra, 
 97 (1985), 14--16.
\item{[M]} D. Madsen. Projective dimensions and Nakayama
   algebras. Fields Institute Communications. 45.
   Amer\. Math\. Soc\., Providence, RI, 2005. 247--265.
\item{[N]} T. Nakayama. On Frobeniusean algebras. I.
    Annals of Math. 40 (1939), 611--633.
\item{[R1]} C. M. Ringel. Representations of $k$-species and bimodules.
  J. Algebra 41 (1976), 269--302.
\item{[R2]} C. M. Ringel.  The Gorenstein projective modules for the Nakayama algebras. I.
    Journal of Algebra (2013), 241--261. DOI 10.1016/j.jalgebra.2013.03.014 
\item{[RS]} M. Rubey, C. Stump, et al. 
  The Ringel permutation of the LNakayama algebra corresponding to a Dyck path.
   Mp00201. 
  In: FindStat --- The combinatorial statistics database (2020),
  www.findstat.org/Mp00201. 
\item{[Sh1]}  D. Shen. A note on resolution quivers. 
   J. Algebra Appl. 13 (2014), 1350120.
\item{[Sh2]}  D. Shen. The singularity category of a Nakayama algebra. J. Algebra
    429 (2015), 1--18.
\item{[Sh3]}  D. Shen. A note on homological properties of Nakayama algebras.
    Archiv der Mathematik. 108 (2017), 251--261.
\item{[S1]} E. Sen. The $\phi$-dimension of cyclic Nakayama algebras. arXiv 1806.01449 (2018). 
\item{[S2]} E. Sen. Syzygy filtrations of cyclic Nakayama algebras. 
   arXiv 1903.04645 (2019).
\item{[S3]} E. Sen. Delooping level of Nakayama algebras. arXiv:2009.08888 (2020).
\item{[S]} J. P. Serre. Sur  la  dimension homologique des  anneaux et  des  modules Noetheriens.
    Proceedings International Symposium on  Algebraic Number Theory, Tokyo and  Nikko (1955),
    175--189.
\item{[RZ]} C. M. Ringel, P. Zhang.  Gorenstein-projective and semi-Gorenstein-projective modules.
    Algebra \& Number Theory. Vol. 14 (2020), No. 1, 1--36.
    \newline DOI: 10.2140/ant.2020.14.1 
	\bigskip\bigskip
{\baselineskip=1pt
\rmk
C. M. Ringel\par
Fakult\"at f\"ur Mathematik, Universit\"at Bielefeld \par
POBox 100131, D-33501 Bielefeld, Germany  \par
ringel\@math.uni-bielefeld.de
\par}

\bye